
  \magnification=\magstep0
  \baselineskip=16pt
  \overfullrule=0pt

  %
%

%
%
\font\twelverm=cmr12
\font\twelvebf=cmbx12
\font\twelveit=cmti12
\font\twelvesl=cmsl12
\font\eightrm=cmr8
\font\eightbf=cmbx8
\font\eightit=cmti8
\font\eightsl=cmsl8

\def\bigtype{\let\rm=\twelverm
             \let\bf=\twelvebf
             \let\it=\twelveit
             \let\sl=\twelvesl}
\def\medtype{\let\rm=\tenrm
             \let\bf=\tenbf
             \let\it=\tenit
             \let\sl=\tensl}
\def\smalltype{\let\rm=\eightrm
             \let\bf=\eightbf
             \let\it=\eightit
             \let\sl=\eightsl}

%
%
\def\chapSkip{\vfill\eject} 
\def\Chapter#1{\rightline{\bigtype\bf #1}}
\def\secSkip{\vskip .3in}
\def\Section#1{\par\noindent{\bf #1}\hskip.1in}
\def\parSkip{\vskip .1in}

\def\Number#1{{\rm #1}\hskip .1in}
\def\Title#1{{\bf #1}\hskip .1in\nobreak}
\def\Subtitle#1{{\rm ( #1 )}\hskip .1in\nobreak}

\def\Proof{\par\parSkip{\sl Proof}\hskip .2in\nobreak}
\def\Qed{\nobreak\hbox{ $\bullet$}}
\def\qed{\nobreak\hbox{ $\circ$}}

%
%
\def\dispM#1{{$$ #1 $$}} 

\def\dispm#1{$ #1 $} 

%
%
\def\Fun#1#2{{}^{#1}{#2}}

%
%
\def\and{{\,\&\,}}

\def\maxs#1{\max\{ #1, \}}         
\def\mod#1{\vert #1 \vert}         
\def\sz#1{\vert #1 \vert}          
\def\cns#1{\vert\{ #1 \}\vert}     
\def\tup#1{\langle #1 \rangle}     
\def\gtup#1{{ #1 }}                
\def\cat{\wedge}                   
\def\wcat{\;\cat\;}                
\def\tac#1{\cat\tup{#1}}           
\def\ens#1{\{ #1 \}}               
\def\rk#1{{\rm rk}(#1)}            
\def\ran{{\rm ran}}                
\def\wran{{\ran\;}}                
\def\dom{{\rm dom}}                
\def\supp{{\rm supp}}              
\def\ric{\lceil}                   
\def\nric#1{{\lceil}_{#1}}         
\def\isom{\cong}                   
\def\aeq{=^{\ast}}                 
\def\const{\equiv}                 
\def\rid{\setminus}                
\def\rids#1{\rid\{ #1 \}}            
\def\diff{\,\triangle\,}           
\def\sub{\subset}                  
\def\subeq{\subseteq}              
\def\supeq{\supseteq}              
\def\larger{\succ}                 
\def\then{\rightarrow}             
\def\halts{\downarrow}
\def\loops{\uparrow}

\def\sgn{{\rm sgn}}
\def\fin{{\rm fin}}

%
%
\def\ONE{{\bf 1}} 
\def\N{{\rm N}} 
\def\Z{{\rm Z}} 
\def\H{{\N}^{\ast}} 

%
%
\def\ss{{\sigma}}                  
\def\ps{{\pi}}                     
\def\tp{{\tau}}                    
\def\osc{\delta}                   


%
%
\def\co#1#2{{[#1]}_{#2}}           
\def\cmp#1#2#3{[#1]_{#2}^{#3}}     
\def\arr#1{[ #1 ]}                 
\def\rarr#1{\{ #1 \}}              

%
%
\def\A{{\sl A}} 

%
%
\def\Hyp#1#2#3{{\sl H}^{#1}_{#2}(#3)} 
\def\Suc{{\sl S}} 
\def\Tl{{\sl T}} 

\def\M{{M}} 
\def\U{{U}} 
\def\V{{V}} 
\def\W{{W}} 

\def\MM{{{\bf M}}} 
\def\UU{{{\bf U}}} 
\def\VV{{{\bf V}}} 
\def\WW{{{\bf W}}} 




\def\UC{{\cal U}} 
\def\UT{{\UC}_0} 

%
%
\def\f{{\rm f}} 
\def\j{{\rm j}} 

\def\SS{{\cal S}} 
\def\CSE{{\cal G}} 
\def\CSI{{\cal G}} 
\def\CSP{{\cal G}^o} 
\def\Dom{{\cal D}} 

\def\Tree{{\cal T}} 
\def\Lev{{\rm Lev}} 
\def\Ext{{\rm Ext}} 
\def\Path{{P}} 

\def\nullseq{\tup{\;}} 
\def\root{\nullseq} 
\def\void{\emptyset} 

\def\lf{{\ell}} 
\def\lfs{{\ell}} 
\def\leaf{{\lambda}} 
\def\Lf{{\cal L}} 
\def\Lfs{\Lf} 


\def\wt{{\psi}} 
\def\mwt{\xi} 
\def\aug#1{{\widehat{#1}}} 
\def\xwt{\aug{\wt}} 
\def\rn{{r}} 
\def\rwt{{w}} 
\def\rxwt{\aug{\rwt}} 
\def\kwt{\zeta} 

\def\wts{\wt} 
\def\mwts{\mwt} 
\def\xwts{\xwt} 
\def\rwts{\rwt} 
\def\rxwts{\rxwt} 
\def\kwts{\kwt} 


\def\rmwts{{X}} 

\def\equipoll{\approx}
\def\smalleq{\preceq}

\def\assoc{\to} 

\def\free{\phi} 
\def\frees{\phi} 

\def\new{{{\nu}}} 
\def\news{{{\nu}}} 

\def\nnw{\eta} 
\def\kkw{\kappa} 

\def\rkkw{{K}} 
\def\rllw{{L}} 

\def\rnnw{{H}} 

\def\tM{\widetilde{\M}}
\def\tMM{\widetilde{\MM}}

  \def\dsc{\delta}

%
%

%
%
\def\shortVert{{{}_{{}_{\vert}}}}
\def\raiseVert{\raise 4pt\hbox{$\vert$}}
\def\obranch{\buildrel\shortVert
                      \over
                      {\overline{\shortVert\qquad\qquad\qquad\qquad
                                 \qquad\qquad\qquad\qquad\shortVert}}}
\def\pbranch{\buildrel\shortVert
                      \over
                      {\overline{\shortVert\qquad\qquad
                                 \qquad\qquad\shortVert}}}
\def\qbranch{\buildrel\shortVert
                      \over
                      {\overline{\shortVert\quad\quad\quad\shortVert}}}
\def\mrx#1#2#3{\matrix{
               \hfill\hrulefill
               &#1&
               \hrulefill\hfill\cr
               \vert&&\vert\cr
               {#2}&&{#3}\cr}}
\def\nrx#1#2{\matrix{#1\cr\vert\cr#2\cr}}
\def\orx#1#2#3{\matrix{
               #1\cr
               \obranch\cr
               #2\qquad#3\cr}}
\def\prx#1#2#3{\matrix{
               #1\cr
               \pbranch\cr
               #2\qquad#3\cr}}
\def\qrx#1#2#3{\matrix{
               #1\cr
               \qbranch\cr
               #2\quad#3\cr}}
\def\trx#1{\matrix{#1\cr\cr\cr}}


  \pageno=0
  \vskip .8in
  {\bigtype
  \centerline{{\bigtype\bf Uniform chains of small transfinite rank}}
  }
  \vskip .8in
  \centerline{\smalltype F.P.Weber}
  \centerline{\smalltype R.M.Vollmer}
  \vskip .2in
  \centerline{\smalltype School of Earth Sciences}
  \centerline{\smalltype University of Leeds}
  \centerline{\smalltype LS2 9JT}
  \centerline{\smalltype United Kingdom}

  \vskip .8in
  \centerline{\bf Abstract}

  \parSkip\noindent
  We extend F.Pastjin's construction of uniform decomposable 
  chains of finite rank to those of rank 'finite over a limit' 
  and investigate infinite (length $\omega$) products and unions 
  of such chains. We derive an extension of Pastjin's 
  characterisation to uniform decomposable chains of small 
  transfinite rank ($\le\omega+\omega$). We conclude by indicating 
  how the resulting descriptions can be iterated to chains of rank 
  $\le\omega^{\omega}$.

  \vskip .8in
  \centerline{\bf Contents}
  \item{} [1] Pastjin's Chains
  \item{} [2] Tails of first limit rank
  \item{} [3] Chains of rank finite over a limit
  \item{} [4] Products of uniform chains
  \item{} [5] Chains of first limit rank
  \item{} [6] Chains of next limit rank
  \item{} [7] Short representations
  \item{} [8] Short representations of tails
  \item{} [9] Extensions

  \pageno=0
  
\chapSkip
\Chapter{[1] Pastjin's Chains}

\secSkip
\Section{[1.1]}
\Title{Decomposable chains}

\parSkip\noindent
Let $A$ be a chain (i.e., a linear, or total, order). 
Define an equivalence relation $\rho = \rho^A$ by:
    $x \rho y$ 
if the interval between $x$ and $y$ in $A$ 
is finite.  Let 
    $\cmp{x}{}{} = \cmp{x}{}{A} = \ens{y\in A: x \rho y}$
(the component of $x$ in $A$) and
    $A/\rho = \ens{\cmp{x}{}{} : x\in A}$.

Iterate the definition through all ordinals:
Set
    $\rho_1 = \rho$;
suppose $\rho_{\alpha}$ defined, set
    $A/\rho_{\alpha}=\ens{\cmp{x}{\alpha}{}:x\in A}$
and
    $x \rho_{\alpha+1} y$ in $A$
if 
    $\cmp{x}{\alpha}{} \rho \cmp{y}{\alpha}{}$ in $A/\rho_{\alpha}$
(where
    $\cmp{x}{\alpha}{} = \cmp{x}{\alpha}{A} 
    = \ens{y\in A: x \rho_{\alpha} y}$
is the {\it $\alpha$-component} of $x$ (in $A$));
suppose $\rho_{\beta}$ defined for all $\beta<\alpha$, for
limit $\alpha$, set
    $A/\rho_{\beta}=\ens{\cmp{x}{\beta}{}:x\in A}$
and set
    $x \rho_{\alpha} y$ in $A$
if for some $\beta<\alpha$ it is the case that
    $\cmp{x}{\beta}{} \rho \cmp{y}{\beta}{}$ in $A/\rho_{\beta}$.
In the following, write 
    $A/\alpha$ for $A/\rho_{\alpha}$.
and 
    $\co{x}{\alpha}$ for $\cmp{x}{\alpha}{}$.
The {\it rank} of $A$ is the least ordinal $\alpha$ such that
$A/{\alpha} \isom 1$ (the singleton chain) and is denoted $\rk{A}$.
A chain is {\it decomposable} if $\rk{A}$ is defined.
Subsequently, all chains are assumed decomposable.

\parSkip
\Number{(1.2)}
\Title{Definition}
\Subtitle{Uniform chains}
Let $A^{\ge x}\isom\ens{y\in A: x\le y}$.
A chain $A$ is (right) {\it uniform} if for all $x,y\in A$ the 
chains $A^{\ge x}$ and $A^{\ge y}$ are order isomorphic. 

In this case denote by $\f A$ (a representative sample of) the 
tail of $A$. (We shall also use $\f_xA$ to denote $A^{\ge x}$ 
and use $\f$ and $\f_x$ interchangably with uniform chains.) 
$\f$ is an idempotent operation, $\f^2=\f$. Let $\UC(\alpha)$
denote the set of uniform chains of rank $\alpha$ and 
$\UT(\alpha)$ the set of tails of chains in $\UC(\alpha)$.
$\UC(<\alpha)=\bigcup_{i<\alpha}\UC(i)$.
If not specified otherwise, letters $A$ and $B$ refer in the
following to uniform (decomposable) chains.

\secSkip
\Section{[1.2]}
\Title{Pastjin's uniform chain generator}

\parSkip\noindent
Let $\N$ be the chain of non-negative integers and $\Z$ the
chain of all integers. There is a scheme to generate all
uniform chains of finite rank as chains of form $A\Z$ or
$B+A\N$ from uniform chains of smaller finite rank, starting
with the singleton chain. ($AB$ is the anti-lexicographically
ordered cartesian product of chains $A$ and $B$, effectively
    $AB = \sum_{b\in B}A$.) 

\parSkip
\Number{(1.3)}
\Title{Definition}
\Subtitle{Sub-diagonal strings}
For an ordinal $\alpha$, let
    $\SS(\alpha)=\ens{ s\in\Fun{\alpha}{\alpha} 
    : s(i)\le i \hbox{ for all } i<\alpha }$.
Also let
    $\SS(<\alpha) = \bigcup_{i<\alpha} \SS(i)$
and, e.g., (for ordinal $\beta$)
    $\SS([\beta,\alpha)) = \bigcup_{\beta\le i<\alpha} \SS(i)$,
etc. (Note that $\SS(1) = \ens{\tup{0}}$.)

\parSkip
\Number{(1.4)}
\Title{Definition}
\Subtitle{Pastjin's finite rank unifom chain constructor}
Set $\A(\tup{0}) = 1 = \tup{0}$,
and for any $s\in\SS$ 
\dispM{
    \A(s\tac{a,b}) 
    = \cases{
    \A(s\tac{a})\Z            
    & if $b = \mod{s\tac{a}}$ \cr
    \A(s\tac{b}) + \A(s\tac{a})\N 
    & else, i.e., if $b < \mod{s\tac{a}}$ \cr
    }
}

(We use '$\tac{}$' to denote concatenation of finite strings,
while retailing the conventional '$+$' for chains.) Letting
    $\A(\nullseq) = \nullseq = \bf{ 0 }$ ($:=$ the empty chain),
by finite induction,
\dispM{
    \A(s) 
    = \A(s\ric k(s))\Z
    + \sum_{k(s)<i<n(s)}\A(s\ric i)\N
}
where 
    $n(s)=\mod{s}-1 $ 
and 
    $k(s)=\max\ens{i\le n:s(n(s))=i}$.

Here is a basic list:
$$\vbox{
\settabs 3\columns
\+$\A(\tup{0}) \isom 1$&& \cr
\+$\A(\tup{0,0}) \isom 1 + 1\N \isom \N$&
  $\A(\tup{0,1}) \isom 1\Z \isom \Z$& \cr
\+$\A(\tup{0,0,0}) \isom 1 + 1\N + \N\N \isom \N^2$&
  $\A(\tup{0,0,1}) \isom \Z + \N\N \isom \Z + \N^2$&
  $\A(\tup{0,0,2}) \isom N\Z$ \cr
\+$\A(\tup{0,1,0}) \isom \N + \Z\N$&
  $\A(\tup{0,1,1}) \isom \Z\N$&
  $\A(\tup{0,1,2}) \isom \Z^2$ \cr
}$$

\parSkip
\Number{(1.5)}
\Title{Theorem}
\Subtitle{Pastjin}
The chains $\A(s),s\in\SS(n+1)$ exhaust all uniform chains
of rank $n$ (up to isomorphism).


\parSkip
\Number{(1.6)}
\Title{Notation}
Define operation $s\to s^+$ and $s\to s^-$ on binary strings
$s$ by:
\dispM{
    (s\tac{a})^-   = s, \qquad
    (s\tac{a,b})^+ = s\tac{b}.
}
Thus, where $n=\mod{s}-1$,
\dispM{
    \A(s) = \cases{
    \A(s^-)\Z           & if $s(n) = n$ \cr
    \A(s^+) + \A(s^-)\N & if $s(n) < n$.\cr
    }\eqno(1)
}
Define iterations of $()^+$ and $()^-$ by
    $s^{+(k+1)} = (s^{+k})^+$
and
    $s^{-(k+1)} = (s^{-k})^-$.
Note that $s^{+k}\in\SS$ just if $s(n)\le n-k$
(where $n=\mod{s}-1$), since $n-k=\mod{s^{+k}}-1$,
and $s^{+k}(n-k) = s(n)$. In particular, $s^{+n}\in\SS$ 
just if $s(n) = 0$, where $n=\mod{s}-1$.
For convenience, let 
    $s^{=} = s^{--} = (s^-)^-$,
and
    $s^{\equiv} = s^{---} = (s^{=})^-$.
(However,
    $s^{\pm}$ 
continues to read ''$s^-$ or $s^+$''.)
Note that
   $s^{+-}=s^{--}=s^=$.

\secSkip
\Section{[1.3]}
\Title{Uniform tails}

\parSkip
\Number{(1.7)}
\Title{Lemma}
\Subtitle{Finite tail lemma}
    $\f\A(s) 
    \isom \A(s^-\tac{0})
    = \sum_{i<\mod{s}}\A(s\ric i)\N$.

\Proof
The second equation is immediate. For the first, use induction. 
The base cases are trivial, e.g., 
\dispM{
    \f\A(\tup{0,1}) = \f\Z \isom \N = \A(\tup{0,0})
    = \A(\tup{0}\tac{0})
}
In the induction step use the following observation:
\dispM{
    \f_{(a,b)}AB = \f_aA + A\f_bB. \eqno(2)
}
If $s(n)=n$ for $n=\mod{s}-1$,
\dispM{
    \f\A(s) 
    = \f\A(s^-) + \A(s^-)\f\Z
    \isom \f\A(s^-) + \A(s^-)\N
    = \A(s^=\tac{0}) + \A(s^-)\N
    = \A(s^-\tac{0}),
}
from the induction hypothesis and the fact that
    $(s^-\tac{0})^+ = s^=\tac{0}$.

If $s(n)<n$ distinguish the case where the cut-off is in
$\A(s^+)$ or in $\A(s^-)\N$. In the respective case we 
get tails of form
\dispM{
    \f\A(s) = \cases{
    \f\A(s^+) + \A(s^-)\N
    = \A(s^=\tac{0}) + \A(s^-)\N, 
    & (cut-off in $\A(s^+)$) \cr
    \f\A(s^-) + \A(s^-)\f\N
    = \A(s^=\tac{0}) + \A(s^-)\N, 
    & (cut-off in some $\A(s^-)\times\ens{n}$), \cr
    }
}
again from the induction hypothesis and the fact that
    $(s^{\pm})^-\tac{0} = s^=\tac{0})$.
As before, either case reduces to
    $\A(s^-\tac{0})$.
\Qed

\parSkip
This suggests extending $\f$ to an operation on $\SS(\alpha)$
for successor ordinals $\alpha$ by
    $\f s = s^-\tac{0}$.
Thus,
    $\f\A(s) = \A(\f s)$
for $s\in \SS(n+1)$, $n\in\omega$. 
As a direct corollary:

\parSkip
\Number{(1.8)}
\Title{Lemma}
\Subtitle{String lemma}
\item{(i)}
    $(\f s)^+ = (s^-\tac{0})^+ = s^=\tac{0} 
    = \f s^- = \f s^+$;
\item{(ii)}
    $(\f s)^- = (s^-\tac{0})^- = s^-$.
\Qed

\parSkip
\Number{(1.9)}
\Title{Remark}
How do we identify elements in chains $\A(s)$? For instance,
in the second case of the decomposition of $A(s)$ in (1),
we could represent elements of $A(s)$ in the form:
\dispM{
    x=\cases{\tup{a,0}   &for $a$ in $\A(s^+)$ \cr
             \tup{b,n+1} &for $b$ in $\A(s^-)$, \cr 
    }
}
and we frequently do. Similarly, we use '$B\times\ens{n}$' to
refer to the $n+1$st hyper-component of $A+B\N$, where
$\rk{A}=\rk{B}$. On the other hand, since, e.g., a hyper-component
of $\A(s)$ is isomorphic to $\A(t)$ for $t=s^{\pm}$, we say
e.g., of $x\in\A(s)=\A(s^+)+\A(s^-)\N$ that it is 'in $\A(s^+)$'
to indicate that it is in the initial $\mod{s}-1$ component of
$\A(s)$, or that it is 'in some $\A(s^-)$' to indicate that
it is in one of the successor $\mod{s}-1$ components of $\A(s)$, 
and that it is 'in $\A(s^-)\times\ens{n}$' to indicate that it 
is in the $n$th successor $\mod{s}$ component of $\A(s)$. I.e.,
when considering embeddings between chains, we assume them to
be 'canonical', as defined below.

\parSkip
\Number{(1.10)}
\Title{Lemma}
\Subtitle{Sub-component lemma}
For $x\in A\in\UC$, each $\cmp{x}{i}{A}$, for $i\le\rk{A}$
is of the form $\A(s_i)$, for some $s_i\in\SS(i+1)$
such that 
\item{(a)}
$s_{i}\in\ens{s_{i+1}^+,s_{i+1}^-}$, hence
\item{(b)}
$s_{i}^-=s_{i+1}^=$, and 
\item{(c)}
$s_i^-\sub s$ for $s\in\SS(<\omega)$ such that $A=\A(s)$ 

\Proof
(a) and (b) are clear from the corollary to Pastjin's definition
and the preceding remark.
Setting $n=\rk{A}=\rk{\A(s)}=\mod{s}-1\in\omega$, we get 
    $s_n=s$, 
and for all $0<i\le n$,
    $s_{i-1}\in\ens{s_i^+,s_i^-}$,
by (a), so
    $s_{i-1}^-\in\ens{s_i^{+-},s_i^{--}}
    =\ens{s_i^=}$,
and so by finite ('downward') induction:
    $s_{i-1}^-\sub s_i^-\sub s_n^-=s^-$.
\Qed

\parSkip
\Number{(1.11)}
\Title{Definition}
\Subtitle{Canonical embeddings}
Let $s\in\SS(n+2),t\in\SS(n+1)$ such that $t\in\ens{s^+,s^-}$.
Call an embedding of $\A(t)$ into $\A(s)$ {\it canonical}, if
\item{(a)}
$\A(t)$ is mapped onto the initial $n$ component of $\A(s)$ of 
form $\A(s^+)$ if $t=s^+$, and 
\item{(b)}
$\A(t)$ is mapped onto one of the successor $n$ components of 
$\A(s)$ of form $\A(s^-)$, if $t=s^-$. 


\chapSkip
\Chapter{[2] Tails of first limit rank}

\secSkip\noindent
The modest goal here is a characterisation of the tails (of 
uniform chains) of rank $\omega$.

\secSkip
\Section{[2.1]}
\Title{Tail-types and successor-types of hyper-components}   

\parSkip\noindent
We start with a basic observation:
Let $A\in\UC$, $x\in A$ and let
    $C_n = \cmp{x}{n}{A}$, $n\in\rk{A}$.
Then
\dispM{
    \cmp{x}{n+1}{\f_x A}
    = \cmp{x}{n}{\f_x A} + A_n\N
    = \f_x C_n + A_n\N
    \eqno{(2.1)}
}
where
    $A_n$ is (a representative of the isomorphism class of)
    the successor $n$-component (to $C_n$) in $C_{n+1}$,
and
    $\f A_n = \f C_n$.

\parSkip
\Number{(2.1)}
\Title{Definition}
\Subtitle{$n$th $\alpha$ hyper-component}
If $\rk{A} = \alpha+1$, and $\A/\alpha\isom\N$ let
\dispM{\hbox{
    $\Hyp{\alpha}{n}{A} :=$ the $n$th $\alpha$ hyper-component
    of $A$.
}}
If $\rk{A} = \alpha+1$, and $A/\alpha\isom\N$ or $A/\alpha\isom\Z$ 
let
\dispM{\hbox{
    $\Hyp{\alpha}{n}{A,x} :=$ the $n$th $\alpha$ hyper-component
    of $A$ relative to $\cmp{x}{\alpha}{A}$
}}
(i.e., such that that
    $\Hyp{\alpha}{0}{A,x} = \cmp{x}{\alpha}{A}$).
Where $x$ is fixed in a context, write
    $\Hyp{\alpha}{n}{A}$ again for $\Hyp{\alpha}{n}{A,x}$.
If $A/\alpha\isom\N$, $x$ can always be arranged to lie in the 
zeroth hyper-component. In particular,
\dispM{
    \Hyp{\alpha}{n}{\cmp{x}{\alpha+1}{A}}
    := \Hyp{\alpha}{n}{\cmp{x}{\alpha+1}{A},x}.
}

\parSkip
\Number{(2.2)}
\Title{Lemma}
\Subtitle{Hyper-component isomorphism}
For $x,y\in A\in\UC$, $\alpha<\rk{A}$,
\dispM{
    (\forall n,m)
    [\f\Hyp{\alpha}{n}{\cmp{x}{\alpha+1}{A}}
     \isom\f\Hyp{\alpha}{m}{\cmp{y}{\alpha+1}{A}}],
    \qquad
    (\forall n,m > 0)
    [\Hyp{\alpha}{n}{\cmp{x}{\alpha+1}{A}}
     \isom\Hyp{\alpha}{m}{\cmp{y}{\alpha+1}{A}}]
}

\Proof
Fix $x$ in $A$. In
    $\f_x(\cmp{x}{\alpha+1}{A})$,
\dispM{
    \eqalign{
    &\hbox{
    $\f_x(\cmp{x}{\alpha}{A})
    = \f\Hyp{\alpha}{0}{\cmp{x}{\alpha+1}{A}}$
    is the zeroth $\alpha$ component,}
    \cr
    &\hbox{
    (for $n>0$) $\Hyp{\alpha}{n}{\cmp{x}{\alpha+1}{A}}$
    is the $n$th $\alpha$ component.}
    \cr
}}

For any $y$ in $A$, since
    $\f_xA \isom \f_yA$
and since isomorphisms preserve components and component
ordering,
\dispM{
    \eqalign{
    &\f\Hyp{\alpha}{0}{\cmp{x}{\alpha+1}{A}}
    \isom \f\Hyp{\alpha}{0}{\cmp{y}{\alpha+1}{A}}
    \cr
    &\Hyp{\alpha}{n}{\cmp{x}{\alpha+1}{A}}
    \isom\Hyp{\alpha}{n}{\cmp{y}{\alpha+1}{A}}
    \cr
    }\eqno{(a)}
}
wherever either is defined.

Fix $n$ such that 
    $\Hyp{\alpha}{n}{\cmp{x}{\alpha+1}{A}}$
is defined and let
    $z\in\Hyp{\alpha}{n}{\cmp{x}{\alpha+1}{A}}
      \subeq\cmp{x}{\alpha+1}{A}$. 
In
    $\f_z(\cmp{x}{\alpha+1}{A})$,
\dispM{
    \eqalign{
    &\hbox{
    $\f_z(\cmp{z}{\alpha}{A})
    = \f\Hyp{\alpha}{0}{\cmp{x}{\alpha+1},z}
    = \f\Hyp{\alpha}{n}{\cmp{x}{\alpha+1}}$
    is the zeroth $\alpha$ component,}
    \cr
    &\hbox{(for $k>0$)
    $\Hyp{\alpha}{k}{\cmp{x}{\alpha+1}{A},z}
    = \Hyp{\alpha}{n+k}{\cmp{x}{\alpha+1}{A}}$
    is the $m$th zeroth $\alpha$ component.}
    \cr
}}

Since $\cmp{x}{\alpha+1}{A}$ is uniform,
    $\f_x(\cmp{x}{\alpha+1}{A})
    \isom\f_z(\cmp{x}{\alpha+1}{A})$,
and since isomorphisms preserve components and component
ordering,
\dispM{
    \eqalign{
    &\f\Hyp{\alpha}{0}{\cmp{x}{\alpha+1}}
    \isom\f\Hyp{\alpha}{n}{\cmp{x}{\alpha+1}}
    \cr
    &\hbox{(for $n,k>0$) }
    \Hyp{\alpha}{n}{\cmp{x}{\alpha+1}{A}}
    \isom\Hyp{\alpha}{n+k}{\cmp{x}{\alpha+1}{A}}.
    \cr
    }\eqno{(b)}
}

Now, (a) allows us to 'vary' $x$ in the expressions
    $\f\Hyp{\alpha}{n}{\cmp{x}{\alpha+1}{A}}$
and
    $\Hyp{\alpha}{n}{\cmp{x}{\alpha+1}{A}}$
for fixed $n$ and $n>0$, respectively, while (b)
allows us to 'vary' $n$ and $n>0$ in them, respectively,
and always produce isomorphic components, respectively.
\Qed

\parSkip
We say a component is a {\it successor component} if it 
has an immediate predecessor in the ordering taken modulo
its rank, i.e., if the component is of form 
      $C = \cmp{x}{\alpha}{A}$
but $C/\alpha$ has an immediate predecessor in 
      $\cmp{x}{\alpha+1}{A}/\alpha$
 -- hence in any super-component of the latter in $A/\alpha$.
(For illustration, a component is a {\it gap component} if it does 
not have an immediate predecessor. We also then say that {\it there 
is a gap at $C$.\/} If a gap-component is of form
      $C = \cmp{x}{\alpha}{A}$,
the {\it depth} of the gap {\it relative to $C$\/} is 
    $\delta=\beta-\alpha$ 
where $\beta$ is minimal such that either
      $B=\cmp{x}{\beta}{A}$
has an immediate successor in $A/\beta$ (hence in
      $\cmp{x}{\beta}{A}$),
or that $\beta=\rk{A}$. (Here
    $\beta-\alpha = \sup\ens{\gamma:\alpha+\gamma\le\beta}$.)
For all $\gamma\in (\alpha,\beta)$, $C$ is the also {\it the
initial} $\alpha$ component in 
      $\cmp{x}{\gamma}{A}$.
The {\it absolute depth} of the gap determined by that component 
is $\beta-\gamma$ where $\gamma$ is minimal $\le\alpha$ such that 
for all $\delta\in [\gamma,\alpha]$, there is a $y$ such that 
    $D=\cmp{y}{\delta}{C}$ is minimal in $C/\delta$ 
 -- and hence has no immediate predecessor, hence is a gap
component. Note that for such $y$, $D$ is a minimal hyper-component
of
    $\cmp{y}{\delta+1}{A}$,
as well a minimal element in 
    $\cmp{y}{\gamma}{A}
    = \cmp{y}{\gamma}{B}$,
for any $\gamma\in[\delta.\beta)$.
Isomorphisms preserve not only immediate successor-ship but
as gaps and gap-depth as well.)
The preceding justifies:

\parSkip
\Number{(2.3)}
\Title{Definition}
\Subtitle{Tail and successor-component}

Let $\Tl^A_{\alpha}$ denote the common isomorphism type of the
tails of components of rank $\alpha$ in $A$, i.e., those 
representable as $\f\cmp{x}{\alpha}{A}$ for some $x\in A$ -- 
and hence as $\f\Hyp{\alpha}{n}{\cmp{x}{\alpha+1}{A}}$ for some 
$x\in A$ and $n$.

Let $\Suc^A_{\alpha}$ denote the common isomorphism type of the
successor components of rank $\alpha$ in $A$, i.e., those 
representable as $\Hyp{\alpha}{n}{\cmp{x}{\alpha+1}{A}}$ for some 
$x\in A$ and $n > 0$.

\parSkip
Thus, 

\parSkip
\Number{(2.4)}
\Title{Lemma}
\Subtitle{Tail lemma}
For $x\in A\in\UC$,
\dispM{
    \f_x A \isom \sum_{n\in\rk{A}} A_n\N,
    \qquad\hbox{where}\qquad
    A_n\isom\Suc^A_{n+1}\isom\Hyp{n}{1}{\cmp{x}{n+1}{A}}.
}

\Proof
Set $C_n = \cmp{x}{n}{A}$. The statement is equivalent to
    $\f_x C_n = \sum_{i\in n} A_i\N$,
for $n=\rk{A}$. This follows by induction over $n<\rk{A}$. 

For $n=1$, 
    $C_1\in\ens{\N,\Z}$
and
    $A_0\isom \ONE$,
so
    $\f C_1\isom\N\isom A_0\N$.

Assume the lemma true for all $i<\alpha\le\rk{A}$. 

First let $\alpha = n+1$. Then by (2.1)
    $\f_x C_{\alpha} = \f_x C_n + A_n\N$,
but by inductive hypothesis,
    $\f_x C_n = \sum_{i\in n} A_i\N$,
giving
     $\f_x C_{\alpha} = \sum_{i\in n} A_i\N + A_n\N$,
hence the required result.

Finally let $\alpha$ be a limit ordinal.
Then for each $n<\alpha$
     $\f_x C_n = \sum_{n\in n} A_n\N$,
and for every pair $n<k<\alpha$, $\f_x C_n$ embeds as 
initial segment into $\f_x C_k$, giving again the required
result.
\Qed

\secSkip
\Section{[2.2]}
\Title{Tails of rank $\omega$}   

\parSkip\noindent
The preceding allows to classify all $A\in\UT(\omega)$.

\parSkip
\Number{(2.5)}
\Title{Definition}
For $\ss\in\SS(\omega)$ set 
\dispm{    
    \A(\ss\tac{0}) = \sum_{i\in\omega}\A(\ss\ric i+1)\N
}

\parSkip
With $\A(\tup{}) = {\bf 0}$ (:= the empty chain), 
the above is tanatamount to setting 
    $\A(s\tac{0})=\sum_{i\in\omega}\A(s\ric i)\N$.

\parSkip
\Number{(2.6)}
\Title{Lemma}
\Subtitle{Rank $\omega$ tail representation lemma}
For every $A\in\UC(\omega)$ there exists 
$\ss\in\SS(\omega)$ such that 
\dispm{
    \f A=\A(\ss\tac{0}).
}

\Proof
Choose any $x\in A$ and let $B=\f_xA$. Fix 
    $s_n\in\SS(n+1)$
such that 
    $\A(s_n) = \cmp{x}{n}{B}$.  
Since 
    $\cmp{x}{n}{B} = \f_x\cmp{x}{n}{B}$,
in fact
    $s_n=s_n^-\tac{0}$.

Define a sequence $x_n\in A$ such that $x_0 = x$ and 
$\cmp{x_{n+1}}{n+1}{B}$ is the successor $n+1$ component
to $\cmp{x}{n+1}{B}$ in $A/\rho_{n+1}$.
Note that 
    $\cmp{x_{n+1}}{n+1}{B} = \cmp{x_{n+1}}{n+1}{A}$.
By uniformity of $A$ (and hence $B$),
    $\cmp{x}{n+1}{B} 
    = \cmp{x}{n}{B} + \cmp{x_n}{n}{B}\N
    = \cmp{x}{n}{B} + \cmp{x_n}{n}{A}\N$,
so
    $\A(s_{n+1}) = \A(s_n) + \A(t_n)\N$,
where
    $s_n = s_{n+1}^+$
and $t_n\in\SS(n+1)$ such that
    $t_n = s_{n+1}^- = s_{n+1}\ric n+1$.
and
    $s_n\ric n = s_n^- = t_n^- = t_n\ric n$.
Setting
    $\ss = \cup_{n\in\N}s_n\ric n \in\SS(\omega)$,
we get
    $t_n = s_{n+1}\ric n+1 = \ss\ric n+1$,
and thus
\dispM{
    \f_xA 
    = B
    = \sum_{n\in\N}\cmp{x_n}{n}{B}\N
    = \sum_{n\in\N}\cmp{x_n}{n}{A}\N
    = \sum_{n\in\N}\A(t_n)\N
    = \sum_{n\in\N}\A(\ss\ric n+1)\N
    = \A(\ss\tac{0}).
}
Let $y\in A$, $x \not= y$. Fix $k$ such that $x\rho_ky$.
Thus for $n\ge k$ the above procedure leads to the same
sequences $t_n$, hence to the same $s$.
\Qed

\parSkip
As an example, for every $\ss\in\SS(\omega)$ and every $n$,
we get
    $\f B(\ss,n) = \A(\ss\tac{0})$
for each of the chains
\dispM{
    B(\ss,n)
    = \A(\ss\ric n)\Z
    + \sum_{i\in\omega}\A(\ss\ric n+i+1)\N.
}

\parSkip
The attempt to characterise uniform chains of rank $\omega$ will
be resumed in section 5.


\chapSkip
\Chapter{[3] Chains of rank finite over a limit}

\secSkip\noindent
This is the main section of the paper. It describes a way to extend 
the Pastjin chain constructors to finite successors of limit ordinals 
$\lambda$ using the chains of rank $\lambda$ as 'building blocks'.

\secSkip
\Section{[3.1]}
\Title{Chains of rank $\omega+$finite}

\parSkip
\Number{(3.1)}
\Title{Definition}
\Subtitle{$\alpha$ successor}
Call $x\in A$ an {\it $\alpha$ successor (in $A$)} if
$\alpha$ is minimal $<\rk{A}$ such that
    $\co{x}{\alpha}$ has an immediate predecessor in 
    $\co{x}{\alpha+1}$.

Say that $x$ is an {\it $\alpha$ successor in $A/{\beta}$}, 
where $\alpha<\rk{A/\beta}$ and $\alpha<\rk{A/{\beta}}$, if
    $\co{x}{\beta}$ is an $\alpha$ successor in $A/{\beta}$.

\parSkip
\Title{Examples:}
\item{(a)}
Every successor element in $A$ (i.e., every element that has an
immmediate predecessor in $\co{x}{1}\in\ens{\N,\Z}$) is a
0-successor.
\item{(b)}
If $A$ has a least element $x$, then $x$ is a $\rk{A}$ successor.
\Qed

\parSkip
\Number{(3.2)}
\Title{Lemma}
\Subtitle{Successor induction}
Suppose $x$ is an $n$ successor in $\A(t)$, $t\in\SS$, where 
$n\in\ens{0,\ldots,\mod{t}-1}$ and let $t\in\ens{s^+,s^-}$
for some $s\in\SS$. Assume a canonical embedding of $\A(t)$ 
into $\A(s)$.   
\item{(a)}
$x$ is an $n+1$ successor in $\A(s)$
iff $n=\mod{t}-1$ and $t=s^+$.
\item{(b)}
$x$ is an $n$ successor in $\A(s)$ in all other cases, i.e.,
iff $n<\mod{t}-1$ or $t=s^-$.
\par\noindent
In other words,
    $x$ is an $n$ successor in $\A(s)$ iff
    $x$ is an $n$ successor in $\A(t)$, 
except when $n=\mod{t}-1$ and $t=s^+$, in which case 
    $x$ is an $n+1$ successor in $\A(s)$ iff
    $x$ is an $n$ successor in $\A(t)$.

\Proof
If $n=\mod{t}-1$ and $t=s^+$, then $x$ is minimal in $A(t)$
and $\A(s) =\A(t)+\A(s^-)\N$; so $x$ is minimal in $\A(s)$,
so $n$ is a $n+1=\rk{\A(s)}$ successor. Conversely, if $x$
is an $n+1$ succssor, then $x$ is minimal in $\A(s)$, and
since it is also minimal in $A(t)$ for $t\in\ens{s^+, s^-}$
we cannot have $t=s^-$ (since otherwise $\A(t)$ would be 
mapped onto one of the successor $n$ components of $\A(s)$
of form $\A(s)$, contradicting minimality of $x$ in $\A(s)$), 
so $t=s^+$.

If $n<\mod{t}-1$ it suffices to note that $\cmp{x}{i}{A(t)}
=\cmp{x}{i}{A(s)}$ for $i=n,n+1$ so $x$ is an $n$ successor
in $\A(t)$ iff it is an $n$ successor in $\A(s)$.

That leaves the case $n=\mod{t}-1$ and $t=s^-$. In this case
if $x$ is an $n$ successor in $\A(t)$ and $\rk{\A(t}=n$, then
$x$ in minimal in $\A(t)$; but $\A(t)$ is embedded as a
successor $n$ component into $\A(s)$, so $n$ is now minimal
such that $x$ is a successor mod $n$, so $x$ is also an $n$
successor in $\A(s)$. Conversely, if $x$ is an $n$ successor
in $\A(s)$, $x$ is minimal in some successor $n$ component
of $\A(s)$, i.e., $x$ is minimal in some $\A(s^-)$, so
$t=s^-$.
\Qed

\parSkip
By the component tail lemma, for any $x,x'\in A\in
\UC(\omega,\omega*2)$, 
    $\f\co{x}{\omega} \isom \f\co{x'}{\omega}$.

\parSkip
\Number{(3.3)}
\Title{Lemma}
\Subtitle{Successor isomorphism}
For $x,x'\in A\in\UC(\omega,\omega*2)$, 
if both $\co{x}{\omega}$ and $\co{x'}{\omega}$ 
are $n$ successors in $A/{\omega}$ for $n<\rk{A/\omega}<\omega$, 
then
    $\co{x}{\omega} \isom \co{x'}{\omega}$.

\Proof
Let $y,y'$ be immediate predecessors of $x, x'$ mod $\omega+n$
Now,
\dispM{
    \f_yA = \f_y\co{y}{\omega+n}
    +\sum_{i=0}^{m}\co{x}{\omega+n+i}\N
}
where $m=\rk{A/{\omega}}$. By the block lemma, the isomorphism
between $\f_yA$ and $\f_{y'}A$ 
maps the second $\omega+n$ component of $\f_yA$ isomorphically 
onto the second $\omega+n$ component of $\f_{y'}A$. 
Thus, 
    $\co{x}{\omega+n}\isom\co{x'}{\omega+n}$.
By finite induction over $i=n-1,\ldots,0$, using that
$x$ and $x'$ have no immediate predecessor mod $\omega+i$ for 
$i=0,\ldots,n-1$, it follows that 
    $\co{x}{\omega}\isom\co{x'}{\omega}$.
The induction step is simply to note that if
    $\co{x}{\omega+i+1}\isom\co{x'}{\omega+i+1}$,
then also for the respective initial components of these,
    $\co{x}{\omega+i}\isom\co{x'}{\omega+i}$,
by another application of the block lemma.
\Qed

\parSkip
\Number{(3.4)}
\Title{Definition}
\Subtitle{$\SS$-domains}
Define $\Dom(s)$ to be the closure of $\ens{s}$, ($s\in
\SS(<\omega)$) under $t\to t^+$ and $t\to t^-$; formally, 
    $\Dom(\tup{}) = \ens{\tup{}}$
and if $s^+\in\SS$,
    $\Dom(s) = \ens{s} \cup \Dom(s^+) \cup \Dom(s^-)$,
else (if $s^+\not\in\SS$)
    $\Dom(s) = \ens{s} \cup \Dom(s^-)$.

\parSkip
\Number{(3.5)}
\Title{Definition}
\Subtitle{$\SS$-trees}
The {\it binary tree $\Tree(s)$} of $s\in\SS(<\omega)$ 
is the tree of height $\mod{s}$ defined by
\item{(a)} [Level 0]
The top-node (root) of $\Tree(s)$ is $\tup{s;0}$, i.e.,
$\Lev_0\Tree(s) = \ens{ \tup{s;0} }.$
\item{(b)} [Level $i$ $\to$ level $i+1$]
Suppose $\sigma=\tup{t,k}\in\Lev_i\Tree(s)$, where 
$i<\mod{s}-1$.
\item{(b.0)} 
If $t^+\not\in\SS$, i.e., if $t(m)=m$ for $m=\mod{t}-1$, then 
$\sigma^0=\tup{t^-;0}$ is the only extension of $\sigma$ in
$\Lev_{i+1}\Tree(s)$.
\item{(b.1)}
If $t^+\in\SS$, i.e., if $t(m)<m$ for $m=\mod{t}-1$, then 
$\sigma^1=\tup{t^+;k+1}$ and $\sigma^0=\tup{t^-,0}$ are the only 
extensions of $\sigma$ in $\Lev_{i+1}\Tree(s)$.
%

\parSkip
Order the extensions of $\sigma\in\Tree(s)$ by setting
    $\sigma^1 < \sigma^0$ 
whenever $\sigma^1$ is defined.
This will induce an ordering on each level of $\Tree(s)$, and
hence on the leaves of $\Tree(s)$, by setting in addition
    $\sigma^{0,1} < \tau^{0,1}$
on level $k+1$ whenever $\sigma < \tau$ on level $k$.

\parSkip
\Title{Notation}
 Denote by
    $\leaf^s(i)$
the $i$th leaf of $\Tree(s)$.

\parSkip
\Title{Corollary}
\item{(a)}
$\Tree(s^-)$ is a sub-tree of $\Tree(s)$. - Hence:
\item{(b)}
$\Tree(s^{-k})$ is a sub-tree of $\Tree(s)$, for every
$k=1,\ldots,\mod{s}-1$. - More specifically:
\item{(c)}
Let $k$ be maximal $\le\mod{s}$ such that $s^{+k}\in\SS$.
Then $\Tree(s^{-i})$ is the sub-tree of $\Tree(s)$ extending
node $\tup{s^{+i};i}$ via node $\tup{s^{(+i)-},0}$ ('to the right'), 
for every $i=1,\ldots,k$, and $\Tree(s^{-(k+1)})$ is the sole
extension of $\Tree(s)$ at node $\tup{s^{+k};k}$.

\Proof
(a) The top-node of $\Tree(s^-)$ is $\tup{s^-;0}$, which an extension
(to the 'right', if $s^+\in\SS$, else the only) of the top-node 
$\tup{s;0}$ of $\Tree(s)$. (b) follows by induction, and (c) follows
from $s^{(+k)-} = s^{-(k+1)}$.
\Qed

\parSkip
\Number{(3.6)}
\Title{Definition}
\Subtitle{Association}
For each $x\in A = \A(s)$ define the {\it association} between
$x$ and nodes $\sigma$ on $\Tree(s)$ as follows:
\item{(a)}
    $x\assoc\tup{s,0}$
\item{(b)}
Supposing $x\assoc\sigma=\tup{t,k}\in\Tree(s)$ for
$\mod{t} > 1$. Then
\itemitem{(b.0)}
$x\assoc\sigma^1=\tup{t^+,k+1}$ if 
$\co{x}{\mod{t}-1}$ is the initial $\mod{t}-1$ 
component of $\co{x}{\mod{t}}$;
\itemitem{(b.1)}
$x\assoc\sigma^0=\tup{t^-,0}$ if 
$\co{x}{\mod{t}-1}$ is a successor $\mod{t}-1$ 
component of $\co{x}{\mod{t}}$.

\parSkip
Also say that $x$ is associated to a branch in $\Tree(s)$
if $x$ is associated to each node of that branch.

\vfill\eject 
\Title{Examples}
%
%

\dispM{
\Tree(\tup{0})
\qquad\left\{\qquad
\tup{\tup{0};0}
\right.
}

\dispM{\eqalign{
\Tree(\tup{0,0})
&
\qquad\left\{\qquad
\mrx{ \tup{\tup{0,0};0} }
    { \tup{\tup{0};1} }
    { \Tree(\tup{0}) }
\right.
\cr&\quad\cr
\Tree(\tup{0,1})
&
\qquad\left\{\qquad
\nrx{ \tup{\tup{0,1};0} }
    { \Tree(\tup{0}) }
\right.
\cr
}}

\dispM{\eqalign{
\Tree(\tup{0,0,0})
&
\qquad\left\{\qquad
\mrx{ \tup{\tup{0,0,0};0} }
    { \mrx{ \tup{\tup{0,0};1} }
          { \tup{\tup{0};2} }
          { \Tree(\tup{0}) } }
    { \trx{ \Tree(\tup{0,0}) } }
\right.
\cr &\quad\cr
\Tree(\tup{0,0,1})
&
\qquad\left\{\qquad
\mrx{ \tup{\tup{0,0,1};0} }
    { \nrx{ \tup{\tup{0,1};1} }
          { \Tree(\tup{0}) } }
    { \trx{ \Tree(\tup{0,0}) } }
\right.
\cr&\quad\cr
\Tree(\tup{0,1,0})
&
\qquad\left\{\qquad
\mrx{ \tup{\tup{0,1,0};0} }
    { \mrx{ \tup{\tup{0,0};1} }
          { \tup{\tup{0};2} }
          { \Tree(\tup{0}) } }
    { \trx{ \Tree(\tup{0,1}) } }
\right.
\cr
}}

\dispM{\eqalign{
\Tree(\tup{0,0,0,1})
&
\quad\left\{\quad
\mrx{ \tup{ \tup{0,0,0,1};0} }
    { \mrx{ \tup{\tup{0,0,1};1} }
          { \nrx{ \tup{\tup{0,1};2} }
                { \Tree(\tup{0}) } }
          { \trx{ \Tree(\tup{0,0}) } } }
    { \trx{ \trx{ \Tree(\tup{0,0,0}) } } }
\right.
\cr&\quad\cr
\Tree(\tup{0,0,1,0})
&
\quad\left\{\quad
\mrx{ \tup{ \tup{0,0,1,0};0} }
    { \mrx{ \tup{\tup{0,0,0};1} }
          { \mrx{ \tup{\tup{0,0};2} }
                { \tup{\tup{0};3} }
                { \Tree(\tup{0}) } }
          { \trx{ \Tree(\tup{0,0}) } } }
    { \trx{ \trx{ \Tree(\tup{0,0,1}) } } }
\right.
\cr
}}

\vfill\eject 
Let $x\in\A(t)$ where $t\in\SS$ and assume $t\in\ens{s^+,s^-}$
for some $s\in\SS$ and suppose $\A(t)$ canonically embedded in
$\A(s)$. Suppose $x$ is associated to the nodes of a branch
in $\Tree(t)$ containing a node $\tup{u,0}$, $u\in\Dom(t)\sub\Dom(s)$.
Then $x$ is also associated to node $\tup{u,0}$ in $\Tree(s)$.
Namely, in $\Tree(s)$, $x$ is associated to $\tup{t,1}$ if
$t=s^+$, else to $\tup{t,0}$. In the latter case it follows
immediately that the branch in $\Tree(s)$ extending $\tup{t,0}$
to which $x$ is associated is the same as the branch in $\Tree(t)$.
If $t=s^+$ these branches differ in the initial segment of nodes
$\sigma$ for which $x$ gets associated to the left successor
$\sigma^1$ of $\sigma$. I.e., let $i$ be such that $x\assoc
\tup{t^{+i};i}$ and $x\assoc\tup{t^{(+i)-};0}$. Such a node must
exist if $x$ is associated to some node of form $\tup{u,0}$.
Then, in $\Tree(t)$, 
\item{(a)}
$x\assoc\tup{t^{+j};j}$ for $j=0,\ldots,i$,
$x\assoc\tup{t^{(+i)-};0}$, 
\par\noindent
while in $\Tree(s)$, 
\item{(b)}
$x\assoc\tup{t^{+j};j+1}=\tup{s^{+(j+1)};j+1}$ for $j=0,\ldots,i$ 
if $t=s^+$;
\item{(c)}
$x\assoc\tup{t^{+j};j}=\tup{s^{-(+j)};j}$ for $j=0,\ldots,i$ 
if $t=s^-$.
\par\noindent
But in either tree, $x\assoc\tup{t^{(+i)-};0}=\tup{s^{\pm(+i)-};0}$,
so the branches associated to $x$ extending the latter nodes in 
$\Tree(t)$ and $\Tree(s)$, respectively, are identical.
In particular, the branches associated to $x$ extending $\tup{u,0}$ 
in either $\Tree(t)$ or $\Tree(s)$ are identical.

\parSkip
\Title{Notation}
\Subtitle{Leaf number}
Let $\lf(s)$ denote the number of leaves in $\Tree(s)$.
Note that
    $\lf(s) = \lf(s^+) + \lf(s^-)$ if $s^+$ is defined
and
    $\lf(s) = \lf(s^-)$ else.
$\lf(s)$ will figure as the length of certain arrays to be
defined.

A basic list of $\lf(s)$ examples is:
$$\vbox{
\settabs 2\columns
\+$\lf(\tup{0}) = 1$&\cr
\+$\lf(\tup{0,0}) = \lf(\tup{0}) + \lf(\tup{0}) = 1+1 = 2$&
  $\lf(\tup{0,1}) = \lf(\tup{0}) = 1$\cr
\+$\lf(\tup{0,0,0}) = \lf(\tup{0,0}) + \lf(\tup{0,0}) = 2+2 = 4$&
  $\lf(\tup{0,0,1}) = \lf(\tup{0,1}) + \lf(\tup{0,0}) = 1+2 = 3$\cr
\+$\lf(\tup{0,0,2}) = \lf(\tup{0,0}) = 2$&
  $\lf(\tup{0,1,0}) = \lf(\tup{0,0}) + \lf(\tup{0,1}) = 2+1 = 3$\cr
\+$\lf(\tup{0,1,1}) = \lf(\tup{0,1}) + \lf(\tup{0,1}) = 1+1 = 2$&
  $\lf(\tup{0,1,2}) = \lf(\tup{0,1}) = 1$\cr
}$$

\parSkip
\Number{(3.7)}
\Title{Definition}
\Subtitle{Weight}
For a node $\sigma\in\Tree(s)$ set $\wt^s(\sigma) = k$
just if $\sigma = \tup{t,k}$ for some $t\in\Dom(s)$. 
For $x\in\A(s)$ set $\wt^s(x)=\wt^s(\sigma)$ for the
leaf $\sigma\in\Lev_{\mod{s}}(\Tree(s))$ associated to $x$.
For $i<\lf(s)$ set $\wt^s(i)=\wt^s(\leaf^s(i))$.

\parSkip
\Number{(3.8)}
\Title{Lemma}
\Subtitle{Weight induction}
Suppose $x\in\A(t)$ has weight $n$, where 
$t\in\SS$ and $n\in\ens{0,\ldots,\mod{t}-1}$ and suppose 
$t\in\ens{s^+,s^-}$ for some $s\in\SS$. Then 
for all $n\in\ens{0,\ldots,\mod{t}-1}$,
    $\wt^s(x) = n$ iff
    $\wt^t(x) = n$, 
except when $n=\mod{t}-1$ and $t=s^+$, in which case 
    $\wt^s(x) = n+1$ iff
    $\wt^t(x) = n$.

\Proof
Let $t\in\ens{s^+,s^-}$, $n=\mod{t}-1$ and $x\in\A(t)$.
Then
\dispM{
    \eqalign{
    \wt^t(x)=n\quad
    &\hbox{iff $x\assoc\sigma$ for some $\sigma\in\Lev_n(\Tree(t))$ 
        such that $\wt^t(\sigma)=n$} \cr
    &\hbox{iff $x\assoc\sigma=\tup{\tup{0};n}$} \cr
    &\hbox{iff $t^{+n}=\tup{0}\in\SS$ and $x\assoc\tup{t^{+n};n}
        =\tup{\tup{0};n}$} \cr
    &\hbox{iff $t=t^-\tac{0}$ and $x\assoc\tup{t^{+n};n}
        =\tup{\tup{0};n}$} \cr
    &\hbox{iff $x\assoc\tup{t^{+i};i}$, for each 
        $i=0,\ldots,n$, } \cr
    }
}
i.e., $x$ is associated to each node of the left-most branch of 
$\Tree(t)$, given by 
    $\tup{t,0},\tup{t^+,1},\ldots,\tup{t^{+i};i},\ldots,
     \tup{\tup{0};n}$.

Case (a): Suppose the $\wt^t(x)=n=\mod{t}-1$.

Case (a.0): If also $t=s^+$, then
$x$ is associated to each node $\tup{s^{+i};i}$, for $i=0,\ldots,n+1$,
of the left-most branch of $\Tree(s)$, given by
    $\tup{s,0},\tup{t,1},\tup{t^+,2},\ldots,\tup{t^{+i},i+1},\ldots,
    \tup{\tup{0};n+1}$, and so $\wt^s(x)=n+1$.

Case (a.1): If on the other hand $t=s^-$, then 
$x$ is associated to nodes 
    $\tup{s,0}, \tup{s^-,0}=\tup{t,0}$,
and hence each node 
    $\tup{s^{-(+i)};i}=\tup{t,i}$ for $i=0,\ldots,n$ 
of the branch
    $\tup{s,0},\tup{t,0},\tup{t^+,1},\ldots,\tup{t^{+i};i},\ldots,
    \tup{\tup{0};n}$
of $\Tree(s)$, and so $\wt^s(x)=n$.

Case (b): Suppose on the other hand that $\wt^t(x)=n<\mod{t}-1$.
Then there is some least $i<n$ such that $x$ is associated 
with each node of the branch
    $\tup{t,0},\ldots\tup{t^{+i};i},\tup{t^{(+i)-};0}$. 
(I.e., for some $i<n$, $x$ must get associated with 
$\tup{t^{(+i)-};0}$, else we are back in case (a.0). For
the least such $i$, $x$ is associated with the branch just
described.) Hence in $\A(s)$, $x$ is associated to each node 
of the branch
    $\tup{s,0},\tup{t,k_0},\ldots\tup{t^{+i};k_i},\tup{t^{(+i)-};0}$,
where 
    $k_j=j+1$ if $t=s^+$ 
and 
    $k_j=j$ if $t=s^-$. 
The remaining nodes of the branches in either $\Tree(t)$ or 
$\Tree(s)$ with which $x$ gets associated are now identical. 
Thus $\wt^t(x)= \wt^s(x)$, as required.
\Qed

\parSkip
\Title{Corollary}
Let $s\in\SS$ and let $n=\mod{s}-1$.
\item{(a)}
\dispM{
    \wt^s(0)=\cases{
    \wt^{s^+}(0)+1 & if $s^{+n}=\tup{0}\in\SS$ 
                     (i.e. $s=s^-\tac{0}$) \cr
    \wt^{s^+}(0)   & else \cr
    }
}
\item{(b)}
For $0<i<\lf(s)$,
\dispM{
    \wt^s(i)=\cases{
    \wt^{s^+}(i) 
    & if $s^+\in\SS$ and $i<\lf(s^+)$ \cr
    \wt^{s^-}(\lf(s^+)+i) 
    & if $s^+\in\SS$ and $\lf(s^+)\le i<\lf(s)=\lf(s^+)+\lf(s^-)$ \cr
    \wt^{s^-}(i) 
    & if $s^+\not\in\SS$ and $i<\lf(s)=\lf(s^-)$ \cr
    }
}

\parSkip
\Title{Corollary}
\item{(a)}
If $s=s^-\tac{0}$ (and so $s^+\in\SS$), then for all $i<\lf(s)
=\lf(s^+)+\lf(s^-)$
\dispM{
    \wt^s(i)=\cases{
    \wt^{s^+}(i) + 1
    & for $i=0$ \cr
    \wt^{s^+}(i) 
    & for $0<i<\lf(s^+)$ \cr
    \wt^{s^-}(\lf(s^+)+i) 
    & for $\lf(s^+)\le i<\lf(s)$ \cr
    }
}
\item{(b)}
If $s^+\in\SS$ but $s\not=s^-\tac{0}$, then for all $i<\lf(s)
=\lf(s^+)+\lf(s^-)$
\dispM{
    \wt^s(i)=\cases{
    \wt^{s^+}(i) 
    & for $i<\lf(s^+)$ \cr
    \wt^{s^-}(\lf(s^+)+i) 
    & for $\lf(s^+)\le i<\lf(s)$ \cr
    }
}
\item{(c)}
If $s^+\not\in\SS$, then for all $i<\lf(s)=\lf(s^-)$
\dispM{
    \wt^s(i)= \wt^{s^-}(i) 
    \hbox{ for $i<\lf(s)=\lf(s^-)$ } 
}

\vfill\eject
\Title{Examples:}

(i) Compare $\Tree(s)$ for 
    $s=\tup{0,0,1}$
or
    $s=\tup{0,0,0,1}$
with those for 
    $s=\tup{0,1,0}$
or
    $s=\tup{0,0,1,0}$
given above. The latter are examples in which case (a) of the 
preceding corollary applies.

(ii) Compare $\Tree(s)$ for 
    $s=\tup{0,0,2,1}$, 
    $s=\tup{0,0,1,2}$ 
and
    $s=\tup{0,0,1,3}$
with 
    $\Tree(\tup{0,0,1})$.

\dispM{
\Tree(\tup{0,0,2,1})
\qquad\left\{\qquad
\mrx{ \tup{\tup{0,0,2,1};0} }
    { \mrx { \tup{\tup{0,0,1};1} }
           { \nrx{ \tup{\tup{0,1};1} }
                 { \Tree(\tup{0}) } }
           { \trx{ \Tree(\tup{0,0}) } } }
    { \trx { \trx{ \Tree(\tup{0,0,2}) } } }
\right.
}

\dispM{
\Tree(\tup{0,0,1,2})
\qquad\left\{\qquad
\mrx{ \tup{\tup{0,0,1,2};0} }
    { \nrx { \tup{\tup{0,0,2};1} }
           { \Tree(\tup{0,0}) } }
    { \trx { \Tree(\tup{0,0,1}) } } 
\right.
}

\dispM{
\Tree(\tup{0,0,1,3})
\qquad\left\{\qquad
\nrx{ \tup{\tup{0,0,1,3};0} }
    { \Tree(\tup{0,0,1}) }
\right.
}

(iii) Compare $\Tree(s)$ for 
    $s=\tup{0,1,0,1}$,
    $s=\tup{0,1,0,2}$,
and 
    $s=\tup{0,1,0,3}$ 
with
    $\Tree(\tup{0,1,0})$.

\dispM{
\Tree(\tup{0,1,0,0})
\qquad\left\{\qquad
\mrx{ \tup{\tup{0,1,0,0};0} }
    { \mrx{ \tup{\tup{0,1,0};1} }
          { \mrx{ \tup{\tup{0,0};2} }
                { \tup{\tup{0};3} }
                { \Tree(\tup{0}) } }
          { \trx{ \Tree(\tup{0,1}) } } }
    { \trx{ \trx{ \Tree(\tup{0,1,0}) } } }
\right.
}

\dispM{
\Tree(\tup{0,1,0,1})
\qquad\left\{\qquad
\mrx{ \tup{\tup{0,1,0,1};0} }
    { \mrx{ \tup{\tup{0,1,1};1} }
          { \nrx{ \tup{\tup{0,1};1} }
                { \Tree(\tup{0}) } }
          { \trx{ \Tree(\tup{0,1}) } } }
    { \trx{ \trx{ \Tree(\tup{0,1,0}) } } }
\right.
}

\dispM{
\Tree(\tup{0,1,0,3})
\qquad\left\{\qquad
\nrx{ \tup{\tup{0,1,0,3};0} }
    { \Tree(\tup{0,1,0}) }
\right.
}

\vfill\eject
\Number{(3.9)}
\Title{Lemma}
\Subtitle{Weight isomorphism}
\item{(a)}
Every $x\in\A(s)$ is a $\wt^s(x)$ successor.
\item{(b)}
Suppose $A/{\omega}\isom\A(s)$. 
Then for all $x,x'\in A$, if
    $\wt^s(\co{x}{\omega}) = \wt^s(\co{x'}{\omega})
    <\rk{A/\omega}=\rk{\A(s)}<\omega$ 
then
    $\co{x}{\omega} \isom \co{x'}{\omega}$.

\Proof
(i)
By induction on $\rk{\A(s)} = \mod{s}-1$.

For $s=\tup{0}$ the claim is vacuously true. 

For $\mod{s}=2$,
$\wt^s(x) = 1$ iff
$x\assoc\tup{\tup{0};1}=\tup{s,0}^+$ iff 
$x$ is minimal in $\A(s)$ iff 
$x$ is a $\rk{\A(s)} = 1$ successor,
and $\wt^s(x) = 0$ iff $x$ is a 0 successor (i.e., has
an immediate predecessor in $\co{x}{1}=\N$).

For $s=\tup{0,0}$, there is exactly one $x\in\A(s)=\N$ for 
which each statement in the first chain of implications holds.
For any other element in $\N$, and any element in $\A(\tup{0,1})
=\Z$, both statements of the second implication apply.

Let $s\in\SS$ and suppose the statement is true 
for $t\in\ens{s^+,s^-}\cap\SS$, i.e., 
$x$ is an $n$ successor in $\A(t)$ iff $\wt^t(x)=n$, 
for all $n\in\ens{0,\ldots,\mod{s}-1}$.
Assume $\A(t)$ embeds canonically into $\A(s)$. 
By the successor induction lemma, 
for all $n\in\ens{0,\ldots,\mod{t}-1}$,
$x$ is an $n$ successor in $\A(s)$ iff
    $x$ is an $n$ successor in $\A(t)$, 
except when $n=\mod{t}-1$ and $t=s^+$, in which case 
    $x$ is an $n+1$ successor in $\A(s)$ iff
    $x$ is an $n$ successor in $\A(t)$.
By the weight induction lemma, for all such $n$,
    $\wt^s(x) = n$ iff
    $\wt^t(x) = n$, 
except when $n=\mod{t}-1$ and $t=s^+$, in which case 
    $\wt^s(x) = n+1$ iff
    $\wt^t(x) = n$.
Thus, by inductive hypothesis, $x$ is an $n$ successor in
$\A(s)$ iff $\wt^s(x)=n$, for all $n\in\ens{0,\ldots,\mod{s}-1}$.

(ii) 
Follows from (i) and the successor isomorphism lemma:
Suppose $A/{\omega}\isom\A(s)$. If
    $\wt^s(\co{x}{\omega}) = \wt^s(\co{x'}{\omega}) = n$ 
for $x,x'\in A$ and $n<\rk{A/\omega}\in\omega$, then by (i),
    $\co{x}{\omega}$ and $\co{x'}{\omega}$ 
    are $n$ successors in $A/\omega$.
Again by the successor isomorphism lemma,
    $\co{x}{\omega} \isom \co{x'}{\omega}$,
as required.
\Qed

\parSkip
\Number{(3.10)}
\Title{Definition}
\Subtitle{Plus finite rank Pastjin uniform chain constructor ---
long version}
Let $\M_i\in\UC(\omega),i\in\omega$. Define
\dispM{
    \A(\tup{0})[\M_0] = 1 [\M_0] = \M_0
}
and for $\lf(s)$ arrays 
    $\tup{\M_0,\ldots,\M_{\lf(s)-1}}\in\Fun{\lf(s)}{\UC(\omega)}$,
\dispM{
    \eqalign{
    \A(s)\arr{\M_0,\ldots,\M_{\lf(s)-1}}
    &= \cases{
    \A(s^-)\arr{\M_0,\ldots,\M_{\lf(s^-)-1}}.\Z
    &if $s(n)=n$ \cr
    \A(s^+)\arr{\M_0,\ldots,\M_{\lf(s^+)-1}}
    + \A(s^-)\arr{\M_{\lf(s^+)},\ldots,\M_{\lf(s)-1}}.\N
    &if $s(n)<n$ \cr} \cr
    }
}
where $n=\mod{s}-1$ (and where the '.' is meant to regulate 
precedence of operations in the usual fashion). 
The sequence 
    $\MM 
    = \MM_{[0,\lf(s))}
    = \tup{\M_0,\ldots,\M_{\lf(s)-1}}$
is subsequently called the {\it array} (of building blocks) of the 
constructor, and 
    $\A(s^+)\arr{\M_0,\ldots,\M_{\lf(s^+)-1}}$,
respectively
    $\A(s^-)\arr{\M_0,\ldots,\M_{\lf(s^-)-1}}$
or
    $\A(s^-)\arr{\M_{\lf(s^+)},\ldots,\M_{\lf(s)-1}}$
(whichever applies) are referred to as the 'plus', respectively 
'minus' sub-components of
    $\A(s)\arr{\M_0,\ldots,\M_{\lf(s)-1}}$.

\parSkip
Note that $\A(s)\arr{\MM}/\omega\isom\A(s)$.

\parSkip
A basic list of such chains of rank finite over $\omega$ 
contructed with $\lf(s)$ arrays is:
$$\vbox{
\settabs 2\columns
\+$\A(\tup{0})[\M_0] 
    = \M_0$&\cr
\+$\A(\tup{0,0})[\M_0,\M_1] 
    = \M_0 + \M_1\N$&
  $\A(\tup{0,1})[\M_0] 
    = \M_0\Z$\cr
\+$\A(\tup{0,0,0})[\M_0,\M_1,\M_2,\M_3]
    = \M_0 + \M_1\N + (\M_2 + \M_3\N)\N$&
  \cr
\+&
  $\A(\tup{0,0,1})[\M_0,\M_1,\M_2]
    = \M_0\Z + (\M_1 + \M_2\N)\N$\cr
\+$\A(\tup{0,0,2})[\M_0,\M_1]  
    = (\M_0 + \M_1\N)\Z$&
  $\A(\tup{0,1,0})[\M_0,\M_1,\M_2]  
    = \M_0 + \M_1\N + \M_2\Z\N$\cr
\+$\A(\tup{0,1,1})[\M_0,\M_1]  
    = \M_0\Z + \M_1\Z\N$&
  $\A(\tup{0,1,2})[\M_0]  
    = \M_0\Z^2$\cr
}$$

\parSkip
Arbitrary chains of the above form will not in general be uniform
again (even with $\M_i\in\UC(\omega)$). 
\item{(1)}
Order-isomorphisms between tails of chains induce (restrict 
to) order isomorphisms between the initial components of these
tails. The set of initial $\omega$ components of tails of
    $A = \A(s)\arr{\M_0,\ldots,\M_{\lf(s)-1}}$
is the set of tails
    $\f\M_0,\ldots,\f\M_{\lf(s)-1}$.
These must be pairwise isomorphic, if $A$ is uniform. 
\item{(2)}
Order-isomorphisms between tails of chains induce isomorphisms 
between the components that are 'successors' to the initial 
components of the tail chains. This yields a partition of the 
set
    $\M_0,\ldots,\M_{\lf(s)-1}$
into subsets of pairwise isomorphic members. For instance,
any two $\M_i,\M_j$ that are be successor $\omega$ components
in (tails) of $A$ must be isomorphic, any two $\M_i,\M_j$ that
can be initial $\omega$ components in successor $\omega+1$ 
components in tails of $A$ must be isomorphic, etc.

\parSkip
\Number{(3.11)}
\Title{Lemma}
\Subtitle{Uniformity over limit rank}
If
    $A = \A(s)[\M_0,\ldots,\M_{\lf(s)-1}]$ 
is uniform, then
\dispM{
    \eqalign{
    \hbox{(a)}\qquad &
    \hbox{$\f\M_i = \f\M_j$ for all $i,j<\lf(s)$, and} \cr
    \hbox{(b)}\qquad &
    \hbox{$\M_i = \M_j$ whenever $\wt^s(i) = \wt^s(j)$.} \cr
    } \eqno{(3.1)}
}

\Proof
(a) is immediate by the fact that for any $x,x'\in A$,
    $\f\co{x}{\omega} \isom \f\co{x'}{\omega}$
by the component tail isomorphism lemma. 

(b) By the weight isomorphism lemma it suffices to show that
    $\wt^s(\M_i) = \wt^s(i)$ for all $i<\lf(s)$.
For this in turn it suffices to show that 
    $\M_i\assoc\leaf^s(i)$ for all $i<\lf(s)$.

This follows by induction on $\rk{\A(s)}$, $s\in\SS(<\omega)$.
By definition of the association of $x\in\A(s)$ to nodes
$\sigma\in\Tree(s)$,
\item{(a)}
    $\M_i\assoc\tup{s,0}\in\Tree(s)$ for all $i<\lf(s)$,
and
\item{(b.0)}
if $s^+\in\SS$, then
    $\M_i\assoc\tup{s^+,1}\in\Tree(s)$ for all $i<\lf(s^+)$,
    $\M_i\assoc\tup{s^-,0}\in\Tree(s)$ for all $\lf(s)^+\le i<\lf(s)$,
i.e.,
    $\M_{lf(s^+)+i}\assoc\tup{s^-,0}\in\Tree(s)$ for all $i<\lf(s^-)$,
while
\item{(b.1)}
if $s^+\not\in\SS$, then
    $\M_i\assoc\tup{s^-,0}\in\Tree(s)$ for all $i<\lf(s)=\lf(s^-)$.

Suppose $s^+\in\SS$, then,
if 
    $\leaf^{s^+}(i)=\tup{\tup{0};k}$ for $i<\lf(s^+)$,
then
    $\leaf^{s}(i)=\tup{\tup{0};k+1}$ for $i=0$,
and
    $\leaf^{s}(i)=\tup{\tup{0};k}$ for $0<i<\lf(s^+)$;
and if 
    $\leaf^{s^-}(i)=\tup{\tup{0};k}$ for $i<\lf(s^-)$,
then
    $\leaf^{s}(\lf(s^+)+i)=\tup{\tup{0};k}$ for $i<\lf(s^-)$.

Suppose $s^+\in\SS$, then,
if 
    $\leaf^{s^-}(i)=\tup{\tup{0};k}$ for $i<\lf(s^-)$,
then
    $\leaf^{s}(i)=\tup{\tup{0};k}$ for $i<\lf(s)$.

In case (b.1) it follows immediately from the inductive hypothesis
that 
    $\wt^s(M_i)=\wt^{s^-}(M_i)=\wt^{s^-}(i)=\wt^s(i)$.

In case (b.0), 
\item{(a)}
    $\wt^s(M_0)=\wt^{s^+}(M_0)+1=\wt^{s^-}(0)+1=\wt^s(0)$,
\item{(b)}
    $\wt^s(M_i)=\wt^{s^+}(M_i)+1=\wt^{s^+}(i)=\wt^s(i)$
    for $0<i<\lf(s^+)$, and
\item{(c)}
    $\wt^s(M_{\lf(s^+)+i})=\wt^{s^+}(M_i)
    =\wt^{s^+}(i)=\wt^s({\lf(s^+)+i})$
    for $\lf(s^+)\le i\le\lf(s)$,
\par\noindent
again by inductive hypothesis.
\Qed

\parSkip
\Title{Examples:}
By earlier calculations, 
\item{(i)}
    $\lf(\tup{0,1,0})
    =\lf(\tup{0,0})+\lf(\tup{0,1})
    =(\lf(\tup{0})+\lf(\tup{0}))+\lf(\tup{0})=3,$
\item{(ii)}
    $\lf(\tup{0,0,1})
    =\lf(\tup{0,1})+\lf(\tup{0,0})
    =\lf(\tup{0})+(\lf(\tup{0})+\lf(\tup{0}))=3.$
\par\noindent
Thus, for (i),
\dispM{
\eqalign{
A(\tup{0,1,0})\arr{\M_0,\M_1,\M_2}
&=\A(\tup{0,0})[\M_0,\M_1] + \A(\tup{0,1})[\M_2].\N \cr
&=(1+\N)[\M_0,\M_1] + \Z[\M_2].\N \cr
&=\M_0 + \M_1\N + \M_2\Z\N,\cr
}}
\dispM{
\Tree(\tup{0,1,0})
\qquad\left\{\qquad
\mrx{ \tup{\tup{0,1,0};0} }
    { \mrx{ \tup{\tup{0,0};1} }
          { \tup{\tup{0};2} }
          { \tup{\tup{0};0} } }
    { \nrx{ \tup{\tup{0,1};0} }
          { \tup{\tup{0};0} } }
\right.
}
So,
    $\wt^{\tup{0,1,0}}(\M_0)=2$,
    $\wt^{\tup{0,1,0}}(\M_1)=\wt^{\tup{0,1,0}}(\M_2)=0$,
so
    $\M_1\isom\M_2$.
Hence, the uniform sets contucted from $\A(\tup{0,1,0})$ over
uniform rank $\omega$ chains are of form
    $A(\tup{0,1,0})\arr{\M_0,\M_1,\M_1}
    =\M_0 + \M_1\N + \M_1\Z\N$.
For (ii):
\dispM{
\eqalign{
A(\tup{0,0,1})\arr{\M_0,\M_1,\M_2}
&=\A(\tup{0,1})[\M_0] + \A(\tup{0,0})[\M_1,\M_2].\N \cr
&=\Z[\M_0] + (1+\N)[\M_1,\M_2].\N \cr
&=\M_0\Z + (\M_1 + \M_2\N)\N,\cr
}}
\dispM{
\Tree(\tup{0,0,1})
\qquad\left\{\qquad
\mrx{ \tup{\tup{0,0,1};0} }
    { \nrx{ \tup{\tup{0,1};1} }
          { \tup{\tup{0};0} } }
    { \mrx{ \tup{\tup{0,0};0} }
          { \tup{\tup{0};1} }
          { \tup{\tup{0};0} } }
\right.
}
So,
    $\wt^{\tup{0,0,1}}(\M_0)=\wt^{\tup{0,0,1}}(\M_2)=0$,
    $\wt^{\tup{0,0,1}}(\M_1)=1$,
so
    $\M_0\isom\M_2$.
Hence, the uniform sets contucted from $\A(\tup{0,0,1})$ over
uniform rank $\omega$ chains are of form
    $A(\tup{0,0,1})\arr{\M_0,\M_1,\M_0}
    =\M_0\Z + (\M_1 + \M_0\N)\N$.
\Qed

\parSkip
To show that condition (a) and (b) is also sufficient
for uniformity of $A$ we need to calculate its tails,
$\f A$.


\secSkip
\Section{[3.2]}
\Title{Tails of rank $\omega+$finite}   


\parSkip
\Number{3.12}
\Title{Definition}
\Subtitle{Tail array}
Given an array 
    $\MM_{\lf(s)}
    =\arr{\M_0,\ldots,\M_{\lf(s)-1}}$,
define the array
    $\f_s(\MM_{\lf(s)})
    $
by induction over $\Dom(s)$:

First, define
\dispM{
\eqalign{
    \f_{\tup{0}}(\MM_{\lf(\tup{0)}}) 
    &= \f(\arr{\M_0}) 
    = \arr{\f\M_0} \cr
    \f_{\tup{0,0}}(\MM_{\lf(\tup{0,0)}}) 
    &= \f(\arr{\M_0, \M_1}) 
    = \arr{\f\M_0, \M_1} \cr
    \f_{\tup{0,1}}(\MM_{\lf(\tup{0,1)}}) 
    &= \f(\arr{\M_0}) 
    = \arr{\f\M_0, \M_0} \cr
}}

Next suppose $\f_t\MM_{\lf(s)}$ is defined for all $t\in\Dom(s)$.
For a given $\lf(s)$ array 
    $\MM_{\lf(s)}=\MM_{[0,\lf(s))}$ 
set
\dispM{
    \f_s(\MM_{\lf(s)}) = \cases{
    \f_{s^+}(\MM_{\lf(s^+)})\wcat\MM_{[\lf(s^+),\lf(s))}
    & if $s^+\in\SS$, \cr
    \f_{s^-}(\MM_{\lf(s)})\wcat\MM_{\lf(s^-)}
    & else (i.e., $s^+\not\in\SS$). \cr
    }
}

Tail-operation shall take priority over concatenation.

Note that the definition is consistent with the fact that
\dispM{
    \MM_{\lf(s)}=\cases{
    \MM_{\lf(s^+)}\wcat\MM_{[\lf(s^+),\lf(s))}
    = \MM_{[0,\lf(s^+))}\wcat\MM_{[\lf(s^+),\lf(s))} 
    & if $s^+\in\SS$ \cr
    \MM_{\lf(s^-)} 
    = \MM_{[0,\lf(s^-))}
    & else. \cr
    }
}

The index '$s$' on the $\f$ operator is necessary in order
to disambiguate definiens $\f_{s^-}\MM_{\lf(s^-)}$ from 
definiendum $\f_{s}\MM_{\lf(s)}$ when $s^+\not\in\SS$, as 
in this case $\MM_{\lf(s)}=\MM_{\lf(s^-)}$. (E.g., we may 
have 
    $\MM_{\lf(\tup{0})} =\arr{\M} = \MM_{\lf(\tup{0,1})}$,
but the array of $\omega$ blocks for 
    $\f\A(\tup{0,1})\MM_{\lf(\tup{0,1})} = \f\M+\M\N$
contains two elements ($\f\M$ and $\M$), while
that for
    $\f\A(\tup{0})\MM_{\lf(\tup{0})} = \f\M$
contains only one ($\f\M$), i.e., $\f$-operations on
$\arr{M}$ depend on the sequence for which $\arr{M}$ is
the array.)

\parSkip
\Title {Example:}
Let's calculate 
    $\f A$
for
    $A=\A(\tup{0,0,1,2})[\M_0,\M_1,\M_2,\M_3,\M_4]$.
Note that 
\dispM{
    \eqalign{
    \lf(\tup{0,0,1,2}) 
    &= \lf(\tup{0,0,2})+\lf(\tup{0,0,1})\cr
    &= \lf(\tup{0,0})+\lf(\tup{0,1})+\lf(\tup{0,0})\cr
    &= \lf(\tup{0})+\lf(\tup{0})+\lf(\tup{0})+\lf(\tup{0})+\lf(\tup{0})
    = 5\cr
    }
}
and
\dispM{
    \eqalign{
    \A(\tup{0,0,1,2})[\M_0,\M_1,\M_2,\M_3,\M_4]
    &= \A(\tup{0,0,2})[\M_0,\M_1]+\A(\tup{0,0,1})[\M_2,\M_3,\M_4].\N\cr
    &= \A(\tup{0,0})[\M_0,\M_1].\Z+
      (\A(\tup{0,1})[\M_2]+\A(\tup{0,0})[\M_3,\M_4].\N).\N\cr
    &= (\M_0+\M_1\N)\Z + (\M_2\Z + (\M_3+\M_4\N)\N)\N.\cr
    }
}
But 
    $\Tree(\tup{0,0,1,2})$
is:
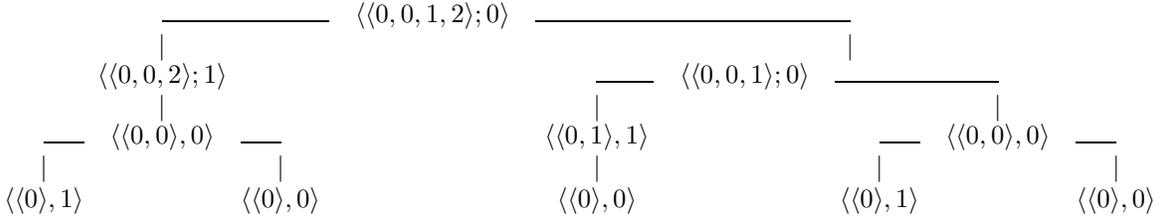
\dispM{
\mrx{ \tup{\tup{0,0,1,2};0} }
    { \nrx{ \tup{\tup{0,0,2};1} }
          { \mrx{ \tup{\tup{0,0},0} }
                { \tup{\tup{0},1} }
                { \tup{\tup{0},0} } } }
    { \mrx{ \tup{\tup{0,0,1};0} }
          { \nrx{ \tup{\tup{0,1},1} }
                { \tup{\tup{0},0} } }
          { \mrx{ \tup{\tup{0,0},0} }
                { \tup{\tup{0},1} }
                { \tup{\tup{0},0} } } }
}
Thus $\wt^s(0)=\wt^s(3)=1$ and $\wt^s(1)=\wt^s(2)=\wt^s(4)=0$
for $s=\tup{0,0,1,2}$ so assuming $A$ satisfies the uniformity
conditions, 
\dispM{
    \f\M_i\isom\f\M_j \hbox{ for all $i,j<5$, } \quad
    \M_0\isom\M_3,\quad
    \M_1\isom\M_2\isom\M_4
}
A cut-off point $x$ in $A$ may be chosen in five locations 
generically described by
\item{(1)} 
   $x=\tup{a,k}\in\M_0\times{k}$ for $k\in\Z$ 
from an $\omega$ sub-component of form $\M_0$ of the first 
$\omega+2$ component of $A$ (of form $(\M_0+\M_1\N)\Z$).
\item{(2)} 
    $x=\tup{\tup{a,n},k}\in\M_1\times{n,k}$ for $n\in\N$
and $k\in\Z$, from an $\omega$ sub-component of form 
$\M_1\times{n}$ of the first $\omega+2$ component of $A$. 
\item{(3)} 
    $x=\tup{\tup{a,k},n}\in\M_2\times{k,n}$ for $n\in\N$
and $k\in\Z$, from an $\omega$ sub-component of form 
$\M_2\times{k}$ of a successor $\omega+2$ component of $A$
(of form $(\M_2\Z + (\M_3+\M_4\N)\N)\times{n}$). 
\item{(4)} 
    $x=\tup{\tup{a,m},n}\in\M_3\times{m,n}$ for $m,n\in\N$,
from a an $\omega$ sub-component of form 
$\M_3\times{m}$ of a successor $\omega+2$ component of $A$.
\item{(5)} 
    $x=\tup{\tup{\tup{a,l},m},n}\in\M_4\times{l,m,n}$ for 
$l,m,n\in\N$, from an $\omega$ sub-component of form 
$\M_4\times{l}$ of a successor $\omega+2$ component of $A$.

Correspondingly,
\dispM{
    \f A = \cases{
    (\f\M_0+\M_1\N)+(\M_0+\M_1\N)\N+(\M_2\Z + (\M_3+\M_4\N)\N)\N
    & cut-off (1) \cr
    (\f\M_1+\M_1\N)+(\M_0+\M_1\N)\N+(\M_2\Z + (\M_3+\M_4\N)\N)\N
    & cut-off (2) \cr
    (\f\M_2+\M_2\N)+(\M_3+\M_4\N)\N+(\M_2\Z + (\M_3+\M_4\N)\N)\N
    & cut-off (3) \cr
    (\f\M_3+\M_4\N)+(\M_3+\M_4\N)\N+(\M_2\Z + (\M_3+\M_4\N)\N)\N
    & cut-off (4) \cr
    (\f\M_4+\M_4\N)+(\M_3+\M_4\N)\N+(\M_2\Z + (\M_3+\M_4\N)\N)\N
    & cut-off (5) \cr
    }
}
which are pairwise isomorphic under the uniformity conditions.


\parSkip
\Number{3.13}
\Title{Lemma}
\Subtitle{Tail uniformity}
If
    $\MM_{\lf(s)}=[\M_0,\ldots,\M_{\lf(s)-1}]$ 
satisfies the conditions of the uniformity over limit rank,
then
    $\f_{s^+}\MM_{[0,\lf(s^+))}
    = \f_{s^-}\MM_{[\lf(s^-),\lf(s))}$.

\Proof
By definition
\dispM{
    \f_{s^+}\MM_{[0,\lf(s^+))} = \cases {
    \f_{s^{++}}(\MM_{[0,\lf(s^{++}))})\wcat\MM_{[\lf(s^{++}),\lf(s^+))}
    & if $s^{++}\in\SS$, \cr
    \f_{s^{+-}}(\MM_{[0,\lf(s^{+-}))})\wcat\MM_{[0,\lf(s^{+-}))}
    & else (i.e., if $s^{++}\not\in\SS$); \cr
    }
}
while
\dispM{
    \f_{s^-}(\MM_{[\lf(s^+),\lf(s))}) = \cases {
    \f_{s^{-+}}(\MM_{[\lf(s^+),\lf(s^+)+\lf(s^{-+}))})
    & if $s^{-+}\in\SS$ \cr
    \f_{s^{--}}(\MM_{[\lf(s^+),\lf(s))})\wcat\MM_{[\lf(s^+),\lf(s))}
        \wcat\MM_{[\lf(s^+)+\lf(s^{-+}),\lf(s))}
    & else (i.e., if $s^{-+}\not\in\SS$), \cr
    }
}
the latter by applying the definition to the array 
    $\MM_{[0,\lf(t))}=\MM_{[\lf(s^+),\lf(s))}$ 
for $t=s^-$.

We shall be tacitly using the fact that $s^{+-}=s^{--}$, and hence
$\lf(s^{+-}) =\lf(s^{--})$ and $\f_{s^{+-}}=\f_{s^{--}}$.
If ${\bf A}=\tup{A_1,\ldots,A_n}$ and ${\bf B}=\tup{B_1,\ldots,B_n}$ 
are arrays of chains, we write ${\bf A}\isom{\bf B}$ to express that
$A_i\isom B_i$ for $i=1,\ldots,n$.

If the uniformity over limit rank conditions apply to $\MM_{\lf(s)}$, 
then the following hold:
\dispM{
    \eqalign{
    \hbox{Case (a)}\quad
    s^{++}\in\SS, s^{-+}\in\SS
    &\; : \;
    \MM_{[\lf(s^{++}),\lf(s^+))}
    \isom \MM_{[\lf(s^+)+\lf(s^{-+}),\lf(s))} \cr
%
%
    \hbox{Case (b)}\quad
    s^{++}\not\in\SS, s^{-+}\not\in\SS
    &\; : \;
    \MM_{[0,\lf(s^{+-}))}
    \isom \MM_{[\lf(s^+),\lf(s^+)+\lf(s^-))} \cr
    \hbox{Case (c)}\quad
    s^{++}\in\SS,s^{-+}\not\in\SS
    &\; : \;
    \MM_{[\lf(s^{++}),\lf(s^+))}
    \isom \MM_{[\lf(s^+),\lf(s))} \cr
    \hbox{Case (d)}\quad
    s^{++}\not\in\SS,s^{-+}\in\SS
    &\; : \;
    \MM_{[0,\lf(s^{+-}))}
    \isom \MM_{[\lf(s^+)+\lf(s^{-+}),\lf(s))} \cr
    }
}

In case (a), note that by leaf-length arithmetic,
\dispM{
    \eqalign{
    &\MM_{[\lf(s^{++}),\lf(s^+))}
    = \MM_{[\lf(s^{++}),\lf(s^{++})+\lf(s^{+-}))} \cr
    &\MM_{[\lf(s^+)+\lf(s^{-+}),\lf(s))}
    = \MM_{[\lf(s^+)+\lf(s^{-+}),
        \lf(s^+)+\lf(s^{-+})+\lf(s^{--}))}. \cr
    }
}

\def\embeds{\hookrightarrow}

This follows from the partial 'tree isomorphisms' between the
subtree of $\Tree(s)$ extending $\tup{s^{+-},0}$ and the one
extending $\tup{s^{--},0}$, which in turn is immediate form
$s^{+-}=s^{--}$. (The aforemetioned subtrees are both isomorphic
to the tree $\Tree(t)$, for $t=s^{+-}=s^{--}$ if we consider 
trees $\Tree(s)$ as partial orders over a set of nodes 
$\tup{t,k}\in\Dom(s) \times\mod{s}$.) Call these subtrees
$\Tree^{+-}=\Ext(\Tree(s,\tup{s^{+-},0})$ and
$\Tree^{--}=\Ext(\Tree(s,\tup{s^{--},0})$, respectively.
Let temporarily '$\embeds$' denote the embeddings of
$\Tree^{+-}$ and $\Tree^{--}$ into $\Tree=\Tree(s)$, and let
$[a,b)$ denote the 'interval' between leaves $a$ (inclusive)
and $b$ (exclusive) of $\Lev_{\mod{s}}(\Tree)$. Then
\dispM{
    \eqalign{
    \hbox{[a]}\quad
    s^{++}\in\SS, s^{-+}\in\SS
    \; : \;
    \Lev_{\mod{s^{+-}}}(\Tree^{+-})
    & \embeds [\leaf(\lf(s^{++})),\leaf(\lf(s^{++})+\lf(s^{+-}))), \cr
    & = [\leaf(\lf(s^{++})),\leaf(\lf(s^+))) \cr
    \Lev_{\mod{s^{--}}}(\Tree^{--})
    & \embeds [\leaf(\lf(s^+)+\lf(s^{-+}))),
        \leaf(\lf(s^+)+\lf(s^{-+})+\lf(s^{--}))), \cr
    & = [\leaf(\lf(s^+)+\lf(s^{-+})),\leaf(\lf(s))) \cr
    \hbox{[b]}\quad
    s^{++}\not\in\SS, s^{-+}\not\in\SS
    \; : \;
    \Lev_{\mod{s^{+-}}}(\Tree^{+-})
    & \embeds [\leaf(0),\leaf(\lf(s^{+-}))), \cr
    \Lev_{\mod{s^{--}}}(\Tree^{--})
    & \embeds [\leaf(\lf(s^+)),\leaf(\lf(s^+)+\lf(s^-))) \cr
    \hbox{[c]}\quad
    s^{++}\in\SS, s^{-+}\not\in\SS
    \; : \;
    \Lev_{\mod{s^{+-}}}(\Tree^{+-})
    & \embeds [\leaf(\lf(s^{++})),\leaf(\lf(s^+))), \cr
    \Lev_{\mod{s^{--}}}(\Tree^{--})
    & \embeds [\leaf(\lf(s^+)),\leaf(\lf(s))) \cr
    \hbox{[d]}\quad
    s^{++}\not\in\SS, s^{-+}\in\SS
    \; : \;
    \Lev_{\mod{s^{+-}}}(\Tree^{+-})
    & \embeds [\leaf(0),\leaf(\lf(s^{+-}))), \cr
    \Lev_{\mod{s^{--}}}(\Tree^{--})
    & \embeds [\leaf(\lf(s^+)+\lf(s^{-+})),\leaf(\lf(s))) \cr
    }
}

Schematically, where the underlines indicate isomorphisms between
leaf-intervals:
\item{}
Case $s^{++},s^{-+}\in\SS$:
\dispM{
\mrx{ \tup{s;0} }
    { \mrx{ \tup{s^+;1} }
          { \tup{s^{++};2} }
          { \underline{\Tree(s^{+-})} } }
    { \mrx{ \tup{s^-;0} }
          { \tup{s^{-+};1} }
          { \underline{\Tree(s^{--})} } }
}
\item{}
Case $s^{++},s^{-+}\not\in\SS$: ($s^{+-}=s^{--}$)
\dispM{
\mrx{ \tup{s;0} }
    { \nrx{ \tup{s^+;1} }
          { \underline{\Tree(s^{--})} } }
    { \nrx{ \tup{s^-;0} }
          { \underline{\Tree(s^{--})} } }
}
\item{}
Case $s^{++}\in\SS,s^{-+}\not\in\SS$:
\dispM{
\mrx{ \tup{s;0} }
    { \mrx{ \tup{s^+;1} }
          { \tup{s^{++};2} }
          { \underline{\Tree(s^{+-})} } }
    { \nrx{ \tup{s^-;0} }
          { \underline{\Tree(s^{--})} } }
}
\item{}
Case $s^{++}\not\in\SS,s^{-+}\in\SS$:
\dispM{
\mrx{ \tup{s;0} }
    { \nrx{ \tup{s^+;1} }
          { \underline{\Tree(s^{+-})} } }
    { \mrx{ \tup{s^-;0} }
          { \tup{s^{-+};1} }
          { \underline{\Tree(s^{--})} } }
}

Using the inductive hypothesis and $s^{+-}=s^{--}$,
\dispM{
    \eqalign{
    \f_{s^{++}}\MM_{[0,\lf(s^{++}))}
    &\isom \f_{s^{+-}}\MM_{[\lf(s^{++}),\lf(s^{+}))} \cr
    &= \f_{s^{+-}}\MM_{[\lf(s^{++}),\lf(s^{++})+\lf(s^{+-}))} \cr
    &= \f_{s^{--}}\MM_{[\lf(s^{++}),\lf(s^{++})+\lf(s^{--}))}, \cr
    \f_{s^{-+}}\MM_{[0,\lf(s^{-+}))}
    &\isom \f_{s^{--}}\MM_{[\lf(s^{-+}),\lf(s^{-}))} \cr
    &= \f_{s^{--}}\MM_{[\lf(s^{-+}),\lf(s^{-+})+\lf(s^{--}))}, \cr
    \f_{s^{-+}}\MM_{[\lf(s^+),\lf(s^+)+\lf(s^{-+}))}
    &\isom \f_{s^{--}}\MM_{[\lf(s^+)+\lf(s^{-+}),
                            \lf(s^+)+\lf(s^{-}))} \cr
    &= \f_{s^{--}}\MM_{[\lf(s^+)+\lf(s^{-+}), \lf(s))}. \cr
    }
}
So,
\dispM{
    \eqalign{
        \f_{s^+}\MM_{[0,\lf(s^+))}
        &= \f_{s^{--}}\MM_{[\lf(s^{++}),\lf(s^+))}
           \wcat\MM_{[\lf(s^{++}),\lf(s^+))}, \cr
        \f_{s^-}\MM_{[\lf(s^+),\lf(s))}
        &= \f_{s^{--}}\MM_{[\lf(s^+)+\lf(s^{-+}),\lf(s))}
           \wcat\MM_{[\lf(s^+)+\lf(s^{-+}),\lf(s))}. \cr
    }
}

But by (a), the 'second halves' are isomorphic
    $\MM_{[\lf(s^{++}),\lf(s^+))}
    \isom \MM_{[\lf(s^+)+\lf(s^{-+})+\lf(s))}$.
In the immediately preceding equations, the respective 'first halves' 
of $\f_{s^+}\MM_{[0,\lf(s^+))}$ and $\f_{s^-}\MM_{[\lf(s^+),\lf(s))}$ 
are isomorphic to tails of the 'second' ones, hence also isomorphic
to each other, hence so are the respective left-hand terms, as 
required.

Case $s^{++},s^{-+}\not\in\SS$: In that case,
\dispM{
    \eqalign{
    \f_{s^+}\MM_{[0,\lf(s^+))}
    &= \f_{s^{+-}}\MM_{[0,\lf(s^{+-}))}\wcat\MM_{[0,\lf(s^{+-}))} \cr
    \f_{s^-}\MM_{[\lf(s^+),\lf(s))}
    &= \f_{s^{--}}\MM_{[\lf(s^+),\lf(s^+)+\lf(s^-))}
        \wcat\MM_{[\lf(s^+),\lf(s^+)+\lf(s^-))} \cr
    &= \f_{s^{--}}\MM_{[\lf(s^+),\lf(s))}
        \wcat\MM_{[\lf(s^+),\lf(s))} \cr
    }
}
Now by (b), 
    $\MM_{[0,\lf(s^{+-}))}
    \isom \MM_{[\lf(s^+),\lf(s^+)+\lf(s^-))}
    = \MM_{[\lf(s^+),\lf(s))}$.
it follows immediately that the respective first and second
'halves' are isomorphic, whence so are the respective left-hand
terms.

Case $s^{++}\in\SS,s^{-+}\not\in\SS$:
Using the inductive hypothesis ans $s^{+-}=s^{--}$,
\dispM{
    \eqalign{
    \f_{s^+}\MM_{[0,\lf(s^+))}
    &= \f_{s^{++}}\MM_{[0,\lf(s^{++}))}
        \wcat\MM_{[\lf(s^{++}),\lf(s^+))} \cr
    &= \f_{s^{+-}}\MM_{[\lf(s^{++}),\lf(s^+)}
        \wcat\MM_{[\lf(s^{++}),\lf(s^+))} \cr
    &= \f_{s^{--}}\MM_{[\lf(s^{++}),\lf(s^+)}
        \wcat\MM_{[\lf(s^{++}),\lf(s^+))} \cr
    \f_{s^-}\MM_{[\lf(s^+),\lf(s))}
    &= \f_{s^{--}}\MM_{[\lf(s^+),\lf(s^+)+\lf(s^-))}
        \wcat\MM_{[\lf(s^+),\lf(s^+)+\lf(s^-))} \cr
    &= \f_{s^{--}}\MM_{[\lf(s^+),\lf(s))}
        \wcat\MM_{[\lf(s^+),\lf(s))} \cr
    }
}
By (c), 
    $\MM_{[\lf(s^{++}),\lf(s^+)} 
    \isom \MM_{[\lf(s^+),\lf(s))}$,
it follows again that the respective first and second
'halves' are isomorphic, whence so are the respective left-hand
terms.

Case $s^{++}\not\in\SS,s^{-+}\in\SS$:
From $s^{+-}=s^{--}$ and the inductive hypothesis,
\dispM{
    \eqalign{
    \f_{s^+}\MM_{[0,\lf(s^+))}
    &= \f_{s^{+-}}\MM_{[0,\lf(s^{+-}))}
        \wcat\MM_{[0,\lf(s^{+-}))} \cr
    &= \f_{s^{--}}\MM_{[0,\lf(s^{+-}))}
        \wcat\MM_{[0,\lf(s^{+-}))} \cr
    \f_{s^-}\MM_{[\lf(s^+),\lf(s))}
    &= \f_{s^{-+}}\MM_{[\lf(s^+),\lf(s^+)+\lf(s^{-+}))}
        \wcat\MM_{[\lf(s^+)+\lf(s^{-+}),\lf(s))} \cr
    &= \f_{s^{--}}\MM_{[\lf(s^+)+\lf(s^{-+}),\lf(s))}
        \wcat\MM_{[\lf(s^+)+\lf(s^{-+}),\lf(s))}, \cr
    }
}
whence (d), 
    $\MM_{[0,\lf(s^{+-}))} 
    \isom \MM_{[\lf(s^+)+\lf(s^{-+},\lf(s))}$
implies that the respective segments, and hence the respective
left-hand terms are isomorphic.
\Qed

As a corollary, we get

\parSkip
\Number{3.14}
\Title{Lemma}
\Subtitle{Uniformity over limit lemma - sufficiency part}
If
    $\MM_{\lf(s)}=[\M_0,\ldots,\M_{\lf(s)-1}]$ 
satisfies the conditions of the uniformity over limit rank,
then
    $A = \A(s)[\M_0,\ldots,\M_{\lf(s)-1}]$ 
is uniform and in fact,
\dispM{
    \f(\A(s)[\MM_{\lf(s)}]) = \A(\f s)[\f_s(\MM_{\lf(s)})].
}

\Proof
Base-cases:
The condition is vacuous for $s=\tup{0}$.
Cleary, $\A(\tup{0,0})[\M_0,\M_1]\in\U(\omega+1)$ just if
$\f\M_0\isom\f\M_1$ (there are no further conditions on
$[\M_0,\M_1]$ as $\wt^{\tup{0,0}}(0)\not=\wt^{\tup{0,0}}(1)$),
$\f_{\tup{0,0}}[\M_0,\M_1] = [\f\M_0,\M_1]$, and
\dispM{
    \f(\A(\tup{0,0})[\M_0,\M_1])
    = \f(\M_0+\M_1\N)
    = \M+\M_1\N
    = \A(\tup{0,0})[\M,\M_1]
}
for $\M\isom\f\M_0\isom\f\M_1$.
For $\A(\tup{0,1})[\M_0]\in\U(\omega+1)$ there are no further 
conditions on $[\M_0]$ (except the underlying that $\M_0\in
\UC(\omega)$), $\f_{\tup{0,1}}[[\M_0] = [\f\M_0,\M_0]$, and
\dispM{
    \f(\A(\tup{0,1})[\M_0])
    = \f(\M_0\Z)
    = \M+\M_0\N
    = \A(\tup{0,0})[\M,\M_0]
    = \A(\f\tup{0,1})[\M,\M_0]
}
for $\M\isom\f\M_0\isom\f\M_1$.

Let now $s\in\SS$ and assume the lemma is true for all $t\in\Dom(s)$.

Case $s^+\not\in\SS$: Here,
\dispM{
    \eqalign{
    \f(\A(s)[\MM_{\lf(s)}])
    &= \f(\A(s^-)[\MM_{[\lf(s^+),\lf(s))}].\Z) \cr
    &= \f(\A(s^-)[\MM_{\lf(s^-)}])
     + \A(s^-)[\MM_{[\lf(s^+),\lf(s))}].\N\cr
    &= \A(\f s^-)[\f_{s^-}(\MM_{\lf(s^-))}]
     + \A(s^-)[\MM_{[\lf(s^+),\lf(s))}].\N\cr
    &= \A(\f s^)[\f_{s^-}(\MM_{\lf(s^-))}\wcat
     + [\MM_{\lf(s^-)}].\N\cr
    &= \A(\f s)[\f_s\MM_{\lf(s)}],\cr 
    }
}
using the inductive hypothesis in step 3, the definition of
$\A(t)[\MM_{\lf(t)}]$ for $t=\f s$ in step 4, and the definition
of $\f_s(\MM_{\lf(s)})$ in step 5.

Case $s^+\in\SS$: Here,
\dispM{
    \eqalign{
    \f(\A(s)[\MM_{\lf(s)}])
    &= \f(\A(s^+)[\MM_{\lf(s^+)}]
     + \A(s^-)[\MM_{[\lf(s^+),\lf(s))}].\N) \cr
    &= \cases{
    \f(\A(s^+)[\MM_{\lf(s^+)}])+\A(s^-)[\MM_{[\lf(s^+),\lf(s))}].\N 
    & cut-off in 
        $\A(s^+)[\MM_{\lf(s^+)}]$ \cr
    \f(\A(s^-)[\MM_{\lf(s^-)}])+\A(s^-)[\MM_{[\lf(s^+),\lf(s))}].\N 
    & cut-off in some 
        $\A(s^-)[\MM_{[\lf(s^-),\lf(s)}]\times\ens{n}$. \cr } \cr
    }
}
By inductive hypothesis,
\dispM{
    \eqalign{
    \f(\A(s^+)[\MM_{\lf(s^+)}])
    &= \A(\f s^+)[\f_{s^+}(\MM_{\lf(s^+)}], \cr
    \f(\A(s^-)[\MM_{[\lf(s^+),\lf(s))}])
    &= \A(\f s^-)[\f_{s^-}(\MM_{[\lf(s^+),\lf(s))}]. \cr
    }
}
But $\f s^+=\f s^-=s^=\tac{0}$ by the tail string lemma, and
    $\f_{s^+}(\MM_{[0,\lf(s^+))}) 
    = \f_{s^-}(\MM_{[\lf(s^+),\lf(s))})$
by the tail uniformity lemma.
\Qed

\parSkip
\Number{3.15}
\Title{Remark}
The following description yields an equivalent definition of $\f_s$. 
The base case is as in definition 3.12. For the inductive step, 
suppose 
    $\f_t(\MM_{\lf(t)})$ is defined for every $t\in\Dom(s)$.
Let $n(s)=\mod{s}-1$ and $k(s)=\max\ens{i\le n(s): s(i)=i}$

\item{(a)}
Case $k(s)=n(s)$. Then $\lf(s)=\lf(s^-)$, while $\lf(\f s)=
\lf((\f s)^+)+\lf((\f s)^-)=\lf(\f s^-)+\lf(s^-)$. 
By hypothesis, $\f(\MM_{\lf(s^-)})$ is defined. 
So
\dispM{
    \f_s(\MM_{\lf(s)})
    = \f_{s^-}(\MM_{\lf(s^-)})\wcat\MM_{\lf(s)}
}

\item{(b)}
Case $k(s)<n(s)$. Then $\lf(s)=\lf(s^+)+\lf(s^-)$, and 
$\lf(\f s)=\lf((\f s)^+)+\lf((\f s)^-)=\lf(\f s^-)+\lf(s^-)$. 
Let $k$ maximal $\le n(s)$ such that $s^{-k}\in\SS$.
By induction
\dispM{
    \lf(s)=\lf(s^{-k})+\lf(s^{-(k-1)})+\ldots+\lf(s^-),
}
and by hypothesis, $\f_{s^{-i}}(\MM_{\lf(s^{-i})})$ are defined
for $i=1,\ldots,k$ and $\f_{s^{-k}}(\MM_{\lf(s^{-k})})$ is defined
by case (a). 
So
\dispM{
    \f_s(\MM_{\lf(s)})
    = \f_{s^{-k}}(\MM_{\lf(s^{-k})})
    \wcat\MM_{[\lf(s^{-(k-1)}),\lf(s^{-(k-2)}))}
    \wcat\ldots
    \wcat\MM_{[\lf(s^{-2}),\lf(s^{-}))}
}

\secSkip
\Section{[3.3]}
\Title{Short representations}

\parSkip\noindent
Set 
\dispm{
    \rwt^s = \ran(\wt^s). 
}
Given $A=\A(s)\arr{\M_i}_{i<\lf^s}$, by the $\omega+$finite 
uniformity lemma, there are chains
    $\U_k\in\UC(\omega)$, $k\in\rwt^s$
such that
\dispM{
    (\forall i\in\lf(s))[\M_i=\U_{\wt^s(i)}].
    \eqno{(3.2)}
}
(i.e., the set 
    $\ens{\U_k:k\in\rwt^s}$
contains just one representative (suitably tagged) under 
isomorphism, for each of the possibly non-isomorphic $\omega$ 
components.) Set
    $\A(s)\rarr{\tup{\U_k,k}}_{k\in\rwt^s}
    := \A(s)\arr{\M_i}_{i<\lf(s)}$
whenever (6.3) holds and refer to this as the {\it tagged short
representation} of $A=\A(s)\arr{M_i}_{i<\lf(s)}$. Also refer to
    $\UU = \ens{\tup{\U_k,k}}_{k\in\rwts^s}$ 
as the {\it tagged short array} of that representation and to 
    $\UU$
as the {\it tagged reduction} of 
    $\MM = \tup{\M_i}_{i<\lf(s)}$. 

\parSkip
\Title{Corollary}
For any $s\in\SS(n+1)$ and any 
    $\ens{\U_k:k\in\rwt^s}$, $\U_k\in\UC(\omega)$ 
such that 
\dispM{
    \hbox{ $\f\U_k\isom\f\U_l$ for any $k,l\in\rwt^s$, }
    \eqno{(3.3)}
}
    $\A(s)\rarr{\tup{\U_k,k}:k\in\rwt^s}$ 
is a uniform chain of rank $\omega+n$ and every
$A\in\UC(\omega,\omega+\omega)$ has a repsentation of
this form.

\parSkip
Whenever referring to a set
    $\ens{\tup{\U_k,k}:k\in\rwt^s}$, $\U_k\in\UC(\omega)$ 
as a tagged short array, condition (3.3) is assumed to
hold.

\parSkip
\Title{Example}
Consider $s=\tup{0,0,1,0}$: $\Tree(s)$, fully expanded, is
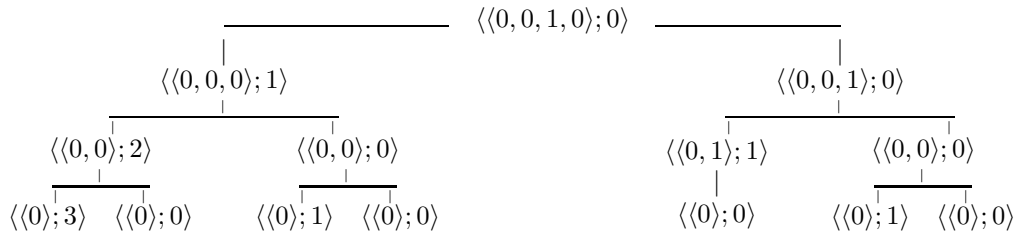
\dispM{
\mrx{ \tup{\tup{0,0,1,0};0} }
    { \prx{ \tup{\tup{0,0,0};1} }
          { \qrx{ \tup{\tup{0,0};2} }
                { \tup{\tup{0};3} }
                { \tup{\tup{0};0} } }
          { \qrx{ \tup{\tup{0,0};0} }
                { \tup{\tup{0};1} }
                { \tup{\tup{0};0} } } }
    { \prx{ \tup{ \tup{0,0,1};0} }
          { \nrx{ \tup{\tup{0,1};1} }
                { \tup{\tup{0};0} } }
          { \qrx{ \tup{\tup{0,0};0} }
                { \tup{\tup{0};1} }
                { \tup{\tup{0};0} } } }
}
so $\wt^s=\tup{3,0,1,0,0,1,0}$ and with a tagged short array 
    $\ens{\tup{\U_0,0},\tup{\U_1,1},\tup{\U_3,3}}$, 
    $\U_0,\U_1,\U_3\in\UC(\omega)$ 
\dispM{
    \A(\tup{0,0,1,0})\rarr{\tup{\U_0,0},\tup{\U_1,1},\tup{\U_3,3}} 
    = \A(\tup{0,0,1,0})\arr{\U_3,\U_0,\U_1,\U_0,\U_0,\U_1,\U_0} 
}

\def\SubLfs{{\cal J}}
 
\parSkip
If it seems that the price in tags seems high compared to the 
apparent gain in length-reduction, observe for example that 
for all
    $s_n = \tup{0}\tac{1^{(n+1)}}$,
$\lf(s_n) = 2^n$, while $\sz{\rwts^{s_n}} = 1$. This will be 
become relevant in situations where sequences $s_n$ arise 
naturally (such as when approximating $A\in\UC(\omega+\omega)$).

Moreover, the tagging can be made dispensible in turn, by
defining a canonical ordering. There is an obvious association
between the nodes of form $\tup{t,k}$ of $\Tree^s$ and components
of 
    $A = \A\arr{\M_i}_{i<\lf(s)}$ 
of form 
    $\A(t)\arr{\M_i}_{i\in\SubLfs}$.
(Here $\SubLfs$ is a sub-interval of $[0,\lf(s))$, i.e.,
    $\tup{\M_i}_{i\in\SubLfs}$ 
is the set of building-blocks
    $\tup{\M_i'}_{i<\lf(t)}$ 
used for that component). Namely, starting with $A$ at the root
and proceeding dow the tree, associate the 'plus'/'minus' 
subcomponent(s) arived at to the 'plus'/'minus' extension(s) of 
the node arrived at. Next, extend that association to the tagged
short arrays of these components. It then turns out that along 
any branch 
    $\sigma_0,\ldots,\sigma_j,\ldots,\sigma_n
    =\tup{\tup{0},k_0},\ldots,\tup{t,k},\ldots,\tup{s,0}$
of $\Tree^s$, proceeding from the leaves to the root, the arrays 
$\UU_{t,k}$ associated to the nodes $\tup{t,k}$ in that way form 
a chain whose members increase in size by at most 1 over their
predecessors. That generates an ordering of the elements in the 
union over such a chain (e.g., the order of 'aquisition'). But that 
union is always $\UU=\UU_{s,0}$, the array for $A$ itself (that gets 
associated with the top node). For a 'canonical' ordering it then 
suffices that we can always single out a unique such branch, say by 
chosing the right-most one (i.e., following consistently the
'minus' extensions as considered from the root).

\parSkip
\Title{Example}
Here is the 'associated' tree of tagged short arrays for
$t\in\Dom(s)$ for $s=\tup{0,0,1,0}$ of the preceding example
 -- note that the tags on the chains in $\UU_{t,k}$ now 
represent the weights from the weight-sequence of sequence 
$t$ of the corresponding node $\tup{t,k}$, not the weight $k$ 
of the node:
\dispM{
\orx{ \ens{\tup{\U_0,0},\tup{\U_1,1},\tup{\U_3,3}} }
    { \prx{ \ens{\tup{\U_0,0},\tup{\U_1,1},\tup{\U_3,2}} }
          { \qrx{ \ens{\tup{\U_0,0},\tup{\U_3,1}} }
                { \ens{\tup{\U_3,0}} }
                { \ens{\tup{\U_0,0}} } }
          { \qrx{ \ens{\tup{\U_0,0},\tup{\U_1,1}} }
                { \ens{\tup{\U_1,0}} }
                { \ens{\tup{\U_0,0}} } } }
    { \prx{ \ens{\tup{\U_0,0},\tup{\U_1,1}} }
          { \nrx{ \ens{\tup{\U_0,0}} }
                { \ens{\tup{\U_0,0}} } }
          { \qrx{ \ens{\tup{\U_0,0},\tup{\U_1,1}} }
                { \ens{\tup{\U_1,0}} }
                { \ens{\tup{\U_0,0}} } } }
}
The 'canonical' ordering of the chains in the tagged short
array would thus be
    $\tup{\U_0,\U_1,\U_3}$,
corresponding to the right-most branch, starting at leaf $i=6$. 
(E.g., using the left-most branch, starting at leaf $i=0$ would 
result in
    $\tup{\U_3,\U_0,\U_1}$
and starting at leaf $i=1$ would result in
    $\tup{\U_0,\U_3,\U_1}$.
Etc.) The canonically ordered array shall be the desired {\it short
array}.
%
%
Since the size of $\UU_{t,k}$ is equal to that of $\rwt^t$, 
some facts are required about how $\sz{\rwt^s}$ relates to 
$\sz{\rwt^{s^{\pm}}}$ in order to implement this.

\parSkip
\Number{(3.16)}
\Title{Definition}
\Subtitle{Leaf-intervals of components and component-node association}
Suppose 
    $A=\A\arr{\M_i}_{i<\lf(s)}\in\UC(\omega+n)$,
where 
    $n=\mod{s}-1$.
\item{(i)}
Set $\Lf_{s,0}=[0,\lf(s))$, the leaf interval of $A$.
Associate $A$ to $\tup{s,0}$.
\item{(ii)}
Let
    $\tup{t,k}\in\Tree^s$ 
and supose 
    $\Lf_{t,k}\subeq\Lf_{s,0}$ 
is the leaf-interval (of length $\lf(t)$) of 
    $A_{t,k}=\A(t)\arr{\M_i}_{i\in\Lf_{t,k}}$ 
and suppose $A_{t,k}$ is associated to that node.
Let 
    $m=\mod{t}-1$.
\itemitem{(a)}
If $t(m)=m$ then
    $\tup{t^-,0}$ 
is the only extension of $\tup{t,k}$ 
    $\lf(t)=\lf(t^-)$,
and 
    $\A(t)\arr{\M_i}_{i\in\Lf_{t,k}} 
    =\A(t^-)\arr{\M_i}_{i\in\Lf_{t,k}}.\Z$,
so set
    $\Lf_{t,k}^- := \Lf_{t^-,0} := \Lf_{t,k}$
and associate
    $\A(t^-)\arr{\M_i}_{i\in\Lf_{t,k}}$ 
to $\tup{t^-,0}$,
\itemitem{(b)}
Else ($t(m)<m$), there are two extensions
    $\tup{t^+,k+1}$ and $\tup{t^-,0}$
and
    $\lf(t)=\lf(t^+)+\lf(t^-)$.
Set
    $\Lf_{t,k}^+ := \Lf_{t^+, k+1} := $ 
    the initial segment of $\Lf_{t,k}$ of length $\lf(t^+)$,
and 
    $\Lf_{t,k}^- := \Lf_{t^-,0} := $
    the final segment of $\Lf_{t,k}$ of length $\lf(t^-)$. 
Then
    $\A(t)\arr{\M_i}_{i\in\Lf_{t,k}} 
    =\A(t^+)\arr{\M_i}_{i\in\Lf_{t,k}^+}
    +\A(t^-)\arr{\M_i}_{i\in\Lf_{t,k}^-}.\N$,
so associate 
    $\A(t^+)\arr{\M_i}_{i\in\Lf_{t,k}^+}$ to $\tup{t^+,k+1}$ 
and 
    $\A(t^-)\arr{\M_i}_{i\in\Lf_{t,k}^-}$ to $\tup{t^-,0}$.

\parSkip
Let
    $\MM=\tup{\M_i}_{i<\lf(s)}
    \in\Fun{\lf(s)}{\UC(\omega)}$, 
be the array in the representation of some 
    $A\in\UC(\omega)$. 
Implicit in the definition of $A$ as
    $\A(s)\arr{\MM}$,
is a mapping
\dispM{
    [0,\lf(s))\to
    \cases{
        [0,\lf(s))^- := [0,\lf(s^-)) = [0,\lf(s))
        & if $s(n)=n$ \cr
        [0,\lf(s))^+, [0,\lf(s))^- 
        := [0,\lf(s^+)), [\lf(s^+),\lf(s))       
        & if $s(n)<n$ \cr
    }
}
and hence
\dispM{
    \MM\to
    \cases{
        \MM^-       
        := \tup{\M_i}_{i<\lf(s^-)} = \MM
        & if $s(n)=n$ \cr
        \MM^+,\MM^- 
        := \tup{\M_i}_{i<\lf(s)}, 
          \tup{\M_i}_{\lf(s^+)\le i<\lf(s)} 
        & if $s(n)<n$. \cr
    }
}
Define the corresponding operations 
    $\UU\to\UU^{\pm}$ 
on short arrays so that
    $\UU^{\pm}$ is the reduction of $\MM^{\pm}$; 
again, these will depend on $s$, so we effectively have an operation
    $\UU\to\UU^{\pm}_s$.
(In case of the arrays $\MM=\tup{\M_i:i<\lf(s)}$ the dependency
on $s$ is implicit in the notation. The operation is assumed to be
defined if 
    $\UU$ is of form $\ens{\U_k}_{k\in\rwt^s}$
for chains $\U_k\in\UC(\omega)$). Thus: 

\parSkip
\Number{(3.17)}
\Title{Definition}
\Subtitle{Decomposition of tagged short arrays}
Where
    $\UU=\ens{\tup{\U_k,k}}_{k\in\rwt^s}$,
set and $n=\mod{s}-1$,
\dispM{
    \eqalign{
    \UU^-_s 
    &= \ens{\tup{\U_k,k}}_{k\in\rwt^{s^-}} \cr
    \UU^+_s 
    &= \cases{
         \ens{\tup{\U_n,n-1}}
         \cup\ens{\tup{\U_k,k}}_{k\in\rwt^{s^+}\rid\ens{n-1}} 
         & if $s(n)=0$ \cr
         \ens{\tup{\U_k,k}}_{k\in\rwt^{s^+}} 
         & if $0<s(n)<n$ \cr
         \hbox{undefined}
         & if $s(n) = n$. \cr 
         } \cr 
    }
}

\parSkip
\noindent
(The change from 
    $\tup{U_n, n}\in\UU$ 
to 
    $\tup{U_n, n-1}\in\UU^+_s$
in case $s(n)=0$ is due to the fact that
    $\wt^s(0) = \wt^{s^+}(0)+1 = n$ just if $s(n)=0$ for $n=\mod{s}-1$,
by the second weight induction corollary).
This allows us to define the chain constructor in terms of the 
short arrays and obtain the 'transfinite' analogue of Pastjin's
uniform chain constructor:


\parSkip
\Number{(3.18)}
\Title{Corollary}
\Subtitle{Plus finite rank Pastjin unifom chain constructor ---
short version}
Where
    $\UU = \ens{\tup{\U_k,k}}_{k\in\rwt^s}
    \subeq\UC(\omega)\times\rwt^s$
is a tagged short array,
\dispM{
    \A(s)\rarr{\UU}
    = \cases{
        \A(s^-)\rarr{\UU^-_s}.\Z
        & if $s(n)=n$ for $n=\mod{s}-1$ \cr
        \A(s^+)\rarr{\UU^+_s} + \A(s^-)\rarr{\UU^-_s}.N
        & if $s(n)<n$ for $n=\mod{s}-1$ \cr
    }
}

\parSkip
\Title{Example}
(Comparing the short arrays 
    $\U_0,\ldots,\U_{\kwt^s}$
with the full array,
    $\M_0,\ldots,\M_{\lfs-1}$.)


 
If $s=\tup{0,0,0}$, then $\tp_1\in\ens{0,1}$ and
$\new^s=\ens{0}$ if $\tp_1=1$ and $\new^s=\ens{2}$ if 
$\tp_0=0$, so
\dispM{
    \eqalign{
    &\A(\gtup{0,0,0})\rarr{\U_0,\U_1,\U_2})\cr
    &= \A(\gtup{0,0,0})
       \arr{\M_0,\M_1,\M_2,\M_3} \cr
    &= \A(\gtup{0,0})\arr{\M_0,\M_1}
        + \A(\gtup{0,0})\arr{\M_2,\M_3}.\N \cr
    &= \cases{
       \A(\gtup{0,0})\arr{\U_0,\U_1} 
        + \A(\gtup{0,0})\arr{U_2,\U_1}.\N 
       & if $\tp_0=0$ and $\tp_1=0$ 
           ($\new^{s^+}=\ens{1}$ and $\new^s=\ens{2}$) \cr
       \A(\gtup{0,0})\arr{\U_2,\U_1} 
        + \A(\gtup{0,0})\arr{U_0,\U_1}.\N 
       & if $\tp_0=0$ and $\tp_1=1$ 
           ($\new^{s^-}=\ens{1}$ and $\new^s=\ens{0}$) \cr
       \A(\gtup{0,0})\arr{\U_1,\U_0} 
        + \A(\gtup{0,0})\arr{U_2,\U_0}.\N 
       & if $\tp_0=1$ and $\tp_1=0$ 
           ($\new^{s^+}=\ens{0}$ and $\new^s=\ens{2}$) \cr
       \A(\gtup{0,0})\arr{\U_2,\U_0} 
        + \A(\gtup{0,0})\arr{U_1,\U_0}.\N 
       & if $\tp_0=1$ and $\tp_1=1$ 
           ($\new^{s^-}=\ens{0}$ and $\new^s=\ens{0}$) \cr} 
    \cr
    &= \cases{
       \A(\gtup{0,0})\rarr{\U_0,\U_1}
        + \A(\gtup{0,0};\rarr{U_2,\U_1}.\N 
       & \cr
       \A(\gtup{0,0})\rarr{\U_2,\U_1}
        + \A(\gtup{0,0})\rarr{U_0,\U_1}.\N 
       & \cr
       \A(\gtup{0,0})\rarr{\U_1,\U_0}
        + \A(\gtup{0,0})\rarr{U_2,\U_0}.\N 
       & \cr
       \A(\gtup{0,0})\rarr{\U_2,\U_0}
        + \A(\gtup{0,0})\rarr{U_1,\U_0}.\N 
       & \cr} 
    \cr
    &= \cases{
       \U_0+\U_1.\N+(\U_2+\U_1.\N).\N
       & \cr
       \U_2+\U_1.\N+(\U_0+\U_1.\N).\N
       & \cr
       \U_1+\U_0.\N+(\U_2+\U_0.\N).\N
       & \cr
       \U_2+\U_0.\N+(\U_1+\U_0.\N).\N
       & \cr}
    \cr
    \cr
    }
}
by following the 'line up' of the $\U_i$ along the
association tree for $\tup{0,0,0}$ (i.e., depending on where
we select the 'base-point'. E.g., if
       $\tp_0=1$ and $\tp_1=0$ 
we start with the second leaf $\leaf^s(1)=\tup{0}$ of $\Tree^s$
picking up $\U_0$, proceed to its parent node $\tup{0,0}$ on
the left branch of $s$ picking up $\U_1$ and arriving at the
top-node $s$, picking up $\U_2$.)

  
\chapSkip
\Chapter{[4] Products of uniform chains}

\secSkip\noindent
The preceding section takes uniform chain constructions
forward by only an 'infinitesimal' step, to ranks finite over
a limit. The crucial step will be to complete the limit steps.
The following discussion considers simple types of uniform 
chains of 'limit' rank, but really serves to illustrate the 
fact that such simple 'algebraic' type constructions 
(involving sums/products) won't get us anywhere near 
exploring the 'limit' explosion in the variety of uniform 
chains.

\secSkip
\Section{[4.1]}
\Title{Atomic chain products}

\parSkip\noindent
If $I$ is a chain (linear order) and for each $i\in I$, $A_i$ is a 
chain, then the expression
    $A = \sum_{i\in I}A_i$
defines a chain in the usual fashion. Sums easily allow to express 
finite products. The (anti-lexicographically ordered cartesian) 
product of two chains,
    $A = B_1B_2$
is expressed by
    $A = \sum_{b_2\in B_2}B_1$
and this definition iterates through all finite products. Some infinite
products can be equally easily expressed by infinite sums (of finite
products) such as e.g.,
    $$\N^{\N} = \prod_{i\in\N}\N = \sum_{i\in\N}\N^i.$$
More revealing is another example,
    $$\Z^{\N} 
    = \prod_{i\in\N}\Z
    = \sum_{-i\in\H}\Z^i + \sum_{i\in\N}\Z^i,$$
which suggests how to express infinite products of chains. Let 
    $H_0=1$ 
and
    $H_i\in\{\N,\Z\}$ for $i\ge 1$.
What should
    $\prod_{i\in\N}H_i$
be? Define $P_0=0$, and for $i>0$,
\item{}
    $P_i = \prod_{k<i}H_k$
\item{}
    $P_{-i} = P_i$ if $H_i = \Z$,
\item{}
    $P_{-i} = 0$ otherwise (if $H_i = \N$).
\par\noindent
Take 
    $A=\prod_{i\in\N}H_i$
to be the chain expressed by
    $A = \sum_{-i\in\H}P_{-i}\H + \sum_{i\in\N}P_{i}\N$
or, compactly
    $A = \sum_{i\in\Z}P_{i}\N^{\sgn(i)}$,
on the understanding that
    $\N^{+} = \N$
    $\N^{-} = \H$.
This definition quickly reproduces the two special cases mentioned; in
$\N^{\N}$ all $P_{-i}$ vanish and $P_{i+1}=\N^i$; in $\Z^{\N}$ for
all $i$, $P_{-i-1}=P_{i+1}=\Z^i$. It also applies in the case of finite
products:

\parSkip
\Number{(4.1)}
\Title{Lemma}
$P_n=\sum_{i\in(-n,n)}P_i\N^{\sgn(i)}$ for $n>1$. 

\Proof
Observe that 
    $P_0=0$, 
    $P_1=H_0=1$, 
    $P_2=H_0H_1=H_1\in\{\N,\Z\}$,
and
    $$
    \sum_{i\in(-2,2)}P_i\N^{\sgn(i)}    
    = P_{-1}\H + 0 + P_1\N,
    $$
which is $\N$ if $H_1=\N$ and $\Z$ if $H_1=\Z$, hence $P_2$. 
Assuming the statement of the lemma true for some $n>1$, we get
    $$
    P_{n+1}
    = \prod_{i\le n}H_i
    =_{\rm dfn} (\prod_{i<n}H_i)H_n
    = P_nH_n
    \in\{P_n\N,P_n\Z\}
    ,$$
while
    $$\eqalign{
    \sum_{i\in[-n,n]}P_i\N^{\sgn(i)} \cr
    &= P_{-n}\H + \sum_{i\in(-n,n)}P_i\N^{\sgn(i)} + P_n\N \cr
    &=_{\rm ind} P_{-n}\H + P_n + P_n\N \cr
    &= P_{-n}\H + P_n\N, \cr
    }$$
which is $P_n\N$ if $H_n=\N$ and $P_n\Z$ if $H_n=\Z$, hence 
equal to $P_{n+1}$.\Qed

Thus, for each $n>1$,
    $$\eqalign{
    \sum_{i\in\Z}P_{i}\N^{\sgn (i)}
    &= \sum_{-i\in\H}P_{-i}\H + \sum_{i\in\N}P_{i}\N \cr
    &= \sum_{-i\le -n}P_{-i}\H + P_n + \sum_{i\ge n}P_{i}\N, \cr
    }$$
and in this representation the middle summand is called the
'central' $P_n$.

\parSkip
\Number{(4.2)}
\Title{Lemma}
$A=\prod_{i\in\N}H_i$ is (right) uniform, and
    $$\f A=\sum_{i\in\N}P_{i}\N.$$


\Proof
By the observation on tails of finite products,
    $$
    \f P_n
    = \f\prod_{i<n}H_i
    = \sum_{i<n}(\prod_{k<i}H_i)\f H_i
    = \sum_{i<n}(\prod_{k<i}H_i)\N
    = \sum_{i<n}P_i\N
    .$$
Let $x\in A$. By (2.3) chose $n>1$ such that $x$ is in the 'central'
$P_n$. Thus, 
    $$
    \f_x A
    = \f P_n + \sum_{i\ge n}P_{i}\N
    = \sum_{i\in\N}P_{i}\N
    $$
as claimed. Alternatively, let $P_n\N^{\sgn(n)}, n\in\Z$ be the 
summand containing $x$ and set $x=(x',x'')$ in the cartesian 
representation of $P_n\N^{\sgn(n)}$. If $n\ge 0$,
    $$
    \f_x A
    = \f_x P_n\N + \sum_{i>n}P_{i}\N
    = \f_{x'} P_n + P_n\f_{x''}\N + \sum_{i>n}P_{i}\N
    = \f P_n + P_n\N + \sum_{i>n}P_{i}\N
    = \sum_{i\in\N}P_{i}\N
    .$$
If $n=-m<0$ (and consequently $P_{-m}=P_m$),
    $$\eqalign{
    \f_x A
    &= \f_x P_{-m}\H + \sum_{i\in (-m,m)}P_{i}\N^{\sgn(i)}
    + \sum_{i\ge m}P_{i}\N \cr
    &= \f_{x'} P_m + P_{-m}\f_{x''}\H + P_m
    + \sum_{i\ge m}P_{i}\N \cr
    &= \f P_m + P_m\times\fin + \sum_{i>m}P_{i}\N \cr
    &= \sum_{i\in\N}P_{i}\N, \cr
    }$$
where $P_m\times\fin$ indicates 'finitely many copies' 
of $P_m$.\Qed

\secSkip
\Section{[4.2]}
\Title{Infinite products of uniform chains}

\parSkip\noindent
Given chains $A_i,i\in\N$, define 
    $A=\prod_{i\in\N}\tup{A_i,x_i}$
by fixing a 'zero' $x_i\in A_i$ and setting
    $A=\sum_{k\in\H}P_k\j_{x_i}A_k
      +\sum_{k\in\N}P_k\f_{x_i}A_k$,
where
    $P_k\prod_{i<k}A_i$.
In general this 'product' depends on the choice of 'zeros'. In
case of finite products it doesn't, e.g. (using right distributivity
of ordered cartesian products),
    $$
    P_2
    = A_1A_2
    = A_1(\j_{x_2}A_2 + \f_{x_2} A_2)
    = A_1\j_{x_2}A_2 + A_1\f_{x_2} A_2,
    $$
and using the 'telescope efect' $\f A = 1+\f A$,
    $$\eqalign{
    P_3 
    &= A_1A_2A_3 \cr
    &= A_1A_2(\j_{x_3}A_3 + 1 + \f_{x_3} A_3) \cr
    &= A_1A_2\j_{x_3}A_3 + A_1A_2 + A_1A_2\f_{x_3} A_3) \cr
    &= A_1A_2\j_{x_3}A_3 
     + A_1(\j_{x_2}A_2 + \f_{x_2} A_2) + A_1A_2\f_{x_3} A_3) \cr
    &= A_1A_2\j_{x_3}A_3 
     + A_1\j_{x_2}A_2 + A_1\f_{x_2} A_2 + A_1A_2\f_{x_3} A_3), \cr
    }$$
etc, for any, hence independently of the, choice of $x_i\in A_i$. 
Thus, assuming for some arbitrary $n$,
    $
    P_n 
    = \prod_{i<n}A_i
    = \sum_{-n<-k<0}P_k\j_{x_i}A_k
    + \sum_{0\le k<n}P_k\f_{x_i}A_k
    $,
it follows that
    $$\eqalign{
    P_{n+1}
    &= \prod_{i\le n}A_i \cr
    &= P_nA_n \cr
    &= P_n(\j_{x_n}A_n + 1 + \f_{x_n} A_n) \cr
    &= P_n\j_{x_n}A_n + P_n + P_n\f_{x_n} A_n) \cr
    &=_{\rm ind} P_n\j_{x_n}A_n 
     + \sum_{-n<-k<0}P_k\j_{x_i}A_k
     + \sum_{0\le k<n}P_k\f_{x_i}A_k
     + P_n\f_{x_n} A_n) \cr
    &= \sum_{-n\le -k<0}P_k\j_{x_i}A_k
     + \sum_{0\le k\le n}P_k\f_{x_i}A_k, \cr
    }$$
for any choice of $x_n\in A_n$. In the case of uniform chains $A_i$ 
there is a 'canonical' choice of zeros:

\parSkip
\Number{(4.3)}
\Title{Lemma}
For every $A\in\UC(<\omega)$ there are $x\in A$ 
and $B_i\in \UC(i)\cup\ens{0}$ such that
    $$\j_xA = \sum_{i<\rk{A}}B_i\H.$$

\Proof
In case $\rk{A} = 1$ ($A\in\ens{\N,\Z}$) the 'ordinary' zeros are 
candidates for $x$ (in case of $A=\Z$, any element will do). Then
$\j_xA\in\ens{0,\H=1\dot\H}$ satisfies the conclusion with
$B_0\in\ens{0,1}$. Assume the statement holds for all $A\in\UC(<n)$.
Given $A\in\UC(n)$ fix $s\in\SS(n)$ such that 
    $A=\A(s)$.

If $s(n)=n$, 
    $\A(s)=\A(s^-)\Z$. 
Set 
    $A'=\A(s^-)$. 
By hypothesis there are $x'\in A'$ and 
$B_i'\in\UC(i)\cup\ens{0}$ for $i<n$ such that 
    $\j_{x'}A'=\sum_{i<n}B_i'\H$. 
Let $x''\in\Z$ (e.g., $x''=0$) and set $x=(x',x'')\in A'\Z$. Then
    $$\eqalign{
    \j_xA
    &= \j_{(x',x''}A'\Z \cr
    &= A'\j_{x''}\Z + \j_{x'}A' \cr
    &= A'\H + \sum_{i<n}B_i'\H \cr
    &= \sum_{i\le n}B_i\H \cr
    }$$
upon setting $B_i=B_i'$ for $i<n$ and $B_n=A'$.

If $s(n)<n$, 
    $\A(s) 
    = \A(s^+) + \A(s^-)\N$. 
Set 
    $A'=\A(s^+)$ 
but now let $x'\in A'$ and $B_i'\in\UC(i)\cup\ens{0}$ for $i<n$ 
such that 
    $\j_{x'}A'=\sum_{i<n}B_i'\H$. 
Then
    $$
    \j_xA
    = \j_xA'
    = \sum_{i<n}B_i'\H
    = \sum_{i\le n}B_i\H,
    $$
upon setting $B_i=B_i'$ for $i<n$ and $B_n=0$.\Qed

This also shows that in the representation of $\j_xA$ we have
    $$B_i = \cases{
    0            & if $s(i)<i$ \cr
    \A(s\ric i) & if $s(i)=i$ \cr
    }$$
for the $s\in\SS(\rk{A})$ such that $A=\A(s)$, so that the
representing $B_i$ are independent of the choice of $x$ (subject
to the condition that they allow such a representation).

\secSkip
\Section{[4.3]}
\Title{Infinite uniform chain products}

\parSkip\noindent
The infinite uniform chain product 
    $A = \prod_{i\in\N}H_i$
is defined to be
    $$A 
    = \sum_{-i\in\H}P_i\j_{x_i}H_i
    + \sum_{i\in\N}P_i\f_{x_i}H_i$$
(with $P_i=\prod_{k<i}H_k$ as before)
for a canonical choice of $x_i\in H_i$. Namely, then,
    $$A 
    = \sum_{-i\in\H}P_iK_i
    + \sum_{i\in\N}P_i\f H_i$$
where the chains
    $K_i=\j_{x_i}H_i=\sum_{k<\rk{H_i}}B_k^{(i)}\H$
are independent of the choice of the $x_i$,
viz,
   $$B_k^{(i)} = \cases{
    0              & if $s_i(k)<k$ \cr
    \A(s_i\ric k) & if $s_i(k)=k$ \cr}$$
for the $s_i\in\SS(\rk(H_i)$ such that $H_i=\A(s_i)$.

\parSkip
\Number{(4.4)}
\Title{Corollary}
    $$\prod_{i\in\N}H_i
    = \sum_{n\in\Z}P_nK_n$$
where for $s_i\in\SS(\rk{H_i}))$ such that $H_i = \A(s_i)$,
\dispM{
    P_{-i} 
    = P_i = \prod_{k<i}H_i, \; 
    K_{-i} 
    = \sum_{k<i,s_i(k)=k}\A(s_i\ric k)\N, \; 
    K_i    
    = \sum_{k<i}\A(s_i\ric k)\N. \; 
}

\parSkip
\Number{(4.5)}
\Title{Lemma}
\Subtitle{$\omega$-product lemma}
Infinite uniform chain products of 'length'
$\N$ are uniform. In fact,
    $$\f\prod_iH_i=\sum_i(\prod_{k<i}H_k)\f H_i.$$

\Proof
Choose a representation
    $$A
    = \prod_iH_i
    = \sum_{-i\in\H}P_iK_i + \sum_{i\in\N}P_i\f H_i$$
where $P_i, K_i$ are as before. For $x\in A$, $x=(x',x'')=(\pm n,x'')$ 
for $x'\in\Z$ and $x''\in P_nK_n$ if $x'=-n$ or $x''\in P_n\f H_n$ if
$x'=n$.

If $z=n\ge 0$, since all terms involved in the right summand are
themselves uniform chains,
    $$\eqalign{
    \f_xA
    &= \f (P_n\f H_n) + \sum_{i>n}P_i\f H_i \cr
    &= \f P_n + P_n\f^2 H_n + \sum_{i>n}P_i\f H_i \cr
    &= \sum_{i<n}P_i\f H_i + P_n\f H_n \sum_{i>n}P_i\f H_i \cr
    &= \sum_{i\in \N}P_i\f H_i.\cr
    }$$
Suppose $z=-n<0$. Then (with canonical zeros $x_i\in H_i$)
    $$\eqalign{
    \f_xA
    &= \f_{x''}P_nK_n
    + \sum_{-n<-k<0}P_iK_i
    + \sum_{i\in\N}P_i\f H_i \cr
    &= \f_{x''}P_nK_n
    + \sum_{-n<-i<0}P_i\j_{x_i}H_i
    + \sum_{0\le i<n}P_i\f H_i 
    + \sum_{i\ge n}P_i\f H_i \cr
    &= \f_{x''}P_nK_n
    + P_n
    + \sum_{i\ge n}P_i\f H_i \cr
    &= \f_{x''}P_nK_n
    + P_n\f H_n 
    + \sum_{i>n}P_i\f H_i \cr
    }$$
Fix $p\in P_n,k\in K_n=\j_{x_n}H_n$ such that $x''=(p,k)$, so that
    $$
    \f_{x''}P_nK_n
    = \f_pP_n + P_n\f_kK_n
    = \sum_{i<n}P_i\f H_i + P_n\f_k\j_{x_n}H_n,
    $$
so
    $$\eqalign{
    \f_xA
    &= \sum_{i<n}P_i\f H_i 
    + P_n\f_k\j_{x_n}H_n
    + P_n\f H_n 
    + \sum_{i>n}P_i\f H_i \cr
    &= \sum_{i<n}P_i\f H_i 
    + P_n(\f_k\j_{x_n}H_n + \f_{x_n} H_n)
    + \sum_{i>n}P_i\f H_i \cr
    &= \sum_{i<n}P_i\f H_i 
    + \sum_{i\ge n}P_i\f H_i \cr
    &= \sum_{i\in \N}P_i\f H_i, \cr
    }$$
using right distributivity and the observation that for any uniform 
chain $X$ and any $a < b\in X$,
    $$\f_a\j_bX + \f_bX = \f_aX = \f_bX.\Qed$$


\chapSkip
\Chapter{[5] Chains of first limit rank}

\secSkip\noindent
We now resume the thread left off in section 2 and proceed to 
characterise (uniform) chains of rank $\omega$.

\secSkip
\Section{[5.1]}
\Title{Chains of rank $\omega$}

\parSkip\noindent
Approximating $A\in\UC(\omega)$ requires consideration of
entire 'inclusion' chains of sub-components, essentially over
the 'length' of the approximation given by the rank of $A$.
I.e., we need to incorporate into any 'limit rank uniform
chain constructor' the information that each such $A$ is 
the union of the 'chain of chains'
    $\ens{\cmp{x}{n}{A}}_{n\in\omega}$
(for any fixed $x\in A$), and hence of an 'embedding' chain of 
chains 
    $\ens{\A(s_n)}_{n\in\omega}$ 
where 
    $\ens{s_n}_{n\in\omega}$ 
is a suitable sequence of $\SS$ strings and where the embedding 
    $\A(s_n)\hookrightarrow\A(s_{n+1})$
is the 'canonical' embedding. To state the obvious fact:

\parSkip
\Number{(5.1)}
\Title{Lemma}
\Subtitle{Copy lemma}
For $A\in\UC(\omega)$, $x\in A$ and $n\in\omega$,
\dispM{
    \cmp{x}{n+1}{A} \isom \cases{
    \cmp{x}{n}{A}+B_{n}.\N 
    & 
    if $\cmp{x}{n+1}{A}/{n} = \N$ 
    and 
    $\cmp{x}{n}{A}$ has no predecessors in
    $\cmp{x}{n+1}{A}$ 
    \cr
    B_{n}+\cmp{x}{n}{A}.\N 
    & 
    if $\cmp{x}{n+1}{A}/{n} = \N$ 
    and 
    $\cmp{x}{n}{A}$ does have predecessors in
    $\cmp{x}{n+1}{A}$ 
    \cr
    \cmp{x}{n}{A}.\Z 
    & 
    if $\cmp{x}{n+1}{A}/{n} = \Z$ 
    \cr
    }
}    
where $B_{n}\in\UC(n)$ such that
    $\f B_{n} \isom \f \cmp{x}{n}{A}$
and where by 'predecessor' in $\cmp{x}{n+1}{A}$ we mean '$n$-component
predecessor'.

\parSkip
We shall now turn this observation into a characterisation of
$A\in\UC(\omega)$.

\parSkip
\Number{(5.2)}
\Title{Definition}

(a)
Call a sequence 
    $\ens{s_{n}}_{n}$
of strings
    $s_{n}\in\SS(n+1)$
{\it realised by} $\tup{A,x}$, where $x\in A\in\UC(\omega)$,
if 
    $\cmp{x}{n}{A} = \A(s_n)$
for all $n$.
Call 
    $\ens{s_{n}}_{n}$
as above {\it realised by} $A\in\UC(\omega)$ if it is
realised by some $\tup{A.x}$, $x\in A$, and {\it realised in} 
$\UC(\omega)$ if it's realised by some $A\in\UC(\omega)$.

(b) 
Define the {\it embedding sequence} $\ps\in\Fun{\omega}{(\omega+1)}$
of $\tup{A,x}$ by
\dispM{
    \ps(n) = \ps^{A,x}(n) = \cases{
    i\in\omega 
    & if $\cmp{x}{n+1}{A}/n\isom\N$ and 
      $\cmp{x}{n}{A}$ has $i$ 
      predecessors in $\cmp{x}{n+1}{A}$, \cr
    \omega     
    & else (i.e., if $\cmp{x}{n+1}{A}/n\isom\Z$). \cr
    }
}

(c)
Define the {\it type sequence} $\tp\in\Fun{\omega}{\ens{0,1,\omega}}$ 
of a sequence $\ps\in\Fun{\omega}{(\omega+1)}$ by
\dispM{
    \tp(n) = \cases{
    \ps(n)     & if $\ps(n)\in\ens{0,\omega}$ \cr
    1          & if $\ps(n)\in (0,\omega)$. \cr
    }
}
Define the type sequence $\tp^{A,x}(n)$ of $\tup{A,x}$ as
the type sequence of $\ps^{A,x}(n)$.

\parSkip
In the following, write $\ps_n,\tp_n$ for $\ps(n),\tp(n)$.

\parSkip
Thus every $A\in\UC(\omega)$ is determined by a pair
    $\tup{\ens{s_n}_{n},\ps}$ 
where $\ens{s_n}_{n}$ is a sequence 
realised by $\tup{A,x}$ for some $s$ in $A$, and where 
$\ps$ is the embedding sequence for that pair $\tup{A,x}$.
Moreover, any two pairs 
    $\tup{\ens{s_n}_{n},\ps}$, $\tup{\ens{s_n'}_{n},\ps'}$ 
determining $A$ in that fashion agree on a final segment: If the 
corresponding elements of $A$ are $x$ and $x'$, then $s_n=s_n'$ 
and $\ps_n=\ps_n'$ for all $n\ge m$ such that $x\rho_m x'$.

\parSkip
\Number{(5.3)}
\Title{Lemma}
\Subtitle{String copy lemma}
Let 
    $\ens{s_n}_{n}$, $s_n\in\SS(n+1)$
be realised by $\tup{A, x}$, i.e., such that
    $\cmp{x}{n}{A} = \A(s_n), n\in\omega$.
Then, 
\dispM{
    s_{n} = \cases{
    s_{n+1}^+ & if $\ps_n=0$ iff $\tp_n = 0$, \cr
    s_{n+1}^- & if $\ps_n>0$ iff $\tp_n > 0$, \cr
    } \eqno{(5.1)}
}
and in either case
    $s_{n}^- = s_{n+1}^=$. 

\Proof
$s_n=s_{n+1}^+$ just if the first of the above cases for
$\cmp{x}{n+1}{A}$ holds, else $s_n=s_{n+1}^-$; the remainder is
the string lemma.\Qed

\parSkip
\Number{(5.4)}
\Title{Definition and Lemma}
\Subtitle{Aligned string sequences}
Call a sequence of strings $\ens{s_n}_{n}$ {\it aligned} 
with $\ps$ (resp., $\tp$), if condition (5.1) holds, and simply
{\it aligned}, if $s_n\in\ens{s_{n+1}^{\pm}}$ for all $n$. 
For such an aligned sequence define 
    $a_n=a^{A,x}_n$,
    $b_n=b^{A,x}_n$, 
$n\in\omega$ such that for all $n>0$,
\dispM{
    s_{n} = s_n^=\tac{a_n,b_n}.
}
(We let $a_0$ undefined and set $b_0=0$, so $s_0=\tup{b_0}$.)
Then,
\dispM{
    b_{n} = \cases{
    b_{n+1} & if $s_{n}=s_{n+1}^+$ iff $\ps_n=0$ iff $\tp_n=0$, \cr
    a_{n+1} & if $s_{n}=s_{n+1}^-$ iff $\ps_n>0$ iff $\tp_n>0$. \cr
    } \eqno{(5.2)}
}
Also, if we set $A_n=\A(s_n)$, then $\ens{A_n}_n$ 
the sequence $\ens{A_n}_{n}$ satisfies
\dispM{
    A_{n+1} \isom \cases{
    A_n\Z 
    & 
    if $\ps_n=\omega$, 
    \cr
    A_n+B_n\N, 
    &
    for $B_n\in\UC(n)$ with $\f B_n = \f A_n$, 
    if $\ps_n=0$,
    \cr
    B_n+A_n\N
    &
    for $B_n\in\UC(n)$ with $\f B_n = \f A_n$, 
    if $\ps_n\in (0,\omega)$. \cr } \eqno{(5.3)}
}

It follows, for instance, that if $\ps_n = 0$ 
for all $n\ge m$ for some $m$, then 
\dispM{
    \eqalign{
    \hbox{(a)} \qquad &
    (\forall n)
    [s_{n}^- \sub s_{n+1}^-], 
    \hbox{ and, for the least such $m$, }
    \cr
    \hbox{(b.1)} \qquad &
    (\forall n < m)[\tp_n=0\then(\exists k, n < k \le m)
    [\tp_k > 0 \hbox{ and hence }
     s_{k} \sub s_{k+1}]],
    \hbox{ i.e., } 
    \cr
    &
    (\forall n < m)[\tp_n=0\then(\exists k, n < k \le m)
    [\tp_k > 0 \hbox{ and hence } 
    s_{n}(n) = b_{n} = a_{k} = s_{k+1}(k)]],
    \hbox{ and } 
    \cr
    \hbox{(b.2)} \qquad &
    (\forall n\ge m)
    [s_{n}(n) = b_{n} = b_{m}], 
    \cr
    }
}
while if $\ps_n > 0$ for infinitely many $n$, then 
\dispM{
    \eqalign{
    \hbox{(a)} \qquad &
    (\forall n)
    [s_{n}^- \sub s_{n+1}^-],
    \hbox{ and } 
    \cr
    \hbox{(c.1)} \qquad &
    (\exists^{\infty} n)
    [\tp_n > 0 \hbox{ and hence }
     s_{n} \sub s_{n+1}],
    \hbox{ i.e., } 
    \cr
    &
    (\forall n)[\tp_n=0\then(\exists k\ge n)
    [\tp_k > 0 \hbox{ and hence } 
    s_{n}(n) = b_{n} = a_{k} = s_{k+1}(k)]].
    \cr
    }
}
Thus, in the first case the entire sequence 
    $\ens{s_n}_{n}$ 
can be recovered from 
\dispM{
    \eqalign{
    \hbox{(b.1)} \qquad &
    \hbox{ $s_m$, in particular $s_m^-$ and $b_{m} = s_m(m)$
    (for the least such $m$ as stipulated), and } 
    \cr
    \hbox{(b.2)} \qquad &
    \hbox{ the chain of strings $s_{m+1}^-\sub s_{m+2}^-\sub\ldots$, }
    \cr
    }
}
while in the second case we can recover the sequence
    $\ens{s_n}_{n}$ 
from
\dispM{
    \eqalign{
    \hbox{(c.1)} &
    \hbox{ the set of $m$ where $\tp_m > 0$ and }
    \cr
    \hbox{(c.2)} &
    \hbox{ the chain of strings $s_{0}^-\sub s_{1}^-\sub\ldots$. }
    \cr
    }
}
In either case, the sequence $\ens{s_n}_{n}$  can
be recovered from the chain $\ens{s_n^-}_{n}$ and
$\tp$ (or $\ps$), plus possibly a single element $b$ (the 
$b_m$ of the first case). This allows to economise the
descriptors of $A\in\UC(\omega)$ (compared with the full 
set $\tup{\ens{s_n}_{n},\ps}$). In fact, we can replace 
$\ens{s_n}_{n}$ by a single string $\ss\in\SS(\omega+1)$.

\parSkip
\Number{(5.5)}
\Title{Definition}
(a) 
\Subtitle{Trunk-string}
Define the {\it trunk} of $\ens{s_n}_{n}$, 
    $\ss = \ss^{A,x}\in\SS(\omega+1)$,
by setting
\dispM{
    \eqalign{
    \ss\ric\omega 
    & 
      = \bigcup_{n}s_{n+1}^- \cr
    \ss(\omega) 
    & 
      = \cases { 
            b_{m} 
            & where $m$ is minimal such that $\ps_n = 0$ 
              for all $n\ge m$, if such exists \cr
            \omega 
            & else. \cr }
    }
}

(b) 
\Subtitle{Type-partition}
Let
\dispm{
    \osc (\tp)
    = \osc (\tup{A,x} )
    = 1+\cns{v\in\omega:\tp_{v+1}\not=\tp_v}.
}
Note that $\osc(\tp)\in\omega\cup\ens{\omega}$.
Define the {\it partition 
    $\omega = \bigcup_{n<\osc(\tp)}[v_{n},v_{n+1})$
induced by} $\tp\in\Fun{\omega}{\ens{0,1,\omega}}$, by
\itemitem{(i)}
    $v_0 = 0$,
\itemitem{(ii)}
    $v_{n+1} = \min\ens{v>v_n:\tp_v\not=\tp_{v_n}}$,
    for $n<\osc(\tp)$  and
\itemitem{(iii)}
    $v_{n+1} = \omega$ 
    if $\osc(\tp)=n+1$, i.e., if $n$ is such 
    that for all $v>v_n$, $\tp_v=\tp_{v_n}$.
\par\noindent
This is the coarsest partition of $\omega$ such that
$\tp$ is constant on each block $[v_{n},v_{n+1})$ (i.e.,
such that also adjacent intervals are of different types). 
Call $[v_{n},v_{n+1})$ {\it of type} $\tp_{v_n}$. 
The partition induced by $\ps$, resp., by $\tup{A,x}$, is 
defined as the partition induced by the type sequence $\tp$ 
of $\ps$, resp., by the type sequence $\tp^{\tup{A,x}}$ of 
$\tup{A,x}$.

\parSkip
Let now $s_k\in\SS(k+1)$ be such that $\cmp{x}{k}{A} = \A(s_k)$,
for $x\in A\in\UC(\omega)$ and all $k\in\omega$.

\parSkip
\Number{(5.6)}
\Title{Lemma}
\item{(i)}
For an interval $[v_{n},v_{n+1})$ of type 0,
    $b_{v_n} = b_{v_n+1} = \ldots = b_{v_{n+1}-1} 
    = a_{v_{n+1}} = \ss(v_{n+1})$;
hence:
\item{(ii)}
For $k$ in an interval $[v_{n},v_{n+1})$ of type 0,
    $s_{k}=\ss^{A,x}\ric k\tac{b_{v_n}}$.
Moreover: 
\item{(iii)}
If $v_{n+1}<\omega$ for an interval $[v_{n},v_{n+1})$ of type 0, 
then also
    $b_{v_n} = s_{v_{n+1}+1}(v_{n+1}))$.
\item{(iv)}
For $k$ in an interval $[v_{n},v_{n+1})$ of type $1$ or $\omega$,
    $s_{k}=\ss^{A,x}\ric k+1$. 

\Proof
By induction over $i\in [v_{n},v_{n+1})$ using the string copy
lemma.
\Qed

\parSkip
\Number{(5.7)}
\Title{Corollary}
\item{(i)}
For $k$ in $[v_{n},v_{n+1})$ of type 0,
    $s_{k+1}=\ss^{A,x}\ric k \tac{\ss(v_{n+1})}$;
and:
\item{(ii)}
If $[v_{n},v_{n+1})$ is of type 0 and finite, then in fact
    $s_{k+1}=\ss^{A,x}\ric k \tac{s_{v_{n+1}+1}(v_{n+1})}$.
\item{(iii)}
For $k$ in $[v_{n},v_{n+1})$ of type 1 or $\omega$,
    $s_{k+1}\sub\ss^{A,x}$.
Moreover:
\item{(iii)}
On $[v_{n},v_{n+1})$ of type $\omega$,
    $s_{k+1}=s_k\tac{k}\sub\ss^{A,x}$ 
\Qed

\parSkip
\Number{(5.8)}
\Title{Definition and Lemma}

(a)
Call
    $\ss\in\SS(\omega+1)$
and 
    $\ps\in\Fun{\omega}{(\omega+1)}$,
{\it compatible} if
    $\ss(k)=k$ iff $\ps_k=\omega$ 
and for every $n$, 
    if $\ps_n=0$ then for the least $m>n$ such that $\ps_m>0$,
    $\ss(m)\le n$. I.e., if 
\dispM{
    \eqalign{
    \ss(k) = k 
    \qquad &
    \hbox{iff $\ps_k = \omega$}
    \cr 
    \ss(v_{k+1})\le v_k 
    \qquad &
    \hbox{whenever $\ps_{v_k}=0$ (and hence $\ps_{v_{k+1}} > 0$)}
    \cr } \eqno{(5.4)}
}
Similarly, call
    $\ss\in\SS(\omega+1)$
and 
    $\tp\in\Fun{\omega}{\ens{0,1,\omega}}$
{\it compatible} if the preceding condition holds
with '$\tp$' in place '$\ps$'.

(b)
Call the sequence
    $\ens{s_n}_n$, $s_n\in\Fun{n+1}{\omega}$,
{\it derived from} $\tup{\ss,\ps}$ if
\dispM{
    s_n = \cases{
    \ss\ric n+1                  
    & if $\ps_n>0$ \cr
    \ss\ric n\tac{\ss(w_n)} 
    & for $w_n=\inf\ens{v\ge n:\ps_v>0}$, 
    if $\ps_n=0$. \cr
    } \eqno{(5.5)}
}
Similarly, call the sequence
    $\ens{s_n}_n$, $s_n\in\Fun{n+1}{\omega}$,
{\it derived from} $\tup{\ss,\tp}$, 
if the preceding condition holds with '$\tp$'
replacing '$\ps$' throughout.

Then 
\item{(i)}
$\tup{\ss,\ps}$ are compatible iff $\tup{\ss,\tp}$ are;
\item{(ii)}
$\ens{s_n}_n$, $s_n\in\Fun{n+1}{\omega}$, is derived from 
$\tup{\ss,\tp}$ iff it is derived from $\tup{\ss,\tp}$; 
\item{(iii)}
$s_n\in\SS(n+1)$ whenever $\tup{\ss,\ps}$ (resp.,
$\tup{\ss,\tp}$) are compatible and $\ens{s_n}_n$ is derived
from $\tup{\ss,\tp}$.
\item{(iv)}
$s_n\in\ens{s_{n+1}^{\pm}}$ whenever $\tup{\ss,\ps}$ (resp.,
$\tup{\ss,\tp}$) are compatible and $\ens{s_n}_n$ is derived
from $\tup{\ss,\tp}$, in fact $\ens{s_n}_n$ is aligned with
$\ps$ (resp., $\tp$).

\Proof

(i) and (ii) are clear, since $\ps_k = \omega$ iff 
$\tp_k = \omega$. and $\ps_k = 0$ iff $\tp_k = 0$. 

(iii) We need to show that $s_n(i)\le i$ for all $i<n+1, n\in
\omega$.  This is clear if $\ps_n>0$ since then $s_n(i)=\ss(i)
\le i$ for all $i<n+1$. If $\ps_n=0$ then $s_n(i)=\ss(i)\le i$
for all $i<n$, so it remains to show that $s_n(n)\le n$ in
that case also. We have $n\in [v_k,v_{k+1})$ of type 0, for
some $k<\osc(\tp)$, and $w_m=v_{k+1}$ for all $m\in [v_k,v_{k+1})$.
By compatibility, $\ss(w_n)=\ss(v_{k+1})\le v_k$, so in
particular $\ss(w_n)\le n$, as required.

(iv) 
From the definition, $s_n\ric n = \ss\ric n$, so
$s_n^- \sub s_{n+1}^-$ for all $n$, and it remains to
show that $b_n\in\ens{a_{n+1}, b_{n+1}}$ for all $n$,
where again, from the preceding, $a_{n+1}=\ss(n)$
for all $n$.
\item{(a)} 
Case $\ps_n=0$ and $\ps_{n+1}=0$. Then $w_n=w_{n+1}>n+1$
and so $b_n=b_{n+1}=\ss(w_n)$, where here, and subsequently,
$w_n$ is as in the definition.
\item{(b)} 
Case $\ps_n=0$ and $\ps_{n+1}>0$. Then $w_n=n+1$,
and $s_{n+1}=\ss\ric n+2$, i.e., $b_{n+1}=\ss(n+1)$
and so $b_n=\ss(w_n)=\ss(n+1)=b_{n+1}$.
\item{(c)} 
Case $\ps_n>0$. Then immediately, $b_n=\ss(n)=a_{n+1}$.
\Qed

\parSkip
\Title{Corollary}
The sequence $\ens{s_n}_n$ derived from a compatible
pair satisfies property (5.1) of realised sequences
(string copy lemma).
Hence $A_n = \A(s_n)$ satisfy (5.3).
    
\parSkip
\Number{(5.9)}
\Title{Definition and Lemma} 
Call
    $\ss,\ps,a,b$,
{\it compatible}, where
    $\ss\in\SS(\omega+1)$,
    $\ps\in\Fun{\omega}{(\omega+1)}$
and
    $a,b\in\SS(\omega)$,
if
    $\ss,\ps$ compatible
(i.e, $\ss(k) = k$ whenever $\tp_k = \omega$),
    $\ss\ric k+1 = \ss\ric k\tac{a_{k}}$
and
    $b_{k+1} = b_{k}$ whenever $\ps_k = 0$,
    $a_{k+1} = b_{k}$ whenever $\ps_k\in\omega\rid{0}$.
Call the sequence
    $\ens{s_n}_n$, $s_n\in\SS(n+1)$,
{\it derived from} $\tup{\ss,\ps,a,b}$ if
    $s_n = \ss\ric n-1\tac{a_{n},b_{n}}$.
If 
    $\ens{s_n}_n$, $s_n\in\SS(n+1)$,
derives from $\tup{\ss,\ps,a,b}$, then it
derives from $\tup{\ss,\ps}$, hence satifies (5.1).

\Proof
If $\ps_n\not=0$, then either $\ps_n = \omega$ and
    $s_n = \ss\ric n+1$
by definition, or $\ps_n \in (0, \omega)$ and
    $s_n = \ss\ric n-1\tac{a_{n},b_{n}}
    = \ss\ric n-1\tac{a_{n},a_{n+1}}
    = \ss\ric n\tac{a_{n+1}}
    = \ss\ric n+1$,
by induction. If $\ps_n=0$ then
    $s_n 
    = \ss\ric n-1\tac{a_{n},b_{n}}
    = \ss\ric n\tac{b_{n}}
    = \ss\ric n\tac{w_{n}}$
since $b_{k} = b_{k+1}$ for all $k\in [n,w_n)$
by compatibility.\Qed

\parSkip
\Number{(5.10)}
\Title{Corollary}
\Subtitle{Summary}
The following summarises the properties of the strings
    $\ss = \ss^{A,x},
     \ps = \ps^{A,x},
     a = a^{A,x},
     b = b^{A,x}$.

\item{(a)}
    $a_0$ is undefined,
    $b_0 = 0$ and
    $\cmp{x}{0}{A} \isom \A(\ss\ric 0\tac{b_0}) = \A(\tup{0}) = \ONE$;

\item{(b)}
If
    $\ps(n) = \omega$ 
then
\itemitem{(i)}
    $\cmp{x}{n+1}{A} \isom \cmp{x}{n}{A}\Z$;
\itemitem{(ii)}
    $a_{n+1} = \ss(n) = b_n$ and $b_{n+1} = n+1$, so
\dispM{
    \eqalign{
    \hbox{(ii.1)} 
    \qquad &
    \cmp{x}{n}{A} 
    \isom \A(\ss\ric n\tac{b_n});
    \cr
    \hbox{(ii.2)} 
    \qquad &
    \cmp{x}{n+1}{A} 
    \isom \A(\ss\ric n\tac{b_n,n+1});
    \cr
    }
}

\item{(c)}
If 
    $\ps(n) = 0$, 
then 
\itemitem{(i)}
    $\cmp{x}{n+1}{A} \isom \cmp{x}{n}{A} + B_{n}\N$
for $B_{n}\in\UC(n)$ such that
    $\f B_{n}\isom\f\cmp{x}{n}{A}\isom\A(\ss\ric n\tac{0})$;
\itemitem{(ii)}
    $a_{n+1} = \ss(n)$ and $b_{n+1} = b_n \le n$, so
\dispM{
    \eqalign{
    \hbox{(ii.1)} 
    \qquad &
    \cmp{x}{n}{A} 
    \isom \A(s_n)
    \isom 
        \A(s_{n+1}^+)
        \isom \A(s_{n+1}^=\tac{b_{n+1}})
        \isom \A(\ss\ric n\tac{b_n}),
    \cr
    &
    B_{n} 
    \isom 
        \A(s_{n+1}^-)
        \isom \A(s_{n+1}^=\tac{a_{n+1}})
        \isom \A(\ss\ric n\tac{\ss(n)}) 
    \cr
    \hbox{(ii.2)} 
    \qquad &
    \cmp{x}{n+1}{A}
    \isom \A(s_{n+1})
    \isom 
        \A(s_{n+1}^=\tac{a_{n+1},b_{n+1}})
        \isom \A(\ss\ric n\tac{\ss(n),b_n});
    }
}
\itemitem{(iii)}
    $\ps(n)=0$ is the number of ($n$-component) predecessors of 
    $\cmp{x}{n}{A}$ in $\cmp{x}{n+1}{A}$, and $\cmp{x}{n}{A}$ is
    an initial $n$ component of $\cmp{x}{n+1}{A}$.

\item{(d)}
If $\ps(n) \in \omega\rid\ens{0}$, then 
\itemitem{(i)}
    $\cmp{x}{n+1}{A} \isom B_{n} + \cmp{x}{n}{A}\N$
for $B_{n}\in\UC(n)$ such that
    $\f B_{n} \isom\f\cmp{x}{n}{A}\isom\A(\ss\ric n\tac{0})$;
\itemitem{(ii)}
    $a_{n+1} = \ss(n)$, but also $a_{n+1} = b_n$, 
    and $b_{n+1} \le n$, so
\dispM{
    \eqalign{
    \hbox{(ii.1)} 
    \qquad &
    B_{n} 
    \isom 
        \A(s_{n+1}^+)
        \isom \A(s_{n+1}^=\tac{b_{n+1}})
        \isom \A(\ss\ric n\tac{b_{n+1}}), 
    \cr
    &
    \cmp{x}{n}{A}
    \isom 
        \A(s_n) 
        \isom \A(s_{n+1}^-) 
        \isom \A(\ss\ric n\tac{a_{n+1}}) 
        \isom \A(\ss\ric n\tac{b_n}) 
    \cr
    \hbox{(ii.2)} 
    \qquad &
    \cmp{x}{n+1}{A}
    \isom \A(\ss\ric n\tac{\ss(n),b_{n+1}}); 
    \isom \A(\ss\ric n\tac{b_n,b_{n+1}}); 
    \cr
    }
}
\itemitem{(iii)}
    $\cmp{x}{n}{A}$ is a successor $n$ component of $\cmp{x}{n+1}{A}$
    and $\ps(n)>0$ is the number of ($n$-component) predecessors of 
    $\cmp{x}{n}{A}$ in $\cmp{x}{n+1}{A}$.

\parSkip
We are now ready to state the characterisation of $A\in\UC(\omega)$.

\parSkip
Let $\CSI(\omega+1)$ denote the set of pairs $\tup{\ss,\ps}$
of compatible
    $\ss\in\SS(\omega+1)$,
    $\ps\in\Fun{\omega}{(\omega+1)}$,
and let $\CSP(\omega+1)$ denote the set of pairs $\tup{\ss,\tp}$ 
of compatible
    $\ss\in\SS(\omega+1)$,
    $\tp\in\Fun{\omega}{\ens{0,1,\omega}}$,
respectively.

\parSkip
\Number{(5.11)}
\Title{Lemma}
\Subtitle{Realisation lemma}
For
    $\tup{\ss,\ps}\in\CSE)\omega+1)$
and let
    $\ens{s_n}_n$, $s_n\in\SS(n+1)$,
be the sequence derived from $\tup{\ss,\ps}$.
Then $\ens{s_n}_n$ is realised in $\UC(\omega)$.

\Proof
We need to construct $A\in\UC(\omega)$ and $x\in A$, 
such that $\ens{s_n}_n$ is realised by $\tup{A,x}$.

It will suffice to define $A$ as the union over the chain 
of chains
    $A = \bigcup_{n\in\omega}\tup{A_n,\phi_n}$
where 
    $\phi_n:A_n \hookrightarrow A_{n+1}$ 
is the inclusion of $A_n$ in $A_{n+1}$ as the $n$-component 
with $\ps_n$ $n$-component predecessors. I.e., where
    $A_{n+1} = \sum_{i\in A_{n+1}/{n}}C_i$
shall be such that
\item{(i)}
if
    $\ps_n\in\omega$
then
    $A_{n+1}/{n} = \N$,
and
    $A_n = C_{\ps_n}$;
\item{(ii)}
if
    $\ps_n=\omega$
then
    $A_{n+1}/{n} = \Z$,
and
    $A_n = C_{0}$.
\par\noindent
The element $x\in A$ for which,
    $A_n = \cmp{x}{n}{A}$ for all $n$
is given as the $x$ for which $\ens{x} = A_0$ 
(really, by the 'thread' $\tup{\phi_n(0)}_n$, since
    $A_0 = \ONE = \ens{0}$).

Set
    $A_n = \A(s_n)$.
Clearly then, 
\item{(a)} 
$A_n\in\UC(n)$;
\item{(b)} 
$A_n$ is an $n$ component of $A_{n+1}$, for each
$n$ (since $\ens{s_n}_n$ are aligned with $\ps$ and
$\ens{A_n}$ satisfy (5.3)), so $A_{n+1}/n$ is as 
required; and 
\item{(c)} 
since $s_n=s_{n+1}^+$ whenever $\ps_n=0$ we can embed
$A_n$ as an initial $n$ component into $A_{n+1}$, while 
$s_n=s_{n+1}^-$ if $\ps_n>0$ allows to embed $A_n$ as a 
successor $n$ component into $A_{n+1}$.

Thus, let $\phi_n$ be the canonical embedding 
    $A_n\hookrightarrow A_{n+1}$ 
that maps $A_n$ isomorphically to the $\ps_n$th $n$ component 
of $A_{n+1}$ if $\ps_n<\omega$, and that so maps $A_n$ to the 
0th $n$ component if $\ps_n=\omega$.


Thus, $\ens{s_n}_n$ is realised by $\tup{A,x}$.
\Qed

\parSkip
Let 
\dispM{
    \A(\ss,\ps)
    \eqno{(5.6)}
} 
denote the element of $\UC(\omega)$ 
represented by $\tup{\ss,\ps}$ via the realisation lemma. 

\parSkip
For such $A=\A(\ss,\ps)$ and the sequence
$\ens{s_n}$ derived from $\tup{\ss,\ps}$: 
\item{(a)}
For some $x\in A$, 
    $\ps^{\tup{A,x}} = \ps$
and for all $n$,
    $\cmp{x}{n+1}{A}/n \isom \N$ if $\ps_n<\omega$,
    $\cmp{x}{n+1}{A}/n \isom \Z$ if $\ps_n=\omega$,
and
    $\cmp{x}{n}{A} = \A(s_n)$.
\item{(b)}
If $\ps_n<\omega$ then the initial $n$ component 
of $\cmp{x}{n+1}{A}$ is (mod $\isom$) given by 
$\A(s_{n+1}^+)$.
The successor $n$-components in $\cmp{x}{n+1}{A}$ 
are given by 
    $A_n \isom \A(s_{n+1}^-) = \A(\ss\ric n+1)$,
so by the tail lemma,
\dispM{
    \f A 
    \isom \sum_{n\in \omega} A_n\N
    \isom \sum_{n\in \omega} \A(\ss\ric n+1)\N.
}

In fact, with $a, b$ defined for 
    $\tup{\ss,\ps}\in\CSI(\omega+1)$ 
as in the aligned string sequence lemma, all of the summary 
applies to $A=\A(\ss,\ps)$. 

\parSkip
\Number{(5.12)}
\Title{Corollary}
\Subtitle{Rank $\omega$ representation lemma}
Every $A\in\UC(\omega)$ has a representation of from (5.6).
In particular, $A\in\UC(\omega)$ is isomorphic to $\A(\ss,\ps)$ 
whenever $\ss=\ss^{\tup{A,x}}$ and $\ps=\ps^{\tup{A,x}}$ for some 
$x\in A$.\Qed

\parSkip
Another corollary is that the representation of $A\in\UC(\omega)$ 
by compatible pairs $\tup{\ss,\ps}$ is unique up to finite initial 
segments of $\ss,\ps$.

\parSkip
\Number{(5.13)}
\Title{Corollary}
\item{(i)}
If $\A(\ss,\ps) = \A(\ss',\ps')$, then $\ss^- = \ss'^-$
and $\ps\aeq\ps'$ (i.e., $\ps\ric(n,\omega)=\ps'\ric(n,\omega)$
for some $n$).
\item{(ii)}
If $\f\A(\ss,\ps) = \f\A(\ss',\ps')$, then $\ss^- = \ss'^-$
(i.e., $\ss\ric\omega = \ss'\ric\omega$).
\Qed

\secSkip
\Section{[5.2]}
\Title{Tails of rank $\omega$}

\parSkip\noindent
We shall verify that by computing tails of $A\in\UC(\omega)$ in 
the representation just given, one obtains the representations 
produced in section 2.
Observe that if
    $\ss\in\SS(\omega+1)$
is compatible with
    $\ps = 0^{\omega}$
(where $0^{\omega}=\tup{0,0,\ldots}$ is a sequence of 
0s of length $\omega$), then it must be the case that
    $\ss(\omega) = 0$
(since for all $n\in\omega$, $w_n=\omega$, i.e, for the partition
induced by $\ps = 0^{\omega}$, 
\dispM{
    \osc(\ps) = 0,v_ 0=0, v_1=\omega, 
    \hbox{ and $[v_0,v_1)$ is of type 0,}
}
so 
    $\ss(\omega)=\ss(v_1)\le v_0 = 0$).  
Hence, for the sequence of strings
    $\ens{s_n}_{n}$
derived from $\tup{\ss,0^{\omega}}\in\CSI(\omega+1)$ we have
    $s_n = \ss\ric n+1\tac{0}$. 

Now observe that for the point
    $x_0\in A_0=\A(\ss,0^{\omega})$
such that $\ps^{\tup{A_0,x_0}} = 0^{\omega}$, we have
    $\cmp{A_0}{n}{x_0} 
    = \A(s_n)
    = \A(\ss\ric n\tac{0})$. 
On the other hand, for the left-most element $x_1$ in the chain
    $A_1= \sum_{n\in\N}\A(\ss\ric n+1)\N$
we have
    $\cmp{A_1}{n}{x_1} 
    = \sum_{i<n}\A(\ss\ric i)\N,
    = \A(\ss\ric n\tac{0})$. 
By the block-lemma it follows that $A_0 \isom A_1$.
Thus,
\dispM{
    \A(\ss,0^{\omega}) 
    = \A(\ss^-\tac{0}, 0^{\omega})
    = \A(\ss\ric\omega\tac{0}, 0^{\omega})
    \isom \sum_{n\in\N}\A(\ss\ric n+1)\N.
}

\parSkip
\Number{(5.14)}
\Title{Corollary}
\Subtitle{Rank $\omega$ tail lemma}
For all compatible $\tup{\ss,\ps}\in\CSI(\omega+1)$,
\dispM{
    \f\A(\ss,\ps) 
    \isom \A(\ss^-\tac{0},0^{\omega}).
}

\Proof
Recall that for $x\in A\in\UC$,
    $\f_xA \isom \sum_{0\le n<\rk{A}}A_n\N$,
where
    $A_n$ is (a representative of the isomorphism class of) the 
    successor $n$-component of $\cmp{x}{n}{A}$ in $\cmp{x}{n+1}{A}$.

But, for $A = \A(\ss,\ps)$ in particular,
\dispM{
    \cmp{x}{n+1}{A} 
    \isom \A(s_{n+1})
    \in \ens{
      \A(s_{n+1}^-)\Z,
      \A(s_{n+1}^+)+\A(s_{n+1}^-)\N
      }
    \hbox{ mod $\isom$, }
}
so
    $A_n = \A(s_{n+1}^-) = \A(\ss\ric n+1)$.

Set now $A=\A(\ss^-\tac{0},0^{\omega})$. Note that the sequence 
derived from $\tup{\ss,0^{\omega}}$ is 
    $\ens{\ss\ric n\tac{0}}_n$. 
Thus for the element $x = 0$ such that 
    $\Hyp{1}{0}{A} = \bf{1} = \ens{x}$, 
it is the case that for all n,
    $\cmp{x}{n}{A} = \A(\ss\ric n\tac{0})$.
Equally, for $B=\sum_{n\in\omega}\A(\ss\ric n+1)\N$, the element
$y = 0$ such that 
    $\Hyp{1}{0}{B} = \bf{1} = \ens{y}$, 
it is the case that for all n,
    $\cmp{y}{n}{B}
    =\sum_{i<n}\A(\ss\ric i+1)\N
    = \A(\ss\ric n\tac{0})$.
So $A \isom B$.

Thus, for the element $x\in A$ giving rise to the embedding
sequence $\ps$,
\dispM{
    \f_xA 
    \isom \sum_{0\le n<\rk{A}}A_n\N
    \isom \sum_{n\in\omega}\A(\ss\ric n+1)\N
    \isom \A(\ss^-\tac{0},0^{\omega}).
\Qed
}    

\parSkip
\Title{Corollary}
\Subtitle{Rank $\omega$ tail representation lemma, mark two}
The following expressions range precisely over the set of all
uniform tails of rank $\omega$.

{(a)}
For $\ss\in\SS(\omega+1)$ set 
\dispM{
    \A(\ss) := \A(\ss, 0^{\omega}).
}

{(b)}
For $\ss\in\SS(\omega)$ set 
\dispM{
    \A(\ss) := \A(\ss\tac{0}),
}
the latter being defined via (a) -- i.e.,
\dispm{
    \AA(\ss) = \A(\ss\tac{0}) := \A(\ss\tac{0}, 0^{\omega}).
}
\par


\chapSkip
\Chapter{[6] Chains of next limit rank}

\secSkip\noindent
The plus finite rank chain constructor of section 3 combined
with the rank $\omega$ representation of section 5 allows us
to obtain approximations to (and a characterisation of) uniform 
chains of rank $\omega+\omega$, but the limitations of this
procedure in having to step through all successor ordinals
will also become apparent.

\secSkip
\Section{[6.1]}
\Title{Chains of rank $\omega+\omega$}

\parSkip\noindent
Given $A\in\UC(\omega+\omega)$, applying the rank $\omega$ 
representation twice, to $A/\omega$ and to the rank $\omega$ 
building blocks of the rank $\omega+n$ sub-components 
    $\cmp{x}{\omega+n}{A}$
(for some fixed $x\in A$), one obtains a crude
representation involving 'arrays' of compatible pairs:
By the $\omega$+finite uniformity lemma the chain constructed as
    $\A(s)\arr{\M_0,\ldots,\M_{\lf(s)-1}}$
is uniform (of rank $\omega+\sz{s}$) if for 
    $\ss_i\in\SS(\omega+1)$ 
and 
    $\ps_i\in\Fun{\omega}{(\omega+1)}$
such that
    $\M_i = \A(\ss_i,\ps_i)$,
the $\omega$-uniformity conditions (3.6) hold: 
\dispM{
    \eqalign{
    &
    \M_i\in\UC(\omega) \hbox{ for all } i<\lf(s)
    \cr 
    &
    \f\M_i\isom\f\M_j \hbox{ for all } i,j<\lf(s)
    \cr 
    &
    \M_i\isom\M_j \hbox{ whenever } \wt^s(i) = \wt^s(j).
    \cr
    }
}
Combining this with the rank $\omega$ representation lemma it 
follows that 
    $\A(s)\arr{\M_0,\ldots,\M_{\lf(s)-1}}$
is uniform if the following derived $\omega$-uniformity conditions 
are satisfied:
\dispM{
    \eqalign{
    &
    \tup{\ss_i,\ps_i}\in\CSE(\omega+1) \hbox{ for all } i<\lf(s)
    \cr 
    &
    \ss_i^- = \ss_j^- \hbox{ for all } i,j<\lf(s)
    \cr 
    &
    \ss_i = \ss_j \hbox{ and } \ps_i\aeq\ps_j
    \hbox{ whenever } \wt^s(i) = \wt^s(j).
    \cr
    }\eqno(6.1)
}

Letting $A\in\UC(\omega+\omega)$, by the rank $\omega$ 
representation lemma there are compatible
    $\ss\in\SS(\omega+1)$ 
and,
    $\ps\in\Fun{\omega}{(\omega+1)}$
such that 
    $A/{\omega} = \A(\ss,\ps)$.
Again, by the rank $\omega$ representation lemma, for
some $x\in A$ and each $n$,
\dispM{
    \cmp{x}{\omega+n}{A}
    \isom \A(s_n)\arr{\M^{(n)}_i}_{0\le i<\lfs_n} 
}
where
\dispM{
    \M^{(n)}_i = \A(\ss^{(n)}_i,\ps^{(n)}_i),\eqno(6.2)
}
such that
\item{(i)}
$\ens{s_n}_{n\in\omega}$ derives from $\tup{\ss,\ps}$,
\item{(ii)}
the $\M^{(n)}_i\in\UC(\omega)$ satisfy the $\omega$-uniformity
conditions, and hence 
\item{(iii)}
the $\tup{\ss^{(n)}_i,\ps^{(n)}_i}$ satisfy the derived 
conditions. 

Consequently, there are $\ss^{(n)}\in\SS(\omega)$ such that
\dispM{
    \eqalign{
    &\ss^{(n)-}_i = \ss^{(n)} \hbox{ for all } i<\lfs_n, 
    \hbox{ and } \cr
    &\ps^{(n)}_i \aeq\ps^{(n)}_j \hbox{ whenever }
    \wts_n(i) = \wts_n(j)
    }
}
(where $\wts_n := \wt^{s_n}$). 

\parSkip
We shall need compact compact for the leaf intervals of the
$\Tree(s_{n+1}^{\pm})$ subtrees of $\Tree(s_{n+1})$ and the
corresponding intervals of building blocks. 

\parSkip
\Title{Notation}
Set
    $\Lf(s) = \lf(s) = [0,\lf(s))$, 
and 
    $\Lf(s)^+ = \Lf(s^+) 
    $,
but set
    $\Lf(s)^- = [\lf(s)-\lf(s^-),\lf(s))
    $.
Thus
\dispM{
    \Lf(s)^-
    = \cases{
       [\lf(s^+), \lf(s)) & if $s^+\in\SS$ \cr
       [0, \lf(s)) = \Lf(s) = \Lf(s^-) & else, \cr
   }
}
If we also convene that $\lf(s) = 0$ whenever $s\notin\SS$, then 
always
    $\Lf(s)^- = [\lf(s^+), \lf(s))$.

Set
    $\Lfs_n 
    = \Lf(s_n)$, 
and
    $\Lfs_n^{\pm} 
    = \Lf(s_n)^{\pm}$i
when dealing with $\ens{s_n}_n$ deriving from
    $\tup{\ss,\ps}\in\CSE(\omega+1)$.
%
%
Note, 
    $\Lfs_n^- 
    = \Lf(\ss\ric{n+1})^-$,
    $\Lfs_{n+1}^+ = \Lfs_n$ and
    $\Lfs_{n+1}^- = [\lfs_n,\lfs_{n+1})$ if $\tp_n = 0$,
    $\Lfs_{n+1}^+ = [0, \lfs_{n+1}-\lfs_n)$ and
    $\Lfs_{n+1}^- = [\lfs_{n+1}^+, \lfs_{n+1})$ if $\tp_n = 1$.

\def\SubLfs{{\cal J}}

\parSkip
Set 
    $\MM_n 
    = \MM^{(n)} 
    = \MM_{\Lfs_n} 
    = \tup{\M^{(n)}_i:i\in\Lfs_n}$
and similarly,
    $\MM_n^{\pm} 
    = \MM^{\pm}_{\Lfs_n} 
    = \MM_{\Lfs_n^{\pm}} 
    = \tup{\M^{(n)}_i:i\in\Lfs_n^{\pm}}$,
and for $\SubLfs\subeq\Lfs_n$,
    $\MM^{(n)}_{\SubLfs} 
    = \tup{\M^{(n)}_i:i\in\SubLfs}$,
etc.
    
\parSkip
Let now $x\in A\in\UC(\omega+\omega)$ such that for each $n$
setting
    $\M^{(n)}_i=\A(\ss^{(n)}_i,\ps^{(n)}_i)$.
Then
\dispm{
    \cmp{x}{\omega+n}{A}
    = \A(s_n)\arr{\MM_n},
}
and by plus finite chain construction, 
\dispM{
    \A(s_{n+1})\arr{\MM_{n+1}}
    = \cases{
      \A(s_{n+1}^-)\arr{\MM_{n+1}^-}.\Z
      & if $\tp_n=\omega$ \cr
      \A(s_{n+1}^+)\arr{\MM_{n+1}^+}
    + \A(s_{n+1}^-)\arr{\MM_{n+1}^-}.\N
      & else, \cr
    }
}
so:

\parSkip
\Number{(6.1)}
\Title{Lemma}
\Subtitle{Copy lemma}
Either 
\dispM{
    \hbox{$\tp_n = 0$, $s_{n}=s_{n+1}^+$, and }
    \A(s_n)\arr{\MM_n}
    = \A(s_{n+1}^+)\arr{\MM_{n+1}^+},
}
or 
\dispM{
    \eqalign{
    \hbox{$\tp_n > 0$, $s_{n}=s_{n+1}^- = \ss\nric{n+1}$, and }
    \A(s_{n})
    \arr{\MM_n}
    &= \A(s_{n+1}^-)
       \arr{\MM_{n+1}^-}, \cr
    }
}
and in particular,
\dispM{
    \MM_{n+1} = \cases{
    \MM_n
    & if $s_{n+1}^+\not\in\SS$ \cr 
    \MM_n
    \wcat\MM_{n+1}^-
    & if $s_{n+1}^+\in\SS$ and $s_{n+1}^+ = s_n$ \cr
    \MM_{n+1}^+
    \wcat\MM_n
    & if $s_{n+1}^+\in\SS$ and $s_{n+1}^- = s_n$ \cr}
}
i.e.,
\dispM{
    \MM_n
    = \cases{
    \MM_{n+1}^- & if $\tp_n > 0$ \cr
    \MM_{n+1}^+ & if $\tp_n = 0$ \cr
    }
}
(i.e., the '$s_n$-array'
     $\MM_n$
gets copied as one of the segments $\Lfs_{n+1}^+$ or $\Lfs_{n+1}^-$
of the '$s_{n+1}$ array' 
    $\MM_{n+1}$.)

\parSkip
From 
\dispM{
    \f\M_i^{(n)}\isom\f\M_0^{(n)}
                 \isom\f\M_0^{(n+1)}
                 \isom\f\M_j^{(n+1)}
}
for all $i<\lfs_n$ and $j<\lfs_{n+1}$ it follows that in fact
there is a single 'trunk' $\ss'\in\SS(\omega)$ such that for all
$n\in\omega$ and $i<\lfs_n$,
\dispM{
    \ss^{(n)-}_i = \ss^{(n)}_i\ric\omega = \ss'\ric{n}.
}
Thus we get:

\vfill\eject
\Number{(6.2)}
\Title{Lemma}

[1]
For every $A\in\UC(\omega+\omega)$,
\item{(i)}
    there exists $\ss'\in\SS(\omega)$ such that
    for all $n\in\omega,i<\lfs_n$, 
    $\ss^{(n)-}_i = \ss'\ric{n}$,
    and $\ss^{(n)}_i = \ss^{(n)}_j$
    whenever 
        $\wts_n(i) = \wts_n(j)$;
\item{(ii)}
    for all $n\in\omega$, 
    for all $i<\lfs_n$,
    $\ps^{(n)}_i$ is compatible with $\ss'$, 
    and for all $n,m\in\omega$, 
    for all $i<\lfs_n,j<\lfs_m$,
    $\ps^{(n)}_i\aeq\ps^{(m)}_j$; 
\item{(iii)}
    for all $n\in\omega$, 
    if $\ps_{n} = 0$ then
    $s_{n}=\ss_{n+1}^+$ and
    $\tup{\ps^{(n)}_i}_{0\le i<\lfs_n}
    = \tup{\ps^{(n+1)}_i}_{0\le i<\lfs_{n+1}^+}$; and
\item{(iv)}
    for all $n\in\omega$, 
    if $\ps_{n} > 0$ then
    $s_{n}=\ss_{n+1}^-$ and 
    $\tup{\ps^{(n)}_i}_{0\le i<\lfs_n}
    = \tup{\ps^{(n+1)}_i}_{\lfs_{n+1}^+\le i<\lfs_{n+1}}$.
\par\noindent
Moreover, the set of
    $\tup{\ss,\ps}\in\CSE(\omega+1)$,
    $\tup{\ss^{(n)},\ps^{(n)}}\in\CSE({\omega+1})$
satisfying (i) through (iv) is uniquely determined by $A$.

[2]
Every $A\in\UC(\omega+\omega)$ is uniquely determined by set of
    $\tup{\ss,\ps}\in\CSE(\omega+1)$,
    $\tup{\ss^{(n)},\ps^{(n)}}\in\CSE({\omega+1})$
satisfying (i) through (iv), and can be reconstructed from this set
by applying the rank $\omega$ representation lemma to each 
$\M^{(n)}_i$ and then to $A/\omega$.
\Qed

\parSkip
\Title{Corollary}
\Subtitle{Crude rank $\omega+\omega$ representation lemma}
Let 
\dispM{
    \A(\ss,\ps)\arr{\MM_n}_{n\in\omega} =
    \A(\ss,\ps,\tup{\ss^{(n)}_i,\ps^{(n)}_i}_{
                    n\in\omega,i\in\Lfs_n})
    \eqno(6.3)
}
denote the chain described in the $\omega+\omega$ representation 
lemma, i.e., where
    $\tup{\ss,\ps}, 
    \tup{\ss^{(n)}_i,\ps^{(n)}_i}
    \in\CSE(\omega+1)$
    ($n\in\omega,i\in\Lfs_n$)
satisfy conditions (i)-(iv) of the preceding lemma, and
where
    $\M^{(n)}_i=\A(\ss^{(n)}_i,\ps^{(n)}_i)$. 
Then, every 
    $A\in\UC(\omega+\omega)$ 
is of either of the forms (6.3), for some
    $\tup{\ss,\ps}, 
    \tup{\ss^{(n)}_i,\ps^{(n)}_i}
    \in\CSE(\omega+1)$
    ($n\in\omega,i\in\Lfs_n$),
or some array $\tup{\MM_n}_{n\in\omega}$ satsfying the uniformity 
conditions (6.1), respectively.

\parSkip
This representation lemma simply summarizes what we get 
on the basis of rank $\omega$ representation and rank 
$\omega$+finite approximation alone.
Removing the redundancy in the $\MM_n$ arrays will allow 
us to stream-line this 'representation' considerably.

\secSkip
\Section{[6.2]}
\Title{Tails of rank $\omega+\omega$}

\parSkip\noindent
To compute the tails of the components 
    $\cmp{x}{\omega+n}{A}$
and of $A$ itself we use the 'usual' strategy (basing ourselves 
on the uniformity over limits lemma, but going back to the tail
lemma for the calculations): 
\dispM{
   \f\cmp{x}{\omega+n}{A}
   = \f_x\cmp{x}{\omega+n}{A}
   = \f_x\cmp{x}{\omega}{A}
     +\sum_{i<n}\Suc^A_{\omega+i}.\N
}
where $\Suc^A_{\alpha}$ is the (isomorphism type) of the sucessor 
component of $\cmp{x}{\alpha}{A}$ in component $\cmp{x}{\alpha+1}{A}$.
Since
\dispM{
    \eqalign{
    \cmp{x}{\omega+n+1}{A}
    &
    = \A(s_{n+1})
      \arr{\MM_{n+1}} 
    \cr
    &
    = \cases{
      \A(s_{n+1}^-)
        \arr{\MM_{n+1}^-}.\Z, 
        & if $s_{n+1}^+\not\in\SS$
        \cr
      \A(s_{n+1}^+)
        \arr{\MM_{n+1}^+}
        + \A(s_{n+1}^-)
          \arr{\MM_{n+1}^-}.\N,
        & else
        \cr
        }
    \cr
    }
}
it follows that the successor $\omega+n$ component of
$\cmp{x}{\omega+n}{A}$ in $\cmp{x}{\omega+n+1}{A}$ is
\dispM{
    \Suc^A_{\omega+n}
    \isom \A(s_{n+1}^-)
       \arr{\MM_{n+1}^-} 
     \isom \A(\ss\nric{n+1})
       \arr{\MM_{n+1}^-}. 
}
(And $\cmp{x}{\omega}{A}=\M^{(0)}_0=\A(\ss^{(0)},\ps^{(0)})$.)
Thus (applying the $\omega$ tail lemma for the initial component),
\dispM{
    \eqalign{
    \f_x A
    &= \f_x\A(\tup{0})\arr{\A(\ss^{(0)}_0,\ps^{(0)}_0)}
    + \sum_{n\in\omega} A_{\omega+n}.\N \cr
    &= \f\A(\ss^{(0)}_0,\ps^{(0)}_0)
    + \sum_{n\in\omega} A_{\omega+n}.\N \cr
    &= \A(\ss'\tac{0})
     + \sum_{n\in\omega}
       \A(\ss\nric{n+1})
       \arr{\MM_{n+1}^-}.\N \cr
    }
}

\parSkip
Thus, for
    $x\in A
    = \A(\ss,\ps)
      \arr{\MM_n}$
we get
\dispM{
    \eqalign{
    \f_xA
    &= \f_x\cmp{x}{\omega}{A}
     + \sum_{n\in\omega}A_{\omega+n+1}\N \cr
    &= \f \M^{(0)}_0
     + \sum_{n\in\omega}\A(s_{n+1}^{(0)-})
       \arr{\MM_{n+1}^-}.\N \cr
    &= \sum_{n\in\omega} \A(\ss'\nric{n+1})\N
     + \sum_{n\in\omega}\A(\ss\nric{n+1})
       \arr{\MM_{n+1}^-}.\N \cr
    }
}
where $\ens{s_{n+1}^0}_{n\in\omega}$ derives from $\ss^{(0)}_0$,
and where $\ss'$ is the common initial $\omega$ segment of all the 
$\ss^{(n)}_i$ (i.e., $\ss'=\ss^{(n)}_i\ric\omega$, for all $n,i$).

\parSkip
\Number{(6.3)}
\Title{Lemma}
\Subtitle{Crude rank $\omega+\omega$ tail lemma}
For every $A\in\UC(\omega+\omega)$, there are 
    $\ss\in\SS(\omega+1)$
and
    $\tup{\ss^{(n)},\ps^{(n)}}\in\CSE(\omega+1)$,
such that
\dispM{
    \f A
    = \sum_{n\in\omega} \A(\ss^{(0)}\nric{n+1})\N
     + \sum_{n\in\omega}\A(\ss\nric{n+1})
       \arr{\MM_{n+1}^-}\N,
}
and where 
    $\ss^{(0)}\ric\omega = \ss^{(n)}\ric\omega$ for all $n$.
\Qed

\parSkip
Given $\tup{\ss,\ps}\in\CSE(\omega+1)$ and the derived sequence
$\ens{s_n}_n$, define sequence $\ens{k_n}_{n\in\omega}$ by
    $k_0 = 0$ 
and
    $k_{n+1} = k_n + \lfs_{n+1} - \lfs_n$.
Set $k^{\ss,\ps} = \sup_{n}k_n$.

\parSkip
Given a double array
    $\MM_{n+1}$
deriving from a chain $A\in\UC(\omega+\omega)$
define the sequence $\tup{\V_i}_{i<k^{\ss,\ps}}$ by
    $\V_0 = \M_0^{(0)}$
and supposing 
    $\ens{\V_i}_{i<k_n}$
defined,
\dispM{
    \eqalign{
    & \VV_{n+1} = \VV_{[k_n,k_{n+1})}
    = \tup{\V_i}_{k_n\le i<k_{n+1}} = \cases{
        \root
        & 
        if $\tp_n = \omega$, 
        \cr
        \MM_{n+1}^-
        & 
        if $\tp_n = 0$ 
        ($s_n = s_{n+1}^+$), 
        \cr
        \MM_{n+1}^+
        & 
        if $\tp_n = 1$ 
        ($s_n = s_{n+1}^-$) 
        \cr
        }
    }
}
I.e., the $\V_i$ solve the following set of equations:
\dispM{
    \eqalign{
    & \V_0 = \M_0^{(0)}, \cr
    & \MM_{n+1} = \cases{
        \MM_n 
        \wcat\VV_{n+1} 
        & 
        if $\tp_n = 0$, 
        \cr
        \VV_{n+1} 
        \wcat\MM_n
        & 
        if $\tp_n = 1$.  
        \cr
        }
    }
}
The $\V_i$, $i\in [k_n,k_{n+1})$ contain the bulding-blocks
that are not already contained in the arrays $\MM_k,k<n$.

\parSkip
Tails of rank $\omega+\omega$ represent a novel challenge:
When in the rank $\omega+(<\omega)$ approximation the lower 
rank component gets embedded with $\tp_n=1$ (i.e., there is 
an initial rank $\omega+n$ component in the rank $\omega+n+1$
approximation not necessarily isomorphic to the rank $\omega+n$ 
approximation, the building-blocks of the former are not all 
'visible' in the tail:

\parSkip
\Number{(6.4)}
\Title{Corollary}
\Subtitle{Tail revelation lemma}
\dispM{
    \f A = \f\A(\ss^{(0)},\ps^{(0)}) 
         + \sum_{n\in\omega}\Suc^{A}_{\omega+n}.\N,
}
where the (isomorphism type of the) successor $\omega+n$ component 
    $\Suc^{A}_{\omega+n}$
of 
    $A
    =\A(\ss,\ps,\tup{\MM_n:n\in\omega})
       \in\UC(\omega+\omega)$,
satisfies
\dispM{
    \Suc^{A}_{\omega+n} 
    = \A(s_{n+1}^-)
      \arr{\MM_{n+1}^-}
    = \cases{
        \A(s_{n+1}^-)
        \arr{\VV_{n+1}}
        & if $\tp_n = 0$, \cr 
        \A(s_{n+1}^-)
        \arr{\MM_n}
        & if $\tp_n \not=0$. \cr 
        }
}

Where 
    $\Suc^{A}_{\omega+n} = \A(s_n)\arr{\MM_n}$,
the 'new' building blocks at level $n+1$ are 'hidden' from the tail, 
but they will be 'revealed' at the next higher level $m>n$ where
    $\Suc^{A}_{\omega+m} = \A(s_n)\arr{\MM_m}$
 --- because
    $\VV_{n+1}$ is a sub-interval of building blocks of $\MM_m$.
As a consequence, there can be at most one level introducing new
building blocks that may stay 'hidden' from the entire approximation
to the tail, namely if for some largest $n$, $\tp_{\omega+n}=1$
and for all $m>n$, $\tp_{\omega+m}=1$. Anticipating that at any
level there is in fact atmost one new building block, it follows 
that at most one such can stay permanently 'hidden' from the 
tail. 

\parSkip
It's now time to prepare the grounds for a managable representation.

\secSkip
\Section{[6.3]}
\Title{Array reduction arithmetic}

\parSkip\noindent
Set 
    $\rwts_n = \rwt^{s_n}$ 
(where
    $\rwt^s = \ran(\wt^s)$).
By the $\omega+$finite uniformity lemma, there is a sequence of
sets
    $\ens{\U^{(n)}_i:i\in\rwts_n} \sub \UC(\omega)$
    ($n\in\omega$),
such that for each $n$
\dispM{
    (\forall i\in\lfs_{n})
     [\M^{(n)}_i=\U^{(n)}_{\wts_n(i)}].
}
%
Let 
    $\UU_n=\ens{\tup{\U^{(n)}_k,k}:k\in\rwts_n}$
be the tagged short array for
    $\cmp{x}{\omega+n}{A}$,
i.e., such that
    $\cmp{x}{\omega+n}{A}
    = \A(s_n)\rarr{\UU_n}$.
By the copy lemma 
\dispM{
    \UU_n = \cases{(\UU_{n+1})^+_{s_{n+1}} & if $\tp_n=0$ \cr
                   (\UU_{n+1})^-_{s_{n+1}} & if $\tp_n>0$ \cr}
}

We'd like to obtain a short array for 
    $A\in\UC(\omega+\omega)$
as a sequence
    $\UU
    =\ens{\U_k:k\in K}$
(for some set $K$) that would allow to 'generate' the $\UU_n$
perhaps resembling the way in which the $s_n$ were generated
from $\ss\in\SS(\omega+1)$ (and compatible $\tp$). 
%
%
%
%
We shall first show that we can in fact obtain a 'canonically 
ordered' short array, and subsequently how to generate 
short arrays for the sub-components of $\A$ in the 
approximating chains
    $\cmp{x}{\omega+n}{A}$.

\parSkip
Set 
    $\mwt^s = \max_{i<\lf(s)}\wt^s(i)$ 
(and again, 
    $\mwts_n = \mwt^{s_n}$).
Define the {\it augmented weight sequence} by:
\dispM{
    \xwt^{s} 
    = \tup{\wt^{s}(0)+1}\cat\tup{\wt^{s}(i):0<i<\lf(s)}
    = \tup{\wt^{s}(0)+1}\wcat\wt^{s}\ric[1,\lf(s)).
}
(Also let 
    $\rxwt^s=\rn\xwt^s$; 
note that 
    $\sz{\rxwt^s} = \sz{\rwt^s} = \kwt^s$.)

From the corollaries to the weight induction lemma, 3.9,
we have:

\parSkip
\Number{(6.5)}
\Title{Lemma}
\Subtitle{Weight sequence lemma, maximum weight lemma}
For $s\in\SS$,

(a)
\dispM{
    \wt^s = \cases{
    \wt^{s^-}
    & if $s^+\not\in\SS$ \cr 
    \wt^{s^+}\wcat\wt^{s^-}
    & if $s^+\in\SS$ and $s\not=s^-\tac{0}$ \cr 
    \xwt^{s^+}\wcat\wt^{s^-}
    & if $s^+\in\SS$ and $s=s^-\tac{0}$ \cr 
    }
}

(b)
\dispM{
    \mwt^s = \cases{
    \mwt^{s^-} 
    & if $s^+\in\SS$ and $\mwt^{s^-} > \mwt^{s^+}$,
      or $s^+\not\in\SS$ \cr
    \mwt^{s^+}
    & if $s^+\in\SS$ and $\mwt^{s^-} \le \mwt^{s^+}$
      and $s\not=s^-\tac{0}$ \cr
    \mwt^{s^+}+1
    & if $s^+\in\SS$ and $\mwt^{s^-} \le \mwt^{s^+}$
      and $s=s^-\tac{0}$. \cr
    }
}

(c)
\dispM{
    \eqalign{
    &
    \mwt^s\le\sz{s}-1
    \qquad\hbox{i.e., } 
    \cr
    &
    \mwt^s\le n \hbox{ whenever } s\in\SS(n+1)
    \cr
    }
}

(d)
\dispM{
    \eqalign{
    &
    \mwt^{s\tac{0}}=\wt^{s\tac{0}}(0)=\sz{s}
    \cr
    &
    \mwt^{s\tac{b}},\wt^{s\tac{b}}(0)<\sz{s}\hbox{ if }0<b\le\sz{s}
    \cr
    }
}
i.e.,
\dispM{
    \eqalign{
    &
    \mwt^s=\wt^{s}(0)=\sz{s}-1 \hbox{ if } s=s^-\tac{0}
    \cr
    &
    \mwt^s,\wt^{s}(0)<\sz{s}-1 \hbox{ if } s=s^-\tac{b}
    \hbox{ for } 0< b\le\sz{s}
    \cr
    }
}

\Proof (a) and (b) follow from 3.9, (c) and (d) follow by induction.
%
%
%
%
%
\Qed

\parSkip
The maximum weight lemma almost seems to suggest that
    $\mwts_{n+1}\in\ens{\mwts_n,\mwts_n+1}$.
It may come as somewhat of a surprise that this does not hold. 
In fact, as the examples in preceding sections suggest quite
the contrary.
%
%
This is because the weight-incrementation of the branching 
nodes proceeds rather indiscriminately. Consider the example 
$\Tree(\tup{0,1,0})$,
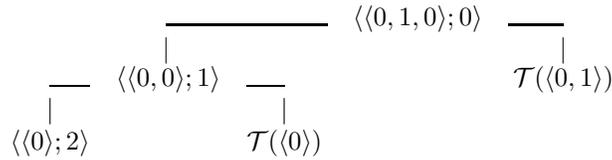
\dispM{
\mrx{ \tup{\tup{0,1,0};0} }
    { \mrx{ \tup{\tup{0,0};1} }
          { \tup{\tup{0};2} }
          { \Tree(\tup{0}) } }
    { \trx{ \Tree(\tup{0,1}) } }
}
where it would have sufficed to 'copy', rather than increment,
the weight from node $\tup{0,0}$ to node $\tup{0}$, in order
to distingish it in the weight-sequence ($\tup{1,0,0}$, rather
than $\tup{2,0,0}$) from the two weight 0 nodes.

\parSkip
Nevertheless, something of this earlier intuition can be
salvaged; we do have that at every stage at most one new 
'weight-element' is introduced --- and consequently at most
one new rank $\omega$ 'building block' in the array of rank
$\omega$ components. More precisely, what we will show is
that
\dispM{
    0 \le \kwt_{n+1} - \kwt_n \le 1,
}
where
\dispm{
    \kwt^s=\sz{\rwt^s}
}
 and $\kwt_n = \kwt^{s_n}$.
(Thus at every rank $\omega+n$ at most one 'new' type of $\omega$ 
component appears in the array of 'building blocks'.)
%

\parSkip
\Number{(6.6)}
\Title{Lemma}
\Subtitle{Plus one range lemma}
For $s\in\SS$ and $0<b\le \sz{s}$,
\dispM{
    \rwt^{s\tac{0}}\diff\rwt^{s\tac{b}}
    = \rwt^{s\tac{0}}\rid\rwt^{s\tac{b}}
    = \ens{\wt^{s\tac{0}}(0)}.
}

\Proof
By induction on $\sz{s}$. 

$\sz{s} = 1$: Then $s\tac{0} = \tup{0,0}$ and the only possiblility
for $b$ is $b=1=\sz{s}$, so $s\tac{b}=\tup{0,0}$. Then
    $\wt^{s\tac{0}} = \tup{1,0}$ 
and 
    $\wt^{s\tac{b}}=\wt^s=\tup{0}$,
and the conclusion holds. Since this only covers the case $b=\sz{s}$,
we give the case $\sz{s}=2$ as well.

$\sz{s} = 2$: There are two cases for $s$ and for each of these two
cases for $b$:

$s=\tup{0,0}$: Then $s\tac{0}=\tup{0,0,0}$, so
    $\wt^{s\tac{0}} = \wt^{\tup{0,0,0}} = \tup{2,0,1,0}$.
For $b=1$, 
    $\wt^{s\tac{1}} = \wt^{\tup{0,0,1}} = \tup{0,1,0}$,
so
\dispM{
    \rwt^{s\tac{0}}\diff\rwt^{s\tac{b}}
    = \rwt^{s\tac{0}}\rid\rwt^{s\tac{b}}
    = \ens{0,1,2} \rid \ens{0,1}
    = \ens{2}
    = \ens{\wt^{s\tac{0}}(0)},
}
as required.
For $b=2$, 
    $\wt^{s\tac{2}} = \wt^{\tup{0,0,2}}
    = \wt^{\tup{0,0}} = \tup{1,0}$,
and the same calculation applies as before.

$s=\tup{0,1}$: Then $s\tac{0}=\tup{0,1,0}$, so
    $\wt^{s\tac{0}} = \wt^{\tup{0,1,0}} = \tup{2,0,0}$.
For $b=1$, 
    $\wt^{s\tac{1}} = \wt^{\tup{0,1,1}} = \tup{0,0}$,
so
\dispM{
    \rwt^{s\tac{0}}\diff\rwt^{s\tac{b}}
    = \rwt^{s\tac{0}}\rid\rwt^{s\tac{b}}
    = \ens{0,2} \rid \ens{0}
    = \ens{2}
    = \ens{\wt^{s\tac{0}}(0)},
}
as required.
For $b=2$, 
    $\wt^{s\tac{2}} = \wt^{\tup{0,1,2}}
    = \wt^{\tup{0}} = \tup{0}$,
and the calculation is as before.


For the induction step, distinguish cases as follows:

{\bf Case} $0<b<\sz{s}$. 
Then
\dispM{
	\wt^{s\tac{0}} 
    = \xwt^{s\tac{0}^+}\wcat\wt^{s\tac{0}^-}
    = \xwt^{s^-\tac{0}}\wcat\wt^{s}.
}
Since $0<b<\sz{s}$, also
\dispM{
	\wt^{s\tac{b}} 
    = \wt^{s\tac{b}^+}\wcat\wt^{s\tac{b}^-}
    = \wt^{s^-\tac{b}}\wcat\wt^s.
}
The respective 'second' half segments being identical, 
and since 
    $\wt^{s\tac{0}}(0) = \wt^{s^-\tac{0}}(0)+1$,
\dispM{
    \eqalign{
    \rwt^{s\tac{0}}\diff\rwt^{s\tac{b}} 
    &= \rxwt^{s^-\tac{0}}\diff\rwt^{s^-\tac{b}} 
    \cr
    &= \rwt^{s^-\tac{0}}\rid\ens{\wt^{s^-\tac{0}}(0)}
                       \cup\ens{\wt^{s^-\tac{0}}(0)+1})
                       \diff\rwt^{s^-\tac{b}}. 
    \cr
    }
}
Now (e.g., by the weight sequence lemma, (c),
    $\wt^{s^-\tac{0}}(0)+1\not\in\rwt^{s^-\tac{b}}$.
Observe that for sets $A,B$ and elements $a\in A\rid B$ and 
$\tilde{a}\not\in B$, 
\dispM{
    (A\rid \ens{a}\cup\ens{\tilde{a}})\rid B 
    = (A\rid B)\rid\ens{a}\cup\ens{\tilde{a}}, 
}
and if also 
    $B\rid A = \emptyset$,
then also
    $B\rid (A\rid \ens{a}\cup\ens{\tilde{a}}) = \emptyset$, 
and so
\dispM{
    (A\rid \ens{a}\cup\ens{\tilde{a}})\diff B 
    = (A\rid B)\rid\ens{a}\cup\ens{\tilde{a}}. 
}
Thus, by inductive hypothesis,
\dispM{
    \eqalign{
    &\rwt^{s^-\tac{0}}\rid\ens{\wt^{s^-\tac{0}}(0)}
                      \cup\ens{\wt^{s^-\tac{0}}(0)+1})
                      \diff\rwt^{s^-\tac{b}} 
    \cr
    &\qquad
    = \rwt^{s^-\tac{0}}\rid\rwt^{s^-\tac{b}}
                       \rid\ens{\wt^{s^-\tac{0}}(0)}
                       \cup\ens{\wt^{s^-\tac{0}}(0)+1}
    \cr
    &\qquad
    = \rwt^{s^-\tac{0}}\diff\rwt^{s^-\tac{b}}
                       \rid\ens{\wt^{s^-\tac{0}}(0)}
                       \cup\ens{\wt^{s^-\tac{0}}(0)+1}
    \cr
    &\qquad
    = \ens{\wt^{s^-\tac{0}}(0)+1)}, 
    \cr
    &\qquad
    = \ens{\wt^{s\tac{0}}(0)}, 
    \cr
    }
}
as required.
    
{\bf Case} $0<b=\sz{s}$. 
Then $\wt^{s\tac{b}} = \wt^s$, so
\dispM{
    \eqalign{
    \rwt^{s\tac{0}}\diff\rwt^{s\tac{b}}
    &= \rxwt^{s^-\tac{0}}\diff\rwt^s 
    \cr
    &= \rwt^{s^-\tac{0}}\rid\ens{\wt^{s^-\tac{0}}(0)}
                       \cup\ens{\wt^{s^-\tac{0}}(0)+1})
                       \diff\rxwt^s 
    \cr
    &= \ens{\wt^{s^-\tac{0}}(0)+1)}, 
    \cr
    }
}
again, as required.
\Qed

\parSkip
\Number{(6.7)}
\Title{Corollary}
For $s\in\SS$ and $0 < a, b \le \sz{s}$,
\item{(i)}
\dispm{
    \rwt^{s\tac{a}}\diff\rwt^{s\tac{b}} = \emptyset.
}
\item{(ii)}
\dispm{
    \sz{\rxwt^{s\tac{0}}\rid\rwt^{s\tac{b}}}
    = \sz{\rxwt^{s\tac{0}}\diff\rwt^{s\tac{b}}}
    = \sz{\rwt^{s\tac{0}}\rid\rwt^{s\tac{b}}}
    = \sz{\rwt^{s\tac{0}}\diff\rwt^{s\tac{b}}}
    \le 1.
}   
\item{(iii)}
\dispm{
    0 \le \kwt^{s\tac{a,b}}-\kwt^{s\tac{a}}\le 1
}
and
\dispm{
    0 \le \kwt^{s\tac{a,b}}-\kwt^{s\tac{b}}\le 1
}

\Proof 
(i)
\dispm{
    \rwt^{s\tac{a}}\diff\rwt^{s\tac{b}} 
    = (\rwt^{s\tac{0}}\rid\ens{\wt^{s\tac{0}}(0)})
      \diff (\rwt^{s\tac{0}}\rid\ens{\wt^{s\tac{0}}(0)})
    = \emptyset.
}

(ii)
By the plus one range lemma,
\dispm{
    \rwt^{s\tac{0}}\rid\rwt^{s\tac{b}}
    = \ens{\wt^{s\tac{0}}(0)}
    = \ens{\mwt^{s\tac{0}}},
}
and by the preceding lemma $\mwt^{s\tac{0}} > \mwt^{s\tac{b}}$,
so
\dispm{
    \rxwt^{s\tac{0}}\rid\rwt^{s\tac{b}}
    = \ens{\xwt^{s\tac{0}}(0)}
    = \ens{\mwt^{s\tac{0}}+1}.
}

(iii)
Distinguish cases, and apply the plus one range lemma to
the ranges compared.

Case $a=b=0$:
\dispm{
    \kwt^{s\tac{a,b}}-\kwt^{s\tac{a}}
    =\sz{\rn(\xwt^{s\tac{0}}\cat\wt^{s\tac{0}})}
     -\kwt^{s\tac{a}}
    =1
}
and ditto
\dispm{
    \kwt^{s\tac{a,b}}-\kwt^{s\tac{b}}=1
}

Case $a>b=0$:
\dispm{
    \kwt^{s\tac{a,b}}-\kwt^{s\tac{a}}
    =\sz{\rn(\xwt^{s\tac{0}}\cat\wt^{s\tac{a}})}
     -\kwt^{s\tac{a}}
    =1
}
and
\dispm{
    \kwt^{s\tac{a,b}}-\kwt^{s\tac{a}}
    =\sz{\rn(\xwt^{s\tac{0}}\cat\wt^{s\tac{a}})}
     -\kwt^{s\tac{0}}
    =1
}

Case $a=0<b$:
\dispm{
    \kwt^{s\tac{a,b}}-\kwt^{s\tac{a}}
    =\sz{\rn(\wt^{s\tac{b}}\cat\wt^{s\tac{0}})}
     -\kwt^{s\tac{0}}
    =0
}
and
\dispm{
    \kwt^{s\tac{a,b}}-\kwt^{s\tac{a}}
    =\sz{\rn(\wt^{s\tac{b}}\cat\wt^{s\tac{0}})}
     -\kwt^{s\tac{b}}
    =1
}

Case $0<a,b$:
\dispm{
    \kwt^{s\tac{a,b}}-\kwt^{s\tac{a}}
    =\sz{\rn(\wt^{s\tac{b}}\cat\wt^{s\tac{a}})}
     -\kwt^{s\tac{a}}
    =0
}
and
\dispm{
    \kwt^{s\tac{a,b}}-\kwt^{s\tac{a}}
    =\sz{\rn(\wt^{s\tac{b}}\cat\wt^{s\tac{a}})}
     -\kwt^{s\tac{b}}
    =0.
\Qed
}

\parSkip
By the preceding lemma, 

\parSkip
\Title{Corollary}
For each $n$,
\item{(i)}
    $\kwt_{n+1}-\kwt_n\le 1$.
\item{(ii)}
    $\rwts_{n+1}$
is an end-extension of
    $\rwts_n$,
since either
    $\rwts_{n+1} = \rwts_n$,
or both
    $\rwts_{n+1}\supeq\rwts_n$
and
    $\min(\rwts_{n+1}\rid\rwts_n)
    =\mwts_{n+1}>\mwts_n$.

\parSkip
\Number{(6.8)}
\Title{Lemma}
\Subtitle{Plus minus lemma}
For all $s\in\SS(<\omega)$, where $s^+\in\SS$,
\dispM{
    \sz{\rwt^{s^+}\diff\rwt^{s^-}} \le 1
}

\Proof
By induction on $\Dom(s)$. Thus, suppose that for all $t\in
\Dom(s)$ with $t^+\in\SS$, 
    $\sz{\rwt^{t^+}\diff\rwt^{t^-}} \le 1$.
By weight sequence lemma, (a),
\dispM{
    \wt^{s^+} = \cases{
    \wt^{s^{+-}}
    & if $s^{++}\not\in\SS$ \cr 
    \wt^{s^{++}}\wcat\wt^{s^{+-}}
    & if $s^{++}\in\SS$ and $s^+\not=s^{--}\tac{0}$ \cr 
    \xwt^{s^{++}}\wcat\wt^{s^{+-}}
    & if $s^{++}\in\SS$ and $s=s^{--}\tac{0}$, \cr 
    }
}
and
\dispM{
    \wt^{s^-} = \cases{
    \wt^{s^{--}}
    & if $s^{-+}\not\in\SS$ \cr 
    \wt^{s^{-+}}\wcat\wt^{s^{--}}
    & if $s^{-+}\in\SS$ and $s^-\not=s^{--}\tac{0}$ \cr 
    \xwt^{s^{-+}}\wcat\wt^{s^{--}}
    & if $s^{-+}\in\SS$ and $s^-=s^{--}\tac{0}$. \cr 
    }
}

Let $s=t\tac{a,b}$ ($0\le a,b\le\sz{t}$). By the preceding,
\dispM{
    \wt^{s^+} = \wt^{t\tac{b}} = \cases{
    \wt^{t}
    & if $b=\sz{t}$ \cr 
    \wt^{t^-\tac{b}}\wcat\wt^{t}
    & if $0<b<\sz{t}$ \cr 
    \xwt^{t^-\tac{0}}\wcat\wt^{t}
    & if $b=0$, \cr 
    }
}
and
\dispM{
    \wt^{s^-} = \wt^{t\tac{a}} = \cases{
    \wt^{t}
    & if $a=\sz{t}$ \cr 
    \wt^{t^-\tac{a}}\wcat\wt^{t}
    & if $0<a<\sz{t}$ \cr 
    \xwt^{t^-\tac{0}}\wcat\wt^{t}
    & if $a=0$. \cr 
    }
}

We shall make heavy use of $s^{+-}=s^{--}$, 
which we refer to below as $t$. 
Write $A\equipoll B$ for $\sz{A} = \sz{B}$ and
$A\smalleq B$ for $\sz{A} \leq \sz{B}$, etc.

\parSkip
{\bf Case 1\hskip .2in} 
$s^{++}\not\in\SS$ (i.e., $b=\sz{t}$).

{\bf Sub-case 1.1\hskip .2in} 
$s^{-+}\not\in\SS$ (i.e., $a=\sz{t}$):
Then
\dispM{
    \rwt^{s^+}\diff\rwt^{s^-}
    = \rwt^{t}\diff\rwt^{t}
    = \emptyset.
}

{\bf Sub-case 1.2\hskip .2in} 
$s^{-+}\in\SS$ and $s^-\not=s^{--}\tac{0}$ (i.e., $0<a\le\sz{t}$):
\dispM{
    \eqalign{
    \rwt^{s^+}\diff\rwt^{s^-}
    &= \rwt^{t}\diff\rn(\wt^{t^-\tac{a}}\wcat\wt^{t})
    \cr
    &= \rwt^{t^-\tac{a}}\rid\rwt^{t}
    \smalleq 1,
    \cr
    }
}
by the induction hypothesis appied to $s^-=t\tac{a}$.

{\bf Sub-case 1.3\hskip .2in} 
$s^{-+}\in\SS$ and $s^-=s^{--}\tac{0}$ (i.e., $a=0$):
Then
\dispM{
    \eqalign{
    \rwt^{s^+}\diff\rwt^{s^-}
    & = \rwt^{t}\diff\rn(\xwt^{s^{-+}}\wcat\wt^{s^{--}})
    \cr
    & = \rwt^{t}\diff\rn(\tup{\wt^{t}(0)+1}
        \wcat\wt^{t}\ric[1,\lf(t))\wcat\wt^{t})
    \cr
    & = \ens{\wt^{t}(0)+1} \equipoll 1.
    \cr
    }
}

\parSkip
{\bf Case 2\hskip .2in} 
$s^{++}\in\SS$ and $s^+\not=s^{--}\tac{0}$ (i.e., $0<b<\sz{t}$).

{\bf Sub-case 2.1\hskip .2in} 
$s^{-+}\not\in\SS$ (i.e., $a=\sz{t}$):
\dispM{
    \eqalign{
    \rwt^{s^+}\diff\rwt^{s^-}
    &= \rn(\wt^{t^-\tac{b}}\wcat\wt^{t})\diff
       \rwt^{t}
    \cr
    &= \rwt^{t^-\tac{b}}\rid\rwt^{t}
    \smalleq 1
    \cr
    }
}
by induction hypothesis applied to $s^+=t\tac{b}$.

{\bf Sub-case 2.2\hskip .2in} 
$s^{-+}\in\SS$ and $s^-\not=s^{--}\tac{0}$ (i.e., $0<a<\sz{t}$):
By the preceding corollary,
\dispM{
    \eqalign{
    \rwt^{s^+}\diff\rwt^{s^-}
    &= \rn(\wt^{t^-\tac{b}}\wcat\wt^{t})\diff
       \rn(\wt^{t^-\tac{a}}\wcat\wt^{t})
    \cr
    &\subeq
       \rwt^{t^-\tac{b}}\diff
       \rwt^{t^-\tac{a}}
    = \emptyset,
    \cr
    }
}
using $(A\cup C) \diff (B\cup C) = (A\diff B)\rid C \subeq A\diff B$.

{\bf Sub-case 2.3\hskip .2in} 
$s^{-+}\in\SS$ and $s^-=s^{--}\tac{0}$ (i.e., $a=0$):
By the preceding lemma,
\dispM{
    \eqalign{
    \rwt^{s^+}\diff\rwt^{s^-}
    &= \rwt^{t^-\tac{0}}\diff
       \rwt^{t^-\tac{b}}
    \equipoll 1
    \cr
    }
}
by the lemma.

\parSkip
{\bf Case 3\hskip .2in} 
$s^{++}\in\SS$ and $s^+=s^{--}\tac{0}$ (i.e., $b=0$).

{\bf Sub-case 3.1\hskip .2in} 
$s^{-+}\not\in\SS$ (i.e., $a=\sz{t}$):
\dispM{
    \eqalign{
    \rwt^{s^+}\diff\rwt^{s^-}
    &= \rn(\xwt^{s^{++}}\cat\wt^{s^{+-}})\diff\rwt^{t}
    \cr
    &= \rn(\tup{\wt^{t^-\tac{0}}(0)+1}
       \wcat\wt^{t^-\tac{0}}\ric[1,\lf(t^-\tac{0})\wcat\wt^{t})\diff 
       \rwt^{t}
    \cr
    &= \rxwt^{t^-\tac{0}}\diff
       \rwt^{t}
    \equipoll 1
    \cr
    }
}
by part (ii) of the preceding corollary.

{\bf Sub-case 3.2\hskip .2in} 
$s^{-+}\in\SS$ and $s^-\not=s^{--}\tac{0}$:
This case is 'symmetric' to case 2.2

{\bf Sub-case 3.3\hskip .2in} 
$s^{-+}\in\SS$ and $s^-=s^{--}\tac{0}$:
Then $s^+=s^-$, and the result is immediate.
\Qed

\parSkip
\Number{(6.9)}
\Title{Corollary}
\Subtitle{Determination lemma}
Let $s = t\tac{a,b}\in\SS$. Then either 
\item{(i)}
$b =\sz{t}+1$ and ($s^+\not\in\SS$ and) 
    $\rwt^{s} = \rwt^{s^-}$; 
or
\item{(ii)}
$0 < a, b \leq \sz{t}$ and 
    $\rwt^{s^+} = \rwt^{s^-}$; 
or
\item{(iii)}
$a = 0 < b \leq \sz{t}$ and 
    $\rwt^{s^+}\sub\rn(\wt^{s^-}\ric[1,\lfs^-))$;
or
\item{(iv)}
$b = 0 < a \leq \sz{t}$ and 
    $\rwt^{s^-}\sub\rn(\wt^{s^+}\ric[1,\lfs^+))$,
hence also
    $\rwt^{s^-}\sub\rn(\xwt^{s^+}\ric[1,\lfs^+))$,
or
\item{(v)}
$0 = a = b$ and 
    $\rwt^{s^+} = \rwt^{s^-}$,
and hence 
    $\rxwt^{s^+}\ric[1,\lfs^+) = \rwt^{s^-}\ric[1,\lfs^-)$.

\Proof
(i) 
Follows since in this case $\lfs=\lfs^-$ and $\wt^s=\wt^{s^-}$.

(ii)
Follows by subcase 2.2 of the preceding lemma.

(iii)
By induction over $\sz{t}$:
For $t = \tup{0}$ the only case to consider is $s = \tup{0,0,1}$.
Then
\dispM{
    \eqalign{
    &\wt^{s^+} = \wt^{\tup{0,1}} = \wt^{\tup{0}} = \tup{0} \cr
    &\wt^{s^-} = \wt^{\tup{0,0}} = \tup{1,0} \cr
    }
}
so $\rwt^{s^+} = \ens{0}\sub\ens{0} = \rn(\wt^{s^-}\ric[1,2))$.
    
So suppose the result holds for all $u\in\SS(\le n)$ for $n\ge 1$ 
and let $t\in\SS(n+1)$. Then
\dispM{
    \eqalign{
        &\wt^{s^+} = \wt^{t\tac{b}}) = \cases{
        \wt^{t} & if $b = \sz{t}$ \cr
        \wt^{t^-\tac{b}}\cat\wt^{t} & else \cr} \cr
        &\wt^{s^-} = \wt^{t\tac{0}} = \xwt^{t^-\tac{0}}\cat\wt^{t} \cr
    }
}
By inductive hypothesis, 
    $\rwt^{t^-\tac{b}}\sub\rn(\wt^{t^-\tac{0}}\ric[1,\lf^{-+}))$,
where
    $\lf^{-+} = \lf(s^{-+}) = \lf(t^-\tac{0})$,
so the result follows.

(iv)
The argument is almost completely symmetric to the one for (iii)
and proceeds again by induction over $\sz{t}$.

For $t = \tup{0}$ the only case to consider is $s = \tup{0,1,0}$.
Then
\dispM{
    \eqalign{
    &\wt^{s^+} = \wt^{\tup{0,0}} = \tup{1,0} \cr
    &\wt^{s^-} = \wt^{\tup{0,1}} = \wt^{\tup{0}} = \tup{0} \cr
    }
}
so $\rwt^{s^-} = \ens{0}\sub\ens{0} = \rn(\wt^{s^+}\ric[1,2))$.
    
Supposing the result holds for all $u\in\SS(\le n)$ for $n\ge 1$ 
and let $t\in\SS(n+1)$. Then
\dispM{
    \eqalign{
        &\wt^{s^-} = \wt^{t\tac{a}}) = \cases{
        \wt^{t} & if $a = \sz{t}$ \cr
        \wt^{t^-\tac{a}}\cat\wt^{t} & else \cr} \cr
        &\wt^{s^+} = \wt^{t\tac{0}} = \xwt^{t^-\tac{0}}\cat\wt^{t} \cr
    }
}
By inductive hypothesis, 
    $\rwt^{t^-\tac{a}}\sub\rn(\wt^{t^-\tac{0}}\ric[1,\lf^{-+}))$,
where
    $\lf^{-+} = \lf(s^{-+}) = \lf(t^-\tac{0})$,
so the result follows.

(v)
Follows again from subcase 3.3 of the preceding lemma.
\Qed

\parSkip
Note that all this implies that $\rwt^{s\tac{a}}=\rwt^{s\tac{b}}$
iff $a=b=0$ or $0<a,b\leq\sz{s}$.

\parSkip
\Title{Corollary}
For $s=t\tac{a,b}\in\SS$ such that $s^+\in\SS$
one of the following applies:
\item{(i)}
if $a=b=0$, then 
    $\kwt^s=\kwt^{s^+}+1=\kwt^{s^-}+1$
and
    $\mwt^s=\mwt^{s^+}+1=\mwt^{s^-}+1$;
\item{(ii)}
if $a>0=b$, then 
    $\kwt^s=\kwt^{s^+}=\kwt^{s^-}+1$
and
    $\mwt^s=\mwt^{s^+}+1>\maxs{\mwt^{s^+},\kwt^{s^-}}$;
\item{(iii)}
if $a=0<b$, then 
    $\kwt^s=\kwt^{s^-}=\kwt^{s^+}+1$
and
    $\mwt^s=\mwt^{s^-}+1>\maxs{\mwt^{s^+},\kwt^{s^-}}$;
\item{(iv)}
if $a,b>0$, then 
    $\kwt^s=\kwt^{s^+}=\kwt^{s^-}$
and
    $\mwt^s=\mwt^{s^+}=\mwt^{s^-}$.

For $s=t\tac{a,b}\in\SS$ such that $s^+i\not\in\SS$
\item{(v)}
    $\kwt^s=\kwt^{s^-}$
and
    $\mwt^s=\mwt^{s^-}$.


\parSkip
\Number{(6.10)}
\Title{Corollary and Definition} 
\Subtitle{Free position, first form}
Let $s=t\tac{a,b}\in\SS(<\omega)$. Then one of the following
applies:

\item{(i)}
    $a=b=0$, so
    $s^+=s^-=t\tac{0}$, hence
\itemitem{}
    $\wt^s=\xwt^{t\tac{0}}\cat\wt^{t\tac{0}}$,
\itemitem{}
$\xwt^{t\tac{0}}\diff\wt^{t\tac{0}}=
    \ens{\wt^{t\tac{0}}(0)+1,\wt^{t\tac{0}}(0)}
    = \ens{\wt^s(0),\wt^{s^-}(0)}
    = \ens{\wt^s(0),\wt^{s^+}(0)}$.
\itemitem{}
Set $\free^s = \ens{0,\lf(s^+)}$.

\item{(ii)}
    $a=0<b<\sz{s}-1$, so
    $s^+=t\tac{b}$, $s^-=t\tac{0}$, hence
\itemitem{}
    $\wt^s=\wt^{t\tac{b}}\cat\wt^{t\tac{0}}$,
\itemitem{}
    $\wt^{t\tac{b}}\diff\wt^{t\tac{0}}
    = \ens{\wt^{t\tac{0}}(0)}=\ens{\wt^{s}(\lf(s^+))}$.
\itemitem{}
Set $\free^s = \ens{\lf(s^+)}$. 

\item{(iii)}
    $a=0<b=\sz{s}-1$, so
    $s^+\notin\SS$, hence
\itemitem{}
    $\wt^{s}=\wt^{s^-}$.
\itemitem{}
Set $\free^s = \void$.

\item{(iv)}
    $b=0<a\le\sz{s}-2$, so
    $s^+=t\tac{0}$, $s^-=t\tac{a}$, hence
\itemitem{}
    $\wt^{s}=\xwt^{t\tac{0}}\cat\wt^{t\tac{a}}$,
\itemitem{}
    $\xwt^{t\tac{0}}\diff\wt^{t\tac{a}}
    = \ens{\wt^{t\tac{0}}(0)+1}=\ens{\wt^s(0)}$.
\itemitem{}
Set $\free^s = \ens{\lf(s^+)}$. 

\item{(v)}
    $0<a,b\le\sz{s}-2$, so
    $s^+=t\tac{b}$, $s^-=t\tac{a}$, hence
\itemitem{}
    $\wt^s=\wt^{t\tac{b}}\cat\wt^{t\tac{a}}$,
\itemitem{}
    $\wt^{t\tac{b}}\diff\wt^{t\tac{a}}=\emptyset$.
\itemitem{}
Set $\free^s = \void$.

\item{(vi)}
    $0<a\le\sz{s}-2<b=\sz{s}-1$, so
    $s^+\not\in\SS$, hence
\itemitem{}
    $\wt^s=\wt^{s^-}$.
\itemitem{}
Set $\free^s = \void$.
\Qed

Contracting this somewhat:

\parSkip
\Number{(6.11)}
\Title{Corollary}
\Subtitle{Free position, second form}
For $s=t\tac{a,b}\in\SS$ such that $s^+\in\SS$
one of the following applies:
\item{(i)}
$a = b = 0$, $\wt^s=\xwt^{s^+}\cat\wt^{s^-}$ and
$\free^s=\ens{0,\lfs_n}$.
\item{(ii)}
$a > 0 = b$, $\wt^s=\xwt^{s^+}\cat\wt^{s^-}$ and
$\free^s=\ens{0}$.
\item{(iii)}
$a = 0 < b$, $\wt^s=\wt^{s^+}\cat\wt^{s^-}$ and
$\free^s=\ens{\lfs_n}$.
\item{(iv)}
$0 < a, b$, $\wt^s=\wt^{s^+}\cat\wt^{s^-}$ and
$\free^s=\void$.

For $s=t\tac{a,b}\in\SS$ such that $s^+\not\in\SS$
one of the following applies:
\item{(v)}
$\wt^s=\wt^{s^-}$ and
$\free^s=\void$.
\Qed

\secSkip
\Title{Summary}

\parSkip\noindent
\item{(i)}
Call $i\in\lf(s)$ {\it free} for $s$ if 
    $\wt^s(i)\not\in\rn(\wt^{s}\ric(\lf(s)\rid\ens{i}))$.
\item{(ii)}
Call $i\in\lf(s)$ {\it free} for $s$ {\it relative to $T$} 
(or simply free for $T$, when $s$ is clear from the context),
where $T$ is a subset of $\lfs=[0,\lfs)$ not containing $i$, if 
    $\wt^s(i)\not\in\rn(\wt^{s}\ric T)$.
(Thus, 'free for $s$' is 'free for $s$ relative to $\lfs\rid
\ens{i}$'.)

\parSkip
Thus, $i$ is free over $s_n$ just if $\wts_n(i)\not=\wts_n(j)$ 
for all $j\in\lfs_n$, $j\not= i$, hence just if the element 
$\M^{n}_i$ is not forced by the $\omega$-uniformity conditions 
to be isomorphic to any of the other elements in the $n$th array, 
	$\tup{\M^{(n)}_i:i\in\lfs_{n}}$.
Similarly, $i\in\Lfs_n^+$ is free (over $s_n$) relative to $T=
\Lfs_n^-$ if $\wts_n(i)\not=\wts_n(j)$ for all $j\in
\Lfs_n^-$, i.e., just if the element $\M^{n}_i$ is not forced 
by the $\omega$-uniformity condition to be isomorphic to any of 
the elements in, 
	$\tup{\M^{(n)}_i:i\in\Lfs_n^-}$.

\parSkip
What the determination lemma shows is that 
\item{(i)}
if $s_{n+1}(n+1)=n+1$, then $\wts_{n+1}(i)=\wts_n(i)$,
and the $\M^{(n+1)}_i$ are determined at the previous level.
\item{(ii)}
if $0<s_{n+1}(n),s_{n+1}(n+1)\le n-1$, then no $i\in\Lfs_{n+1}^+$
is free (for $s_{n+1}$ relative to $\Lfs_{n+1}^-$ (and vice versa
with the roles of '$+$' and '$-$' interchanged), so the 
$\M^{(n+1)}_i$ in the 'plus' half of $\Lfs_{n+1}$ are determined
by the $\M^{(n+1)}$ in the 'minus' half, and vice versa; but one of
these halfs is already determined by the previous level, so none
of them are 'new'.
\item{(iii)}
if $0=s_{n+1}(n)<s_{n+1}(n+1)\le n-1$, then no $i\in\Lfs_{n+1}^+$,
(in the 'plus' half) is free relative to $\Lfs_{n+1}^-\rid\ens{
\lfs_{n+1}^+}$ (i.e., the 'minus' half without its first element).
\item{(iv)}
if $0=s_{n+1}(n+1)<s_{n+1}(n)\le n-1$, then no $i\in\Lfs_{n+1}^-$,
(in the 'minus' half) is free relative to $\Lfs_{n+1}^+\rid\ens{0}$ 
(i.e., the 'plus' half without its first element).
\item{(v)}
if $0=s_{n+1}(n+1)=s_{n+1}(n)$, then no $i\in\Lfs_{n+1}^-\rid\ens{
\lfs_{n+1}^+}$ is free relative to $\Lfs_{n+1}^+\rid\ens{0}$ and
vice versa.


  

\chapSkip
\Chapter{[7] Short representations}

\secSkip\noindent
In this section and the next we give the stream-lined
representations promised earlier.

\secSkip
\Section{[7.1]}
\Title{Short representations of rank $\omega+\omega$}

\parSkip\noindent
A 'new' array element in the approximation of $A\in\UC(\omega+\omega)$
by the chain of components of some arbitrary (but fixed) 'basepoint'
$x$ in $A$, may only be introduced into the 'level' $n+1$ array of 
rank $\omega$ building blocks (i.e., the array used for the rank
$\omega+n$ approximation $\cmp{x}{\omega+n}{A}$,
    $\tup{\MM_n}_{n\in\omega}$),
if 
    $\frees_n = \free^{s_{n}}\not=\void$ 
(free position lemma), in which case the latter also specifies 
the possible positions of the 'new' element in the interval 
$[0,\lfs_n)$.  However, a 'new' element can appear in position 
$0\in\frees_{{n+1}}$ only if $\tp_n>0$, since for $\tp_n=0$ the 
level $n$ array 
    $\MM_n$
gets copied into the $n+1$ array
    $\MM_{n+1}^+$.
in which case, the 'new' element can therefore only appear in the 
'other' half of the array, at position $\lfs_{n+1}^+=\lfs_n$. 
(Obviously, the 'new' element may happen to be isomorphic with an 
'existing' one;  
the $\news_n$ indicate the positions of $\omega$ components not 
necessarily isomorphic to any of 
    $\M^{(k)}_i$, $k\le n$, $i<\lfs_k$,
hence potentially 'new', compared with the arrays for the components
of rank $\le\omega+n$ containing the 'basepoint') 

Thus, the 'new' position at level $n+1$ (rank $\omega+n+1$)
in the construction of a rank $\omega+\omega$ chain is therefore 
unqiuely determined whenever it exists, and this will determine 
the 'canonical ordering' of the short array for such $A$.

\parSkip
\Number{(7.1)}
\Title{Definition and Lemma}
Let $\ens{s_n}_n$ be derived from $\tup{\ss,\ps}\in\CSE(\omega+1)$
Set
\dispM{
    \news_0 := \ens{0}.
}
Suppose $s_{n+1}=t_{n+1}\tac{a_{n+1},b_{n+1}}$ such that 
$s_{n+1}^+\in\SS,$ i.e., $\tp_n < \omega$. Then one of the 
following applies, where we write $\news_n$ for $\new^{s_n}$:

{(i)}
    $a_{n+1} = b_{n+1} = 0$; 
then 
    $\wts_{n+1} = \xwts_{n+1}^+\cat\wts_{n+1}^-$.
Set
\dispM{
    \news_{n+1} := \cases{
    \ens{\lf_{n+1}^+} = \ens{\lf_n} & if $\tp_n = 0$ \cr
    \ens{0} & else (if $\tp_n = 1$). \cr}
}

{(ii)}
    $a_{n+1} > b_{n+1} = 0$; 
then, again,
    $\wts_{n+1} = \xwts_{n+1}^+\cat\wts_{n+1}^-$.
Set
\dispM{
    \news_{n+1} := \cases{
    \void & if $\tp_n = 0$ \cr
    \ens{0} & else (if $\tp_n = 1$). \cr}
}
(By the determination lemma, 
    $\rwt^{t\tac{a}}\sub\rwt^{t\tac{0}}\ric[1,\lf(t\tac{0}))$.)

{(iii)}
    $a_{n+1} = 0 < b_{n+1}$; 
then 
    $\wts_{n+1} = \wts_{n+1}^+\cat\wts_{n+1}^-$.
Set
\dispM{
    \news_{n+1} := \cases{
    \ens{\lfs_{n+1}^+} = \ens{\lfs_n} & if $\tp_n = 0$ \cr
    \void & else (if $\tp_n = 1$). \cr}
}

{(iv)}
    $a_{n+1}, b_{n+1} > 0$; 
then, again
    $\wts_{n+1} = \wts_{n+1}^+\cat\wts_{n+1}^-$.
Set
\dispM{
    \news_{n+1} := \void.
\Qed
}

\parSkip
\Title{Corollary}
\dispM{
    \eqalign{
    &\news_0 = \ens{0} 
    \hskip .2 in\hbox{(note,
    $\free^{\tup{0}} = \frees_0 = \ens{0}$)} \cr
    &\news_{n+1} = \cases{
    \void &
    if $\frees_{n+1} = \void$  \cr
    \ens{0} &
    if $\frees_{n+1}\supeq\ens{0}$ and $\tp_n = 1$ \cr
    \ens{\lfs_{n+1}^+} &
    if $\frees_{n+1}\supeq\ens{\lfs_{n+1}^+}$ and $\tp_n = 0$ \cr
    } \cr
    }
}

\parSkip
Note that $\news_0 \not= \void$ and that $\news_{n+1} = \void$ 
whenever $\tp_n=\omega$. Set
\dispM{
    \eqalign{
    & \rnnw 
    = \ens{n\in\omega : \news_n \not= \void} \cr
    & \rkkw = \sz{\rnnw} \cr
    }
}
Let $\ens{\nnw_k}_{k<\rkkw}$ enumerate $\rnnw$, i.e.,
\dispM{
    \eqalign{
    & \nnw_0 = 0 \cr
    & \nnw_{k+1} = \min\ens{n>\nnw_k:\news_n\not=\void}.\cr
    }
}
By convention, if $\rkkw<\omega$,
    $\nnw_k = \omega$,
for all $k\ge\rkkw$. Set 
    $\supp(\nnw)=\ens{k:\nnw_k<\omega}$, 
and 
    $\ran(\nnw)'=\ran(\nnw\ric\supp(\nnw))=\ran(\nnw)\cap\omega$.
Observe that $\rnnw$ $\nnw$ and $\rkkw$ are determined by 
    $\tup{\ss,\tp}$.

\parSkip
\Number{(7.2)}
\Title{Definition}
\Subtitle{Short arrays -- compressed version}
Define the {\it (compressed) short array} for 
\dispm{
    A
    = \A(\ss,\ps,\tup{\ss^{(n)}_i,\ps^{(n)}_i}_{
                      n\in\omega,0\le i<\lfs_n})
                           \in\UC(\omega+\omega)
}
by
\dispM{
    \U_k := \U^{(\nnw_k)}_{\news_{\nnw_k}}.
}
where the $\U^{(n)}_k$, $k\in\rwts_n$ are defined to 
satisfy
\dispm{
    \U^{(n)}_{\wts_n(i)} 
    = \M^{(n)}_i
    = \A(\ss^{(n)}_i,\ps^{(n)}_i)
} 
via the $\omega$-uniformity conditions. 
(Here, by abuse of notation,
    $\U^{(n)}_{\ens{k}} = \U^{(n)}_{k}$.)

\parSkip
Expanding this,
\dispM{
    \eqalign{
    \U_0 
     & = \U^{(0)}_0 \cr
    \U_{k+1} 
    & = \cases{
        \U^{(\nnw_{k+1})}_{\wts_{\nnw_{k+1}}(\lfs_{\nnw_{k+1}}^+)}
        & if $\tp_{\nnw_{k+1}-1} = 0$ \cr
        \U^{(\nnw_{k+1})}_{\wts_{\nnw_{k+1}}(0)}
      & if $\tp_{\nnw_{k+1}-1} = 1$ \cr
    } \cr
    }
}
which in turn gives
\dispM{
    \eqalign{
    \U_0 
     & = \M^{(0)}_0, 
    \cr
    \U_{k+1} 
    & = \cases{
        \M^{(\nnw_{k+1})}_{\lfs_{\nnw_{k+1}}^+}
         = \M^{(\nnw_{k+1})}_{\lfs_{\nnw_{k+1}-1}}
        & if $\tp_{\nnw_{k+1}-1} = 0$ \cr
        \M^{(\nnw_{k+1})}_{0}
        & if $\tp_{\nnw_{k+1}-1} = 1$ \cr
    } \cr
    }
}

\parSkip
Define a left-inverse of $\nnw$ on $\ran(\nnw)'$, 
by setting
\dispM{
    \kkw_n = \max\ens{k:\nnw_k\le n}.
}
Thus,
    $\kkw_n = k$ 
whenever $n$ is the $k$th integer such that
    $\news_n\not=\void$.
Note that $\rkkw = \ran(\kkw) = \supp(\nnw)$. 
%


\parSkip
Thus:
\item{(i)}
The $\news_{n+1}$ specify where at the level $n+1$ array a 
potentially 'new' element 
    $\U_{\kkw_n}=M\in\UC(\omega)$ ($\kkw_n\le n$)
can get 'added' to the array of building blocks --- either as 
    $\M^{(n+1)}_0$ 
if $\tp_n=1$, or as 
    $\M^{(n+1)}_{\lfs_{n+1}^+}=\M^{(n+1)}_{\lfs_n}$ 
if $\tp_n=0$. 

\item{(ii)}
Furthermore, level $n$ array
    $\MM_n$
gets 'copied' as
    $\MM_{n+1}^-$
if $\tp_n=1$, and as
    $\MM_{n+1}^+$
if $\tp_n=0$. 

\item{(iii)}
The remainder of the level $n+1$ array (i.e., the blocks not
in the 'copied from level $n$'-half --- hence in the 'half'
containing the 'new' position) is determined by the 
uniformity condition from the level $\omega+n$ array $\MM_n$. 

\parSkip
It remains to lay out the details on how the full double array 
    $\tup{\MM_n}_n$ 
is derived from the short array 
    $\tup{\U_k}_{k\in\rkkw}$ 
and vice versa.

\parSkip
\Title{Remark}
For all $k\in\ran(\kkw)=\dom(\nnw)$,
\item{(i)}
    $\tp_{\nnw_k} < \omega$.
\item{(ii)}
    $\kkw_{\nnw_k} = k$;
\item{(iii)}
    $\kkw_m\le k$ iff $m\le \nnw_k$;
\item{(iv)}
if $k+1<\ran(\kkw)$, then
    $\kkw_m=k$ iff $m\in [\nnw_k,\nnw_{k+1})$,
\item{}
if $k+1=\ran(\kkw)<\omega$, then
    $\kkw_m=k$ iff $m\in [\nnw_k,\omega)$.

\Proof
(i) holds, since by definition (e.g., (free position lemma,(vi)).
    $\ens{n:\news_n\not=\void}\subeq\ens{n:\tp_n<\omega}$.
(ii) follows from (iii), and (iii) follows from (iv), so we show (iv): 
But (iv) is true for $k=0$ ($\nnw_0=0=\kkw_0$). Supposing (iv) true 
for all $l\le k$, we get
    $\kkw_m=\kkw_{\nnw_k} = k$ 
    for $m\in [\nnw_k,\nnw_{k+1})$
if 
    $k+1\in\ran(\kkw)$,
while 
    $\kkw_m=\kkw_{\nnw_k}=k$ 
    for $m\in [\nnw_k,\omega)$,
if
    $k+1=\ran(\kkw)<\omega$.
\Qed 

\parSkip
\Number{(7.3)}
\Title{Lemma}
\item{(i)}
For all $n$,
    $\mwts_{n+1} > \mwts_n$
iff
    $\kwts_{n+1} = \kwts_n+1$
iff
    $\news_{n+1} \not= \void$.

\item{(ii)}
For all $n\in\ran(\nnw)'$,
    $\kkw_n = \kwts_n-1$ 
and so
    $\dom(\nnw) = \ran(\kkw) 
    = \sup_{n\in\omega}\kwts_n-1$;

\item{(iii)}
\dispM{
    \dom(\nnw)=\ran(\kkw) = \cases{
    \kwts_n
    & if $n$ is minimal such that for all $m\ge n$, 
      $\news_m=\void$ \cr
    \omega
    & if $\news_n\not=\void$ for infinitely many $n$, \cr
    }
}

\Proof
(i) Fix $n$ and set $s_{n+1} = s\tac{a,b}$. We may assume
$\tp_n<\omega$ as $\tp_n=\omega$ implies that 
    $\mwts_{n+1}= \mwts_n$ 
and 
    $\news_{n+1}=\void$. 
Thus, now, $s_{n+1}^+=s\tac{b}$ and $s_{n+1}^-=s\tac{a}$.
From the definition 
    $\mwts_{n+1} > \mwts_n$
iff
    $\kwts_{n+1} = \kwts_n+1$.

Set $\rwts_n^{\pm} = \rwt^{s_n^{\pm}}$
and $\mwts_n^{\pm} = \mwt^{s_n^{\pm}}$.

{\bf Case $a=b=0$.}
Then 
    $s_{n+1}=s\tac{0,0}$, 
and
    $\rwts_{n+1}^+=\rwt^{s\tac{0}}=\rwts_{n+1}^-$.
and
    $\mwts_{n+1}^+=\mwt^{s\tac{0}}=\mwts_{n+1}^-$.

{\bf Subcase $\tp_n=0$.}
Then
\item{}
    $s_n=s_{n+1}^+=s_{n+1}^-$,
\item{}
\dispm{
    \rwts_{n+1} 
    =\rxwts_{n+1}^+\cup\rwts_{n+1}^- 
    =\rxwt^{s\tac{0}}\cup\rwt^{s\tac{0}} 
    =\rxwts_n\cup\rwts_n 
    =\ens{\wts_n(0)+1}\cup\rwts_n 
    \supset\rwts_n, 
}
\item{}
    $\news_{n+1} = \ens{\lfs_{n+1}^+} = \ens{\lfs_n}$.

{\bf Subcase $\tp_n=1$.}
Then
\item{}
    $s_n=s_{n+1}^-=s_{n+1}^+$,
\item{}
\dispm{
    \rwts_{n+1}
    =\ens{\wts_n(0)+1}\cup\rwts_n
    \larger\rwts_n
}
(as before),
\item{}
    $\news_{n+1} = \ens{0}$.

{\bf Case $a>b=0$.}
Then 
    $s_{n+1}=s\tac{a,0}$, 
and
    $\rwts_{n+1}^+=\rwt^{s\tac{0}}\supset\rwts_{n+1}^-$.

{\bf Subcase $\tp_n=0$.}
Then
\item{}
    $s_n=s_{n+1}^+=s\tac{0}$,
\item{}
\dispm{
    \kwts_{n+1} 
    =\sz{\rxwts_{n+1}^+\cup\rwts_{n+1}^-} 
    =\sz{\rxwt^{s\tac{0}}\cup\rwt^{s\tac{a}}} 
    =\sz{\rwt^{s\tac{0}}\rid\ens{\wt^{s\tac{0}}(0)}
         \cup\ens{\xwt^{s\tac{0}}(0)}\cup\rwt^{s\tac{a}}} 
    =\sz{\rwt^{s\tac{0}}\cup\rwt^{s\tac{a}}} 
    =\kwt^{s\tac{0}} 
    =\kwts_n 
}
(since 
    $\ens{\wt^{s\tac{0}}(0),\wt^{s\tac{0}}(0)+1}
     \cap\rwt^{s\tac{a}}=\void$ 
and 
    $\rwt^{s\tac{0}}\supset\rwt^{s\tac{a}}$),
\item{}
    $\news_{n+1} = \void$.

{\bf Subcase $\tp_n=1$.}
Then
\item{}
    $s_n=s_{n+1}^-=s\tac{a}$,
\item{}
\dispm{
    \kwts_{n+1}
    =\sz{\rwt^{s\tac{0}}\cup\rwt^{s\tac{a}}}
    >\sz{\rwt^{s\tac{a}}}
    =\kwts_n
}
(following the preceding case),
\item{}
    $\news_{n+1} = \ens{0}$.

{\bf Case $a=0<b$.}
Then 
    $s_{n+1}=s\tac{0,b}$, 

{\bf Subcase $\tp_n=0$.}
Then
\item{}
    $s_n=s_{n+1}^+=s\tac{b}$,
\item{}
\dispm{
    \kwts_{n+1} 
    =\sz{\rwts_{n+1}^+\cup\rwts_{n+1}^-} 
    =\sz{\rwt^{s\tac{b}}\cup\rwt^{s\tac{0}}} 
    =\sz{\rwt^{s\tac{0}}} 
    >\sz{\rwt^{s\tac{b}}} 
    =\kwts_n, 
}
\item{}
    $\news_{n+1} = \ens{\lfs_{n+1}^+} = \ens{\lfs_n}$.

{\bf Subcase $\tp_n=1$.}
Then
\item{}
    $s_n=s_{n+1}^-=s\tac{0}$,
\item{}
\dispm{
    \kwts_{n+1}
    =\sz{\rwt^{s\tac{0}}}
    =\kwts_n
}
(again, following the preceding case),
\item{}
    $\news_{n+1} = \void$.

{\bf Case $a,b>0$.}
%
%
%
In that case 
    $\rwt^{s\tac{a}}=\rwt^{s\tac{b}}$ 
by the plus one range corollary, so 
    $\mwts_{n+1}
    =\mwt^{s\tac{a,b}}
    =\mwt^{s\tac{a}}
    =\mwt^{s\tac{b}}
    =\mwts_n$, 
but also $\news_{n+1}=\void$. 

(ii) 
We have 
    $s_0=\tup{0}$, 
    $\rwt^{s_\tup{0}}=\ens{0}$ 
and 
    $\kkw_0=0$.
By the preceding remark, part (iii),
    $\kkw_{n+1}=\kkw_n$ if $\news_{n+1}=\void$, 
    else $\kkw_{n+1} =\kkw_n+1$, 
hence 
    $\kkw_{n+1}=\kkw_n$ if 
        $\kwts_{n+1} = \kwts_n$, 
    else $\kkw_{n+1}=\kkw_n+1$ 
and 
    $\kwts_{n+1}=\kwts_n+1$.

(iii)
Follows from (ii).
\Qed

%
%

\parSkip
The preceding lemma permits the following alternative indexing 
of the short array, which will be exploited in the next 
section.

\parSkip
\Number{7.4}
\Title{Definiton}
\Subtitle{Short arrays -- uncompressed version}
Define the {\it uncompressed short array} 
    $\tup{\V_k}_{k\in\rmwts}$
for $A\in\UC(\omega+\omega)$ as above by
\dispM{
    \V_{\mwts_{\nnw_k}} = \U^{(\nnw_k)}_{\news_{\nnw_k}},
}
where 
    $\rmwts = \ens{\mwts_n:n\in\rnnw'}$.

\parSkip
\Title{Remark} 
The double arrays
    $\tup{\MM^{(n)}}_{n\in\omega}
    = \tup{\M^{(n)}_i}_{i\in\Lfs_n,n\in\omega}$
and
    $\tup{\UU^{(n)}}_{n\in\omega}
    = \tup{\U^{(n)}_k}_{k\in\rwt_n,n\in\omega}$,
and the short arrays
    $\UU= \tup{\U_k}_{k<\rkkw}$ 
and
    $\VV=\tup{\V_k}_{k\in\rmwts}
    \subeq\UC(\omega)$
deriving from a chain $A\in\UC(\omega+\omega)$ as given above
solve the following set of equations:
\dispM{
    \eqalign{
    & \M^{(0)}_0 = \U^{(0)}_0 = \U_0 = \V_0 \cr
    & \U_{\kkw_{n+1}}
    = \V_{\mwts_{n+1}}
    = \cases{
      \M^{(n+1)}_{\lfs_n}
      = \U^{(n+1)}_{\wts_{n+1}(\lfs_n)}
      & if $\news_{n+1}\not=\void$ and $\tp_n = 0$ \cr 
      \M^{(n+1)}_0
      = \U^{(n+1)}_{\wts_{n+1}(0)}
      & if $\news_{n+1}\not=\void$ and $\tp_n = 1$ \cr }
    } 
}
and
\dispM{
    \eqalign{
    \MM_{n+1} 
    & = \cases{
        \MM^{(n)}
        \wcat\tup{\U_{\kkw_{n+1}}}
        \wcat\MM^{(n+1)}_{(\lfs_n,\lfs_{n+1})} 
        & if $\news_{n+1}\not=\void$ and $\tp_n = 0$, \cr
        \tup{\U_{\kkw_{n+1}}}
        \wcat\MM^{(n+1)}_{(0,\lfs_{n+1}^+)} 
        \wcat\MM^{(n)} 
        & if $\news_{n+1}\not=\void$ and $\tp_n = 1$.  \cr
        } \cr
    }
}

\parSkip
Note that we may pass from the uncompressed to the compressed short 
array by applying the order-preserving mapping ('collapse') from 
$\rmwts$ to $\rkkw$, $\mwts_n\mapsto\kkw_n$ (effectively, a 
bounded search operation).

\parSkip
Note that generally $\V_{\mwts_{n+1}}$ and $\M^{(n+1)}_{\mwts_{n+1}}
=\U^{(n+1)}_{\mwts_{n+1}}$ refer to elements in different positions
of the level $n+1$ array (and may differ).


\parSkip
Call
    $\MM^{(+)}\in\Fun{\lfs^+}{\UC(\omega)}$
and
    $\MM^{(-)}\in\Fun{\lfs^-}{\UC(\omega)}$
{\it $s$-compatible} if
    $\MM = \MM^{(+)}\cat\MM^{(-)}$
satisfies the uniformity condition for $s$.

\parSkip
\Number{(7.5)}
\Title{Lemma}
    $\MM^{(+)}\in\Fun{\lfs^+}{\UC(\omega)}$ 
and 
    $\MM^{(-)}\in\Fun{\lfs^-}{\UC(\omega)}$ 
are $s$-compatible if the following apparently weaker condition 
holds:
\dispM{
    (\forall i<\lfs^+)
    (\forall k<\lfs^-)
    [\wt^s(i)=\wt^s(\lfs^++k)\;\then\;
     \M^{(+)}_i\isom\M^{(-)}_k],
    \eqno{(7.1)}
}
where 
    $\MM^{(\pm)}=\tup{\M^{({\pm})}_k}_{k<\lfs^{\pm}}$,
and where, as usual,
    $\lfs=\lf(s)$, $\lfs^{\pm}=\lf(s^{\pm})$.

\Proof
For the 'easy' direction note that for
    $i\in\Lfs^+$ and $j=\lfs^++k\in\Lfs^-$,
if
    $\wt^s(i) = \wt^s(j)$, 
then by the uniformity conditions, 
    $\M^{(+)}_i 
    =\M_i
    \isom\M_j
    =\M^{(-)}_k$,
as required.

So suppose 
    $\MM=\MM^{(+)}\cat\MM^{(-)}$
satisfies condition (7.1) for $s$ and let $i,j\in\Lfs$.
We have to show that either 
    $\wt^s(i)\not=\wt^s(j)$ 
or that
    $\M_i\isom\M_j$.

Fix $t\in\SS$ and $a,b\leq\sz{t}$ such that $s=t\tac{a,b}$.
Since $s^+$ is assumed in $\SS$, either 
    $b=0$ and $\wt^s=\xwt^{s^+}\cat\wt^{s^-}$
or 
    $b>0$ and $\wt^s=\wt^{s^+}\cat\wt^{s^-}$.

Thus, if $0<i<\lfs^+$ and $j=\lfs^++k\in\Lfs^-$, then by assumption
    $\wts(i)=\wts(j)$ 
implies that
    $\M_i=\M^{(+)}_i\isom\M^{(-)}_k=\M_j$.
Since we may assume $i<j$, we are left with the cases 
\item{(i)}
$i=0$,$j<\lfs$, 
\item{(ii)}
$0<i,j<\lfs^+$, and
\item{(iii)}
$\lfs^+\le i,j<\lfs$.

\parSkip
{\bf Case} $i = 0$: 

{\bf Subcase} $b=0$ and $\wt^s=\xwt^{s^+}\cat\wt^{s^-}$:
Then for $k\le\lfs^-$ and $j=\lfs^++k$,
\dispM{
    \wts(0) = \wts^+(0)+1
            = \wt^{t\tac{0}}(0)+1
            = \sz{t}+1
            > \wt^{t\tac{a}}(k)
            = \wt^{s^-}(k)
            = \wt^{s}(j),
}
by weight sequence lemma. Likewise, for $0<j'<\lfs^+$,
\dispM{
    \wts(0) = \wts^+(0)+1
             = \wt^{t\tac{0}}(0)+1
             = \sz{t}+1
             > \wt^{t\tac{0}}(j')
             = \wt^{s^+}(j')
             = \wt^{s}(j').
}
(Thus in this case the antecedent $\wts(0)=\wts(j)$ never applies,
except for $0=j$, when $\M_0\isom\M_j$ holds trivially.)

{\bf Subcase} $b>0$ and $\wt^s=\wt^{s^+}\cat\wt^{s^-}$:
Thus, for $k\le\lfs^-$, $j=\lfs^++k$
\dispM{
    \eqalign{
    &
    \wt^s(0) 
    = \wt^{s^+}(0)
    = \wt^{t\tac{b}}(0),
    \cr &
    \wt^s(j) 
    = \wt^{s^-}(k)
    = \wt^{t\tac{a}}(k),
    \cr }
}
So if $\wt^s(0) = \wt^s(j)$, then by (7.1)
\dispM{
    \M_0 
    = \M^{(+)}_0 
    \isom \M^{(-)}_k
    = M_j,
}
as required.

For $0<j'<\lfs^+$,
\dispM{
    \wt^s(j') 
    = \wt^{s^+}(j')
    = \wt^{t\tac{b}}(j').
}
By the plus one range lemma,
\dispM{
    \rwt^{s^+} 
    = \rwt^{t\tac{b}} 
    = \rwt^{t\tac{b}}
    = \rwt^{s^-},
}
so fix $k''\in\lfs^-$ such that
    $\wt^{s^+}(0) = \wt^{s^-}(k'')$
and let $j''=\lfs^++k''$, so
    $\wt^s(j')=\wt^s(0)=\wt^{s^+}(0)=\wt^{s^-}(k'')=\wt^s(j'')$. 
So if $\wt^s(0) = \wt^s(j')$, then in fact
    $\wt^s(0) = \wt^s(j'') = \wt^s(j')$,
so by (7.1) (applied twice to the second hypothesis),
\dispM{
    \M_0 
    = \M^{(+)}_0 
    \isom \M^{(-)}_{k''}
    \isom \M^{(+)}_{k'}
    = M_{j'},
}
as required.

{\bf Case} $0<i,j<\lfs^+$:
Then
\dispM{
    \eqalign{
    &
    \wt^s(i) 
    = \wt^{s^+}(i)
    = \wt^{t\tac{0}}(i),
    \cr &
    \wt^s(j) 
    = \wt^{s^+}(j)
    = \wt^{t\tac{0}}(j),
    \cr }
}
Again, by plus one range lemma,
\dispM{
    \wt^{t\tac{0}}(i) 
    = \wt^{t\tac{0}}(j)
    \in\rwt^{t\tac{a}},
}
so fix $k<\lfs^-$ such that
    $\wt^{t\tac{0}}(i) 
    = \wt^{t\tac{0}}(j)
    = \wt^{t\tac{a}}(k)$.
By (7.1),
\dispM{
    \M_i = \M^{(+)}_i 
    = \M^{(-)}_k 
    =\M^{(+)}_j 
    = \M_j.
}

{\bf Case} $\lfs^+\le i,j<\lfs$.
Then
\dispM{
    \eqalign{
    &
    \wt^s(i) 
    = \wt^{s^-}(i)
    = \wt^{t\tac{a}}(i),
    \cr &
    \wt^s(j) 
    = \wt^{s^-}(j)
    = \wt^{t\tac{a}}(j),
    \cr }
}
Again, by plus one range lemma,
\dispM{
    \wt^{t\tac{a}}(i) 
    = \wt^{t\tac{a}}(j)
    \in\rwt^{t\tac{b}},
}
so fix $k<\lfs^+$ such that
    $\wt^{t\tac{a}}(i) 
    = \wt^{t\tac{a}}(j)
    = \wt^{t\tac{b}}(k)$.
By two applications of(7.1),
\dispM{
    \M_i 
    = \M^{(-)}_i 
    = \M^{(+)}_k 
    =\M^{(-)}_j 
    = \M_j.
\Qed
}

\def\SubLfs{{\cal J}}

\parSkip
Make the obvious modifications to the above terminology when
dealing with 'shifted' or truncated sequences. E.g.:
Call
    $\MM^{(+)}\in\Fun{\Lf(s)^+}{\UC(\omega)}$
and
    $\MM^{(-)}\in\Fun{\Lf(s)^-}{\UC(\omega)}$
    {\it $s$-compatible} 
if $\MM^{(+)}$ and 
    $\tMM\in\Fun{\lf(s^-)}{\UC(\omega)}$
are compatible, where
    $\tM_i=\M^{(-)}_{i-\lfs^+}$.

\parSkip
Let 
    $\SubLfs^{\pm}\subeq\Lfs^{\pm}$.
Call
    $\MM^{(+)}\in\Fun{\SubLfs^+}{\UC(\omega)}$
and
    $\MM^{(-)}\in\Fun{\SubLfs^-}{\UC(\omega)}$
    {\it $s$-compatible} 
if
\dispM{
    (\forall i\in\SubLfs^+)
    (\forall j\in\SubLfs^-)
    [\wt^s(i) = \wt^s(j)\;\then\;
     \M^{(+)}_i\isom\M^{(-)}_j],
}
i.e., if
    $\MM = \tMM^{(+)}\cat\tMM^{(-)}$
satisfies the uniformity conditions for $s$ restricted to
$\SubLfs^+\cup\SubLfs^-$, where $\tMM^{(\pm)}$ are arbitrary
extensions of $\MM^{(\pm)}$ onto $\SubLfs^{\pm}$. Etc.

\parSkip
Let 
    $\SubLfs\subeq\Lfs^+$.
Call 
    $\tMM\in\Fun{\Lfs^+}{\UC(\omega)}$ 
{\it $s$-determinable on} $\SubLfs$ from 
    $\MM^{(-)}\in\Fun{\Lfs^-}{\UC(\omega)}$,
if 
    $\rn(\wts^+\ric\SubLfs)\subeq\rwts^-$,
and for all $i\in\SubLfs$,
    $\tMM_i = \MM^{(-)}_j$ 
for {\it every} $j$ such that
    $\wts^+(i) = \wts^-(j)$.
Likewise for 
    $\SubLfs\subeq\Lfs^-$.
Call 
    $\tMM\in\Fun{\Lfs^-}{\UC(\omega)}$ 
{\it $s$-determinable on} $\SubLfs$ from 
    $\MM^{(+)}\in\Fun{\Lfs^+}{\UC(\omega)}$,
if 
    $\rn(\wts^-\ric\SubLfs)\subeq\rwts^+$,
and for all $i\in\SubLfs$,
    $\tMM_i = \MM^{(-)}_j$ 
for {\it every} $j$ such that
    $\wts^-(i) = \wts^+(j)$.
If an $\tMM$ is given as $\tMM=\tMM_{\SubLfs}$ for 
$\SubLfs\subeq\Lfs^+$ or $\SubLfs\subeq\Lfs^-$, we 
shall simply say $\MM_{\SubLfs}$ is $s$-determinable
from $\MM^{(-)}$, respectively $\MM^{(+)}$.
If, e.g.,
    $\MM^{(-)}\in\Fun{\Lfs^-}{\UC(\omega)}$,
satisfies the uniformity conditions (on $\Lfs^-$), then
it suffices for
    $\tMM\in\Fun{\Lfs^+}{\UC(\omega)}$ 
to be $s$-determinable on $\SubLfs$ from $\MM^{(-)}$ that 
    $\rn(\wts^+\ric\SubLfs)\subeq\rwts^-$,
and for all $i\in\SubLfs$,
    $\tMM_i = \MM^{(-)}_j$ 
there is {\it some} $j$ such that
    $\wts^+(i) = \wts^-(j)$.

\parSkip
\Number{(7.6)}
\Title{Lemma}
\Subtitle{Array equivalence lemma}
Let
    $A 
    = \A(\ss,\ps)\arr{\MM_n}_{n\in\omega}
    \in\UC(\omega+\omega)$.
Then $\MM_0=\tup{\U_0}$ and for each $n$, one of (i), (ii) or (iii) 
applies:

{(i)}
{\bf Case} 
$\tp_n = \omega$: 

Then
    $\MM_{n+1}
    = \MM_n$.

{(ii)}
{\bf Case} 
$\tp_n = 1$: 

{\bf Subcase} 
$b_{n+1} = s_{n+1}(n+1) = 0$:

Then 
for
    $\tMM_{(0,\lfs_{n+1}^+)}$
$s_{n+1}$-determinable from
    $\MM_{n}$,
\dispM{
    \MM_{n+1}
    = \tup{\U_{\kkw_{n+1}}}
      \wcat\tMM_{(0,\lfs_{n+1}^-)}
      \wcat\MM_n.
}

{\bf Subcase} 
$b_{n+1} > 0$.

Then 
    $\MM_{n+1}$ 
is determined by $\MM_n$ in the sense that
\dispM{
    \MM_{n+1}^-
    = \MM_n
    \hbox{ and 
    $\MM_{n+1}^+$
    is $s_{n+1}$-determinable from
    $\MM_{n+1}^-$.
    }
}

{(iii)}
{\bf Case} 
$\tp_n = 0$: 

{\bf Subcase} 
$a_{n+1} = s_{n+1}(n) = 0$:

Then 
for
    $\tMM_{(0,\lfs_{n+1}^-)}$
$s_{n+1}$-determinable from
    $\MM_{n}$,
\dispM{
    \MM_{n+1}
    = \MM_n
      \wcat\tup{\U_{\kkw_{n+1}}}
      \wcat\tMM_{(0,\lfs_{n+1}^-)}.
}

{\bf Subcase} 
$a_{n+1} > 0$.

Then 
    $\MM_{n+1}$ 
is determined by $\MM_n$ in the sense that
\dispM{
    \MM_{n+1}^+
    = \MM_n
    \hbox{ and 
    $\MM_{n+1}^-$
    is $s_{n+1}$-determinable from
    $\MM_{n+1}^+$}
}

\Proof
The proof consists merely in pulling together earlier observations
and contains nothing new.

By definition
\dispm{
    \M^{(0)}_0 
    = \U_0 
    = \V_0.
}

If 
    $\tp_n = \omega$, 
then $s_{n+1}=s_n$, and per force
\dispm{
    \MM_{n+1}
    = \MM_n.
}

So suppose for the remainder that
    $\tp_n\in\omega$. 
Let $s=s_{n+1}=t\tac{a,b}$ (and
    $\lfs=\lfs_{n+1}=\lf(s_{n+1})$, 
    $\lfs^{\pm}=\lfs_{n+1}^{\pm}=\lf(s_{n+1}^{\pm})$, 
etc.).

{\bf Case} $a=b=0$. Thus $\news_n\subeq\ens{0,\lfs^+}$,
and that $\news_n\not=\void$. Thus, let 
    $k+1=\kkw_{n+1}$
and
    $l+1=\mwts_{n+1}$
and note that 
    $\kkw_{n+1}>k\ge\kkw_m$ for $m\le n$, 
    $\mwts_{n+1}>l\ge\mwts_m$ for $m\le n$
so that $\U_{k+1}=\V_{l+1}$ has not already been 'used' at any 
previous level in the 'approximation' to $A$ by the chain of
components $\A(s_n)\arr{\MM_n}$. Note also that now 
$s^+=t\tac{0}=s^-$, $\lfs^+=\lfs_n=\lfs^-$, etc.

{\bf Subcase} $\tp_n=0$: Then $\news_n=\ens{\lfs^+}$ 
and
\dispM{
    \MM_{n+1}
    = \MM_n
      \wcat\tup{\U_{k+1}}
      \wcat\tMM_{(0,\lfs_{n+1}^-)},
}
where $\tMM_{(0,\lfs_{n+1}^-)}$ is $s_{n+1}$ determinable
from $\MM_n$ (since
    $\rwt_{n+1}^+\rid\ens{\wt_{n+1}(0)}
    =\rwt^{t\tac{0}}\rid\ens{\wt^{t\tac{0}}(0)}
    =\rwt_{n+1}^-\rid\ens{\wt_{n+1}(\lfs_n)}$).

{\bf Subcase} $\tp_n=1$: Then $\news_n=\ens{0}$ 
and
\dispM{
    \MM_{n+1}
    = \tup{\U_{k+1}}
      \wcat\tMM_{(0,\lfs_{n+1}^+)}
      \wcat\MM_n,
}
where $\tMM_{(0,\lfs_{n+1}^+)}$ is $s_{n+1}$ determinable 
from $\MM_n$ (for reasons as before).

{\bf Case} $a>b=0$. Thus $\news_n\subeq\ens{0}$ and
that $\news_n=\void$ iff $\tp_n=0$. Also, 
    $\rwt^{s^-}
    =\rwt^{t\tac{a}}
    \subeq\rn(\wt^{t\tac{0}}\ric[1,\lf(t\tac{0}))
    =\rn(\wt^{s^+}\ric[1,\lfs^+))$.  

{\bf Subcase} $\tp_n=0$: Then $s^+=s_n$, $\lfs^+=\lf_n$
and $\news_n=\void$. By the preceding,
    $\rwt^{s^-}
    \subeq\rn(\wts_n\ric[1,\lfs_n))$,
so every entry
    $\tM^{n+1}_i, i<\lfs^-$
is determined (up to isomorphism) an entry
    $\M^{(n+1)}_i, i\in\Lfs_{n+1}^+\rid\ens{0}$
via $s_{n+1}$-compatibility. Thus,
\dispM{
    \MM_{n+1}
    = \MM_n
      \wcat\tMM_{\Lfs^-}
}
where the array 
    $\tMM_{\Lfs^-}$
is $s_{n+1}$-determinable from
    $\MM_n$
(since 
    $\rwt^{s^-}\subeq\rn(\wt^{s^+}\ric[1,\lfs^+))$, 
implies that no position in $\Lfs^-=[\lfs^+,\lfs)$ is 'free' 
relative to $\Lfs^+=[0,\lfs^+)$.)

{\bf Subcase} $\tp_n=1$: Then $s^-=s_n$, $\lfs^-=\lf_n$
and $\news_n=\ens{0}$. Let 
     $k+1=\kkw_{n+1}$ 
and
     $l+1=\mwts_{n+1}$ 
and note again that $\U_{k+1}=\V_{l+1}$ has not already been
'used' at any previous 'level' in the approximation of $A$. So,
\dispM{
    \MM_{n+1}
    = \tup{\U_{k+1}}
      \wcat\tMM_{(0,\lfs_{n+1}^+)}
      \wcat\MM_n
}
where the array 
    $\tMM_{\Lfs^+}$
is $s_{n+1}$-determinable from
    $\MM_n$
(since
    $\rxwt^{s^+}\rid\rwt^{s^-}=\ens{\wt^s(0)}$, 
i.e., no position in $(0,\lfs^+)$ being 'free' relative to 
    $\Lfs^-=[\lf^+,\lfs)$).

{\bf Case} $a=0<b$. Thus $\news_n\subeq\ens{\lfs^+}$ and
that $\news_n=\void$ iff $\tp_n=1$. Also, 
    $\rwt^{s^+}=\rwt^{t\tac{b}}
    \subeq\rn(\wt^{t\tac{0}}\ric[1,\lf(t\tac{0})))
    =\rn(\wt^{s^-}\ric[1,\lfs^-)))$.  

{\bf Subcase} $\tp_n=0$: Then $s^+=s_n$, $\lfs^+=\lf_n$
and $\news_n=\ens{\lf^+}$. Let 
    $k+1=\kkw_{n+1}$,
and
    $l+1=\mwts_{n+1}$,
and note again that $\U_{k+1}=\V_{l+1}$ is not already 'used' at 
any previous 'level'.  Also,
    $\rwts_n=\rn(\wts^-\ric[1,\lfs^-)$,
so all the entries 
    $\tM^{n+1}_i, i\in\Lfs^-\rid\ens{\lfs_{n+1}^+}$
are determined (up to isomorphism) by the entries in
    $\M^{(n+1)}_i, i<\lfs_{n+1}^+$
via $s_{n+1}$-compatibility. Thus,
\dispM{
    \MM_{n+1}
    = \MM_n
      \wcat\tup{\U_{k+1}}
      \wcat\tMM_{(\lfs_{n+1}^+,\lfs_{n+1}^-)}
}
where the array 
    $\tMM_{(\lfs_{n+1}^+,\lfs_{n+1})}$
is $s_{n+1}$-determinable from
    $\MM_n$.

{\bf Subcase} $\tp_n=1$: Then $s^-=s_n$, $\lfs^-=\lf_n$
and $\news_n=\void$. Now,
    $\rwt^{s^+}
    =\rwt^{t\tac{b}}
    \subeq\rwt^{t\tac{0}}
    =\rwt^{s^-}$,
so all the entries 
    $\tM^{n+1}_i, i<\lfs^+$
are determined (up to isomorphism) by the entries in
    $\M^{(n+1)}_i, i\in\Lfs_{n+1}^-$
via $s_{n+1}$-compatibility. So,
\dispM{
    \MM_{n+1}
    = \tMM_{\Lfs^+}
      \wcat\MM_n
}
where the array 
    $\tMM_{\Lfs^+}$
is $s_{n+1}$-determinable from
    $\MM_n$.

{\bf Case} $0<a,b$. Thus $\news_n=\void$ but also
    $\rwt^{s^+}=\rwt^{s^-}$, 
which implies that each 'half' of 
    $\MM_{n+1}$ 
is determined by the other via $s_{n+1}$ compatibility.

{\bf Subcase} $\tp_n=0$: Then $s^+=s_n$, $\lfs^+=\lf_n$,
and
\dispM{
    \MM_{n+1}
    = \MM_n
      \wcat\tMM_{\Lfs^-}
}
where the array 
    $\tMM_{\Lfs^-}$
is $s_{n+1}$-determinable from
    $\MM_n$.

{\bf Subcase} $\tp_n=1$: Then $s^-=s_n$, $\lfs^-=\lf_n$,
and
\dispM{
    \MM_{n+1}
    = \tMM_{\Lfs^+}
      \wcat\MM_n
}
where the array 
    $\tMM_{\Lfs^+}$
is $s_{n+1}$-determinable from
    $\MM_n$\Qed



\parSkip
The preceding lemma can be read in two ways:

Given
    $A 
    = \A(\ss,\ps)\arr{\MM_n}_{n\in\omega}
    \in\UC(\omega+\omega)$,
the preceding proof describes a procedure to derive the
short arrays
    $\UU=\tup{\U_k}_{k\in\rkkw}$
or
    $\VV = \tup{\V_l}_{l\in\rmwts}$
    $\subeq\UC(\omega)$
from the double array
    $\tup{\MM_n}_{n\in\omega}$
(and $\tup{\ss,\ps}\in\CSE(\omega+1)$).

Conversely, given
    $\UU = \tup{\U_k}_{k\in\rkkw}$
or
    $\VV = \tup{\V_l}_{k\in\rmwts}$
    $\subeq\UC(\omega)$,
the proof gives a procedure to derive the double-array
    $\tup{\MM_n}_{n\in\omega}$
from $\UU$, respectively $\VV$ (and $\tup{\ss,\ps}$).

\parSkip
Note the following 'procedural' definition of the short
arrays:

\item{(i)}
Set
    $\U_0=\V_0=\M^{0)}_0$.

\item{(ii)}
Suppose 
    $\tup{\U_i}_{i<k+1}$ 
is defined, let 
    $n = \nnw_{k+1}$,
and set
\dispM{
    \U_{k+1} = \cases{
    \M^{(n+1)}_0
    & if $\tp_n = 1$ \cr 
    \M^{(n+1)}_{\lfs_n}
    & if $\tp_n = 0$. \cr 
    }
}

\item{(ii)}
Suppose 
    $\tup{\V_l}_{l\in\rwts_n}$ 
is defined, let $n'$ minimal $\ge n$ such that
      $\news_{n'+1} \not= \void$
and set for $l'+1 = \mwts_{n+1} = $ the unique element in 
    $\rwts_{n'+1}\rid\rwts_n$
\dispM{
    \V_{l'+1} = \cases{
    \M^{(n'+1)}_0
    & if $\tp_{n'} = 1$ \cr 
    \M^{(n'+1)}_{\lfs_{n'}}
    & if $\tp_{n'} = 0$. \cr 
    }
}


\parSkip
Thus, a double array
    $\MM = \tup{\M^{(n)}_i}_{n\in\omega,i\in\lfs_n}$
and the 'short' arrays
    $\UU = \tup{\U_k}_{k\in\rkkw}$
or
    $\VV = \tup{\V_l}_{l\in\rmwts}$
related as in the array equivalence lemma can be used interchangably 
in the characterisation of $A\in\UC(\omega+\omega)$. Set
\dispM{
    \A(\ss,\ps)\rarr{\U_k}_{k\in\rkkw}
    =\A(\ss,\ps)\rarr{\V_l}_{l\in\rmwts}
    = \A(\ss,\ps)\arr{\MM_n}_{n\in\omega}
    \eqno{(7.2)}
}
whenever
    $\tup{\MM^{(n)}}_{n\in\omega}$
reduces to (and is derivable from)
    $\tup{\U_k}_{k\in\rkkw}$,
or
    $\tup{\V_l}_{l\in\rmwts}$,
and call the first two expression a {\it short representation}.

\parSkip
\Number{(7.7)}
\Title{Corollary} 
\Subtitle{Short rank $\omega+\omega$ representation lemma}

(a)
Every 
    $A\in\UC(\omega+\omega)$
has a representation of the short forms (7.2).

(b)
For every short array
    $\UU\in\Fun{\rkkw}{\UC\omega)}$,
respectively
    $\VV\in\Fun{\rmwts}{\UC\omega)}$,
the chain constructed as
    $\A(\ss,\ps)\rarr{\U_k}_{k\in\rkkw})$,
respectively
    $\A(\ss,\ps)\rarr{\V_l}_{l\in\rmwts})$
is in $\UC(\omega+\omega)$.
%
\Qed

\vfill\eject
\Section{[7.2]}
\Title{Generating local short arrays}

\parSkip\noindent
It remains to exhibit a procedure to extract (via 
    $\tup{\ss,\ps}\in\CSE(\omega+1)$)
the 'local' tagged short arrays
    $\UU_n
    =\tup{\tup{\U^{(n)}_k,k}:k\in\rwts_n}$,
from the 'global' short
    $\rkkw$-array $\UU$,
respectively
    $\rmwts$-array $\VV$,
for $A\in\UC(\omega+\omega)$.  (Here,
    $\U^{(n)}_{\wts_n(i)}=\M^{(n)}_i$
for $i\in\Lfs_n$.)
This will permit us to define
    $A = \A(\ss,\ps)\rarr{\UU}$,
respectively
    $A = \A(\ss,\ps)\rarr{\VV}$,
just if for some $x\in A$ and all $n$,
    $\cmp{x}{\omega+n}{A} 
    = \A(s_n)\rarr{\UU_n}$,
and vice-versa derive the short representation
from the short 'local' approximations ---
but without having to take the detour through the 
$\MM_n$ arrays.




We agreed to call a sequence
    $\UU=\tup{\U_k:k\in I}$,
(where $I$ is an initial segment of $\omega$)
a short array if for all $k,k'\in K$,
    $\UU_k\in\UC(\omega)$
and
    $\f\U_k\isom\f\U_{k'}$.
Call the short array $\UU$ {\it compatible} with
    $\tup{\ss,\ps}\in\CSE(\omega+1)$
if 
    $I = \rkkw$,
and call it {\it realisable} if for some
    $A\in\UC(\omega+\omega)$,
and some compatible
    $\tup{\ss,\ps}\in\CSE(\omega+1)$,
    $A = \A(\ss,\tp)\rarr{\UU}$.
(By the short $\omega+\omega$ representation theorem,
every short array is realisable with every compatible
    $\tup{\ss,\ps}\in\CSE(\omega+1)$.)
We also said we could dispense with tagging the local 
short arrays if we can provide a 'canonical' choice of
paths through $\Tree(s)$. We do not propose such a choice
'unifomly' for all $s\in\SS(<\omega)$, but there is a natural
such choice for each set
    $\ens{s_n}_n$
derived from $\tup{\ss,\ps}\in\UC(\omega)$, namely the
'embedding' path $\Path_n=\Path_n(\ss,\ps)$ determined by the 
$\tp_m,m<n$, defined as follows:

\parSkip
\Number{(7.8)}
\Title{Definition}
\item{(i)}
The top node is in $\Path_n$:
    $\tup{s_n;0}\in\Path_n$.
\item{(ii)}
Suppose $\sigma=\tup{t,k}\in\Path_n$ and let $m=\mod{t}-1$.
If $m=1$, the definition of $\Path_n$ is complete. Else set
\dispM{
    \displaylines{
    \tup{t^+,k+1}\in\Path_n\hbox{ if $\tp_{m-1}=0$,} \cr
    \tup{t^-,0}\in\Path_n\hbox{ if $\tp_{m-1}>0$.} \cr
    }
}

\parSkip
Immediately from the definition, for all $m<n$,
\dispM{
    \tup{s_m,k_m}\in\Path_n\hbox{ for some $k_m\le m$}.
}
This is true for $m=n$ with $k_0=0$. Supposing
$\tup{s_m,k_m}\in\Path_n$ and $m>0$,
\dispM{
    \displaylines{
    \tup{s_m^+,k+1}\in\Path_n
    \hbox{ and }
    s_{m-1} = s_m^+
    \hbox{ if $\tp_{m-1}=0$}, \cr
    \tup{s_m^-,0}\in\Path_n
    \hbox{ and }
    s_{m-1} = s_m^-
    \hbox{ if $\tp_{m-1}>0$}, \cr
    }
}
so $\tup{s_{m-1},k_{m-1}}\in\Path_n$ for $k_{m-1}=k_m+1$
if $\tp_{m-1}=0$, else $k_{m-1}=0$. Hence we are justified
in calling the paths $\Path_n$ 'determined by $\tup{\ss,\tp}$.

\parSkip
Observe that the canonical ordering imposed on the local arrays 
by taking this path is an initial segment of the ordering
imposed by 'novelty' on the short array for $A$. We can
therefore define local array (untagged) short array as the
initial segment of the (global) short array $\UU$ 'used'
at the given level:

\parSkip
\Number{(7.9)}
\Title{Definition}
\dispM{
    \UU_n'
    =\tup{\U_k:k\le\kkw_n}
    =\tup{\U^{(\nnw_k)}_{\news_{\nnw_k}}:k\le\kkw_n}
}
I.e., 
    $\UU_n'=\UU\ric\kkw_n$. 

\parSkip
Since $\kkw_n=\kwts_n-1$ (lemma 7.3), we can pass from the untagged
to the tagged array by ordering 
    $\rwts_n=\ens{v_0<\ldots,v_{\kkw_n}}$
and setting
    $\tup{\U_{v_k},v_k} = \tup{\U^{(\nnw_k)}_{\news_{\nnw_k}},v_k}$.
I.e.,
    $\U^{(n)}_{\wts_n(i)} = \U_k$
just if $\wts_n(i)$ is the $k$th largest element in $\rwts_n$.
Being able to dispense the tags, we shall now use '$\UU_n$' to
refer to the untagged local short arrays.

\parSkip
\Number{(7.10)}
\Title{Definition}
Define $\UU_n^{\pm}$ (relative to $\ss,\ps$ and \UU), 
so that $\UU_n^{\pm}$ is the 'untagged' version of the
original tagged $(\UU_n)_{s_n}^{\pm}$ of section 3. 
(We may assume $n>0$.)
\item{}
If $\news_n=\void$, set
$\UU_n^-=\UU_n^+=\UU_n$.
\item{}
If $\news_n=\ens{0}$ and $\tp_{n-1}=\lfs_n^+$ set
$\UU_n^-=\UU_n\ric\kkw_{n-1}$ and 
$\UU_n^+=\UU_n$. 
\item{}
If $\news_n=\ens{\lfs_n^+}$ and $\tp_{n-1}=0$ set
$\UU_n^+=\UU_n\ric\kkw_{n-1}$ and 
$\UU_n^-=\UU_n$. 

\parSkip
In set notation:
\item{}
If $\news_n=\ens{0}$ and $\tp_{n-1}=\lfs_n^+$,
$\wran\UU_n^-=\wran\UU_n\rids{\U_{\news_n}}$ and 
$\wran\UU_n^+=\wran\UU_n$. 
\item{}
If $\news_n=\ens{\lfs_n^+}$ and $\tp_{n-1}=0$,
$\wran\UU_n^+=\wran\UU_n\rids{\U_{\news_n}}$ and 
$\wran\UU_n^-=\wran\UU_n$. 

\parSkip
\Number{(7.11)}
\Title{Corollary}

(1)
The sequence of string $\UU_n$ over $\UC(\omega)$ is {\it generated
by (derived from)} $\tup{\ss,\ps}$ and $\UU$ in the follwing sense:
\dispM{
    \UU_0=\ens{\U_0}
}
and for all $n$,
\dispM{
    \UU_{n+1} = \cases{
    \UU_n
    & if $\news_{n+1}=\void$ \cr
    \UU_n\tac{\U_{\news_{n+1}}}
    & if $\news_{n+1}\not=\void$. \cr
    }
}

(2)
%
Let 
    $\ens{\UU_n}_{n\in\omega}$
where
    $\UU_n = \tup{\U_k}_{k\le\kkw_n}$,
be the sequence generated by
    $\tup{\ss,\ps}\in\CSE(\omega+1)$ 
and
    $\UU=\tup{U_k}_{k\in\rkkw}\in\Fun{\rkkw}{\UC(\omega)}$,
Then for each $n$,
    $\tup{\U_k:k\in\rwts_{n+1}}$
is an end-extension
    $\tup{\U_k:k\in\rwts_n}$,
and the 'full' arrays 
    $\MM_n
    =\MM_{\Lfs_n}
    =\tup{\M^{(n)}_i}_{i\in\Lfs_n}$ 
are determined by
    $\M^{(n)}_i 
    = \U_k$
just if $\wts_n(i)$ is the $k$th largest element in $\rwts_n$.

(3)
If 
    $A=\A(\ss,\ps)\rarr{\UU}$ 
and
    $\UU_n$ are derived from 
    $\tup{\ss,\ps}\in\CSE(\omega+1)$ 
and
    $\UU\in\Fun{\rkkw}{\UC(\omega)}$,
then for some $x\in A$ and all $n$,
    $\cmp{x}{\omega+n}{A}
    = \A(s_n)\rarr{\UU_n} 
    = \A(s_n)\arr{\MM_n}$,
and the approximating chain of sub-components can be generated
by the following variant of Pastjin's uniform chain constructor:
\dispM{
    \eqalign{
    \A(\gtup{0})\rarr{\U_0} &= \U_0,\cr
    \A(s_n)\rarr{\UU_n}
    &= \cases{
    \A(s_n^-)\rarr{\UU_n^-}.\Z = \A(s_n^-)\rarr{\UU_n}.\Z 
    &if $s_n(n)=n$ \cr
    \A(s_n^+)\rarr{\UU_n^+} + \A(s^-)\rarr{\UU_n^-}.\N
    &if $s_n(n)<n$. \cr
    } \cr
    }
}
Note that the $\UU_n$ 'argument' does not get truncated in the 
'$\new^s=\void$' cases arising both under $s(n)=n$ and $s(n)<n$.
(Of course in the '$s(n)=n$' case $\news_n$ is always $\void$.)
\Qed

\chapSkip
\Chapter{[8] Short representations of tails}

\secSkip\noindent
Passing from the short array for a uniform chain of
rank $\omega+\omega$ to that of its tail involves an
additional twist, hinted at in the tail revelation lemma.

\secSkip
\Section{[8.1]}
\Title{Short representations of rank $\omega+\omega$}

\parSkip\noindent
The preceding shows how to replace the arrays
    $\MM_n$
in the description of $A\in\UC(\omega+\omega)$.
Now
\dispM{
    \f\A(\ss,\pi)\arr{\MM_n}_{n\in\omega} 
    = \A(\f\ss,\f\pi)\arr{\f_{s_n}\MM_n}_{n\in\omega}
}
and we would like similarly to replace the arrays
    $\f_{s_n}\MM_n$
in the description of tails $\f A\in\UT(\omega+\omega)$
(subsequently referred to as 'tail'-arrays). But note that 
the tail-arrays 'telescope' (since for
tails, $\tp_n\const 0$); i.e.,
    $\f_{s_{n+1}}\MM_{n+1} 
    = \f_{s_{n+1}^{\pm}}\MM_{n+1}^{\pm}\wcat\MM_{n+1}^-$
where the 'plus half',
    $\f_{s_{n+1}^{\pm}}\MM_{n+1}^{\pm}$
belongs to the preceding level. Thus the building blocks 
'new' in the tail array at level $n+1$ are determined by
the 'minus half', $\MM_{n+1}^-$. This is also clear by
evaluating the expression for tail-chains (section 6):
\dispM{
    \eqalign{
    \f A(\ss,\ps)\arr{\MM_n}_{n\in\omega}
    &\isom \f\M^{(0)}_0 
    + \sum_{n\in\omega}\A(\ss\ric{n+1})
      \arr{\MM_{n+1}^-}.\N \cr
    }
}

Now compare how the arrays
    $\tup{\MM_{n+1}}$
and 
    $\tup{\MM_{n+1}^-}$
relate to
    $\tup{\MM_n}$
and $\UU$ (via $\tup{\ss,\ps}$). Here, the 'range' of an
array is the set of possibly non-isomorphic building-blocks.
Namely,

\parSkip

\dispM{
    \eqalign{
    \MM_0 
    & = \tup{\M^{(0)}_0} = \tup{\U_0} \cr
    \wran\MM_{n+1}
    & = \cases{
      \wran\MM_n
      & if $\news_{n+1} = \void$ \cr
      \ens{\U_{\wts_{n+1}(0)}}
      \cup\wran\MM_{n}
      & if $\news_{n+1} \not= \void$
        and $\tp_n = 1$ \cr
      \wran\MM_n
      \cup\ens{\U_{\wts_{n+1}(\lfs_n)}}
      & if $\news_{n+1} \not= \void$
        and $\tp_n = 0$ \cr
    }\cr
    }
} 
while
\dispM{
    \eqalign{
    \f\MM_0 
    & = \tup{\f\M^{(0)}_0} = \tup{\f\U_0} \cr
    \wran\MM_{n+1}^-
    & = \cases{
      \wran\MM_n
      & if $\tp_n > 0$ or $\news_{n+1} = \void$ \cr
      \ens{\U_{\wts_{n+1}(\lfs_n)}}
      \cup\wran\MM_n
      & if $\tp_n = 0$ and $\news_{n+1} \not= \void$ \cr
    } \cr
    }
} 

%

Thus, in the 'tail case' one of the following holds:

{(i)}
    $\news_{n+1} = \void$:
In this case
    $\MM_{n+1}^-$
is determined by $s_{n+1}$ compatibility with
    $\MM_n$.
Or,

{(ii)}
    $\new_{n+1} \not= \void$,
but
    $\tp_n = 1$:
In this case
    $\MM_{n+1}^-
    = \MM_n$;
thus, in this case, the 'new' $\omega$ component is
'hidden', or 'invisible' to the tail. Or,

{(iii)}
    $\new_{n+1} \not= \void$
and
    $\tp_n = 0$:
Then 
    $\MM_{n+1}^-
    = \tup{\U_{\news_{n+1}}}\cat
      \MM^{(n)}_{(\lfs_{n+1}^+,\lfs_{n+1})}$
where
    $\U_{\news_{n+1}} = \U_{\kkw_{n+1}}$
and where the array 
    $\MM^{(n)}_{(\lfs_{n+1}^+,\lfs_{n+1})}$
$s_{n+1}$-determinable from
    $\MM_n$.
In this case, therefore, the 'new' $\omega$ component is
'visible' to the tail.

%
%

\parSkip
Hence, a 'new' rank $\omega$ building block enters the $n+1$st
'tail' array only if $\tp_n = 0$. If $\tp_n > 0$, the 'new' 
element may have to be traced back several levels back, namely 
to the largest level $k<n$ such that $\tp_k = 0$. 
E.g., if $n>0$ and $\tp_{n-1} = 0$,
    $\M^{(n)}_{\lfs_{n-1}} 
    = \U_{\mwts_n}$.
If $n>0$ and $\tp_{n-1} > 0$ also, then the last 'new' $\omega$ 
component was introduced at an even earlier level. 
More precisely, a 'new' element $\U_k$ enters the 'tail' array,
essentially
    $\MM_{n+1}^-$,
just if either
\item{(i)}
    $\tp_n=0$ and $n=\nnw_k$ 
or
\item{(ii)}
    $n$ is minimal such that for some $m<n$,
    $\tp_m=1$ and $m=\nnw_k$, 
    $\tp_{n'}=0$ for ${n'}\in (m,n)$, 
    and $\tp_n=1$ 
\par\noindent
(and while $\tp_{n'}=0$ in the second scenario, some 'later new' 
elements may have entered the tail-array, namely $\U_{k'}$ 
where $n'=\nnw_{k'}$ for $n'\in(m,n))$.

\parSkip
\Number{(8.1)}
\Title{Definition}
\Subtitle{$\WW_{\rllw} = \f\UU_{\rkkw}$}

We define the map that completes the bottom arrow in the following 
diagramme:
\dispM{
\matrix{
&
& \f
&
\cr
& \A(\ss,\pi)\arr{\MM_{n}}_n 
& \longrightarrow 
& \A(\f\ss,\f\pi)\arr{\f_{s_n}\MM_{n}}_n
\cr
& \downarrow
&
& \downarrow
\cr
& \tup{\MM_{n}} _{n\in\omega}
& \longrightarrow 
& \tup{\f_{s_n}\MM_{n}}_{n\in\omega}
\cr
& \downarrow
&
& \downarrow
\cr
& \tup{U_k}_{k\in\rkkw} 
& \longrightarrow 
& \tup{\W_l}_{l\in\rllw}
\cr
}}
I.e., set
    $\WW_{\rllw} = \f\UU_{\rkkw}$
whenever
\item{}
    $\UU_{\rkkw} = \tup{U_k}_{k\in\rkkw}$ 
is the short array for some
    $\A(\ss,\pi)\arr{\MM_{n}}_n$, 
and
\item{}
    $\WW_{\rllw} = \tup{\W_l}_{l\in\rllw}$
is the short array for 
    $\f A=\A(\f\ss,\f\pi)\arr{\f_{s_n}\MM_{n}}_n$, 

\noindent
where
    $\f\ss=\ss\ric\omega\tac{0}$, 
where
    $\f\pi=\tup{0^{\omega}}$,
and where
    $\f_{s_n}(\MM_{n})$ 
is as defined in section 3.

\parSkip
{\it Outline:}
Set $\W_0=\f\U_0$.

While $\tp_n=0$ ($n=0,1,2,\dots$), set $\W_l=\U_k$ whenever $n$ 
is such that $\nu(n+1)\not=\void$ (by implication,
    $\nu(n+1)=\ens{\ell_n}=\ens{\ell_{n+1}^+}$), 
and $l$ is minimal such that $\W_l$ is not yet defined.

For the first $n_1\ge 0=n_0+1$ (i.e., we set $n_0=-1$) such that 
$\tp_{n_1}>0$, set $\W_l=\U_0$ ($l$ as before).

While $\tp_n=0$ ($n=n_1+1,n_1+2,\ldots$), set $\W_l=\U_k$ whenever 
$n$ is such that $\nu_n\not=\void$ (and $l$ as before).

For the second $n_2>n_1$ such that $\tp_{n_1}>0$, set
$\W_l=\U_{k'}$ ($l$ as before), provided it was the case
that $n_1+1=\nnw_{k'}$ (i.e., $\nu(n_1+1)\not=\void$ - else
$\W_l$ remains undefined).
Etc.

\parSkip
{\it Detail:}
Let $\ens{n_i:i\in I}$ be the set of $n$ such that $n=0$ 
or $n>0$ and $\tp_{n-1}>0$. 
(Here $I$ is an initial segment of $\omega$, i.e.,
$I=\omega$ or $I\in\omega$.) Set $n_{I}=\omega$ if
$I$ is finite (i.e., $I\in\omega$). Thus, $\omega=
\bigcup_{i\in I}[n_i,n_{i+1})$.

\def\kMax{\overline{k}}
\def\lMax{\overline{l}}

Set 
    $\U_k^{(n)}:=\loops$ (to mean: 'is undefined'), 
if $n<\nnw_k$, and set 
    $\U_k^{(n)}:=\U_k$, 
if $n\ge\nnw_k$ ($k\in\rkkw$). Say 
    $\U_k^{(n)}:=\halts$ (to mean 'is defined'), 
if 
    $\U_k^{(n)}=\U_k\not=\loops$. 
Let 
    $\kMax=\max\ens{k:\nnw_k<n}
    =\max\ens{k:\U_k^{(n-1)}\halts}$. 
(Note $\kMax(0)=\max\void=0$.)

Define 
    $\W_0^{(0)}:=\W_0:=\f\U_0$ 
and 
    $\W_l^{(0)}:=\loops$ for
all $l>0$. Let 
    $\lMax=\max\ens{l:\W_l^{(n-1)}\halts}$. 
(So, $\lMax(0)=0$.)

Let now $n>0$.

Suppose $n\in(n_i,n_{i+1})$ for some $i$. 
If $n\not=\nnw_k$, for any $k$, set 
    $\lMax(n):=\lMax(n-1)$ 
and 
    $\W_l^{(n)}:=\W_l^{(n-1)}$
for all $l$. Else, for the unique $k$ such that $n=\nnw_k$, 
set
    $\lMax(n):=\lMax(n-1)+1$, 
    $\W_{\lMax(n)}^{(n)}:=\U_k$ 
and 
    $\W_l^{(n)}:=\W_l^{(n-1)}$ 
for all $l\not=\lMax(n)$.

Suppose $n=n_i$ for some $i$ (since $n>0$, $i>0$). If
$n_{i-1}\not=\nnw_k$, for any $k$, set 
    $\lMax(n):=\lMax(n-1)$ 
and
    $\W_l^{(n)}:=\W_l^{(n-1)}$.  
Else, for the unique $k$ such that $n_{i-1}=\nnw_k$, set 
    $\lMax(n):=\lMax(n-1)+1$, 
    $\W_{\lMax(n)}^{(n)}:=\U_k$ 
and 
    $\W_l^{(n)}:=\W_l^{(n-1)}$ 
for all $l\not=\lMax(n)$.

{\it End of definition.}

\parSkip
\Number{(8.2)}
\Title{Corollary}
For all $n$, $\lMax(n)\in\ens{\kMax(n),\kMax(n)+1}$ and
for all $i\in I$, if $n_i\not=\nnw_k$, for any $k$,
then for all $n\in[n_i,n_{i+1})$, $\lMax(n)=\kMax(n)+1$;
if $n_i=\nnw_k$ for some unique $k$, then for all
$n\in[n_i,n_{i+1})$, $\lMax(n)=\kMax(n)$.

\Proof
By induction.

For $n=0$, $\lMax(0)=0=\kMax(0)$.

Let $n>0$, and suppose the first claim is true for all $n'<n$, 
and if $n=n_{i+1}$, $i\in I$, the second claim is true for $i$.

If $n\in(n_i,n_{i+1})$, then from the definiton,
$\lMax(n)-\lMax(n-1)=\kMax(n)-\kMax(n-1)$, so the claim
continues to hold, while the second claim holds
by continuing the induction over $n''\in (n_i,n_{i+1})$.

Suppose then $n=n_i$ for some $i$, and we may assume
$i>0$. If either both $n_{i}$ and $n_{i-1}$ are in
$\ens{\nnw_k}_k$ or neither is, then again 
$\lMax(n)-\lMax(n-1)=\kMax(n)-\kMax(n-1)$, for all $n\in
(n_i,n_{i+1})$ by hypothesis, and for $n=n_i$ by
definition, so the induction proceeds as in the previous 
paragraph. Else one of the following holds:

Either $n_{i-1}\not=\nnw_k$, for any $k$, while $n_i=
\nnw_{k'}$ for some (unique) $k'$; then $\lMax(n_i)=
\lMax(n_i-1)=\kMax(n_i-1)+1=\kMax(n_i)$.

Or $n_{i-1}=\nnw_k$, for some (unique) $k$, while $n_i\not
=\nnw_{k'}$ for any $k'$; then $\lMax(n_i)=\lMax(n_i-1)+1=
\kMax(n_i-1)+1=\kMax(n_i)+1$.
\Qed

\parSkip
\Number{(8.3)}
\Title{Corollary}
\dispM{
    \f\A(\ss,\ps)\rarr{\U_i}_{i\in\rkkw}
    = \A(\ss\ric\omega\tac{0})\rarr{\W_i}_{i\in\rllw}
    = \f\M^{(0)}_0 
    + \sum_{n\in\omega} \A(\ss\ric{n+1})
      \arr{\MM_{n+1}^-}.\N.
}
%

\parSkip
The next lemma tells where (and if) a given building block
$\U_k$ appears in the tail.

\parSkip
\Number{(8.4)}
\Title{Lemma}
\Subtitle{Tail relocation lemma}
For $0<k<\rkkw$, let 
    $n_k=\nnw_k-1$ 
(being thus also the least $n$ such that 
    $\kkw_{n+1}=k$). 

\item{(i)}
Either for all $0<k<\rkkw$,

\itemitem{(a)}
    $\tp_{n_k}=0$,
in which case 
    $\U_k=\M^{(n_k+1)}_{\lfs_{n_k}}$
gets inserted into $\WW_{\rllw}$ at 'stage $n_k$' (i.e., as 
    $\W_k$ or $\W_{k+1}$); 

\item{}
or 

\itemitem{(b)}
    $\tp_{n_k} = 1$,
but there is some $m\ge n_k+1$ such that again
    $\tp_m = 1$,
in which case
    $\U_k=\M^{n_k+1})_0=\M^{(m)}_0=\M^{(m+1)}_{\lfs_m}$
gets inserted into $\WW_{\rllw}$ at stage $m$, for the 
least such $m$. 

\item{(ii)}
Or there is a unique $0<k<\rkkw$ such that
    $\tp_{n_k} = 1$,
in which case 
    $\U_k=\M^{(n_k+1)}_0$, 
but for all $m\ge n_k+1$,
    $\tp_m\in\ens{0,\omega}$,
and the copy-definition
    $\U_k=\M^{(m)}_0=\M^{(m+1)}_{\lfs_m}$
is never made, whence $\U_k$ does not get inserted into
$\WW_{\rllw}$ (although this may happen to some other
elements $\U_{k'}$ isomorphic to $\U_k$.)

\Proof
Suppose (i) fails and let $k_0$ be the least $k$ such 
that (i) fails for $k$, i.e.,
    $\tp_{n_{k_0}} = 1$,
but there is no $m\ge n_{k_0}$ such that 
    $\tp_m = 1$, 
i.e., for all $m\ge n_{k_0}$, 
    $\tp_m\in\ens{0,\omega}$. 
For any 
    $k\in (k_0,\rkkw)$
we have 
    $\tp_{n_k}\not=\omega$, 
hence by the assumption that (i) fails at $k_0$, 
    $\tp_{n_k}=0$; 
so these $k$ all satisfy (i) via the first disjunct. Thus,
$k_0$ uniquely satisfies (ii).
\Qed

\parSkip
\Title{Remark}
\Subtitle{Scenarios for $\tup{W_l:l\in\rllw}$}

\parSkip(a)
Suppose first that for all $m$, 
    $\tp_m = 0$.
In this case $\U_0$ always gets 'copied' as 
    $\M^{(m+1)}_0$
into the 'plus' half of the level $m$ building blocks, 
hence never appears (as 
    $\M^{(m+1)}_{\lfs_m}$) 
in a 'tail'-array
    $\MM_{m+1}^-$,
while every $\U_k$ for $k > 0$ gets 'inserted' 
(as $\M^{(n_k+1)}_{\lfs_{n_k}}$) in the 'tail'-array
    $\MM_{n_k+1}^-$.
at the same level at which it appears in the full array
    $\MM_{n_k+1}$.
In particular, 
    $\W_0 = \f\U_0$,
and
    $\W_{k+1} = \U_{k+1}$
for $k\in\rkkw$, 
and also
    $\rllw=\rkkw$.
$\U_0$ never appears in the tail array. 

Note that if $\tp_m=\omega$, then 
    $\f_{s_{m+1}}\MM_{m+1} 
    = \f_{s_{m+1}^-}\MM_{m+1}\wcat\MM_{m+1}$,
so $\U_0$ does get copied into the tail-array at level $m+1$.
If $\tp_m=1$ for some least $m$, then $U_0=\M^{(m)}_0$ gets
copied as $\M^{(m+1)}_{\lfs_{m+1}^+}$ into the 'minus half'
of the level $m+1$ building blocks, hence also into the 
tail-array
    $\f_{s_{m+1}}\MM_{m+1} 
    = \f_{s_{m+1}^+}\MM_{m+1}^+\wcat\MM_{m+1}^-$
of level $m+1$.

\parSkip(b)
So suppose that 
    $\tp_m = 1$ 
for some $m\in\omega$. Let $m_0$ be the least $m$ such that 
    $\tp_m = 1$; 
then for some $k=k_1\in\rkkw$, $m_0$ also is the 
least $n$ such that $\kkw_n=k_1$. Then for all 
$n\le m_0$ such that 
    $\tp_n < \omega$ 
we have in fact 
    $\ps_n = 0$, 
and so
    $\U_k 
    = \M^{(n_k+1)}_{\lfs_{n_k+1}^+}$
for $k\in (0,k_1)$, while level $m_0+1$ represents the
earliest opportunity to copy $\U_0$ into the 'minus half',
and hence the tail-array,
    $\U_0 
    = \M^{(m_0+1)}_{\lfs_{m_0+1}^+}$.
(The extreme case $m_0=k_1 = 0$ is possible, in which case the
prescription for $k\in (0,k_1)$ will be vacuous.)

\parSkip(b.1)
Suppose next and that the first case of the tail relocation
lemma applies, i.e., that for all $k\in\rllw$, (i) holds. 
Let $k_{2i},i=1,2\ldots$ such that 
    $\tp_{n_{k_{2i}}} = 1$ 
and $m_i,i=1,2,\ldots$ such that $m_i$ is the least 
$m\ge n_{k_{2i}}+1$ such that 
    $\tp_{m_i} = 1$. 
Let $k_{2i+1},i=1,2\ldots$ such that $\kkw_{m_i}=k_{2i+1}$. Then
    $\U_k 
    = \M^{(n_k+1)}_{\lfs_{n_k+1)}^-}$
for $k\in (k_{2i},k_{2i+1})$, while
    $\U_{k_{2i}} 
    = \M^{(m_i+1)}_{\lfs_{m_i+1}^-}$.
Thus,
    $\W_k = \U_{k+1}$
for
    $k\in\bigcup_i [k_{2i},k_{2i+1})$,
where $k_0:=0$, and
    $\W_{k_{2i+1}} = \U_{k_{2i}}$;
in particular
    $\W_{k_1} = \U_{k_0} = \U_0$.

If 
    $\rkkw = n+1\in\omega$, 
then $n$ is maximal such that 
    $\news_{\nnw_n} \not= \void$.
If $n=m_{i^0}$ for some last $i^0$, the definition is again 
complete, with
    $\ens{\W_k:k\in\rkkw} = \ens{\U_k:k\in\rkkw}$.
Else $m_{i^0}<n$ and 
    $\tp_{n_k}=0$ 
for all 
    $k\in (k_{2i^0+1},\rkkw)$
so
    $\U_k 
    = \M^{(n_k+1)}_{\lfs_{n_k+1)}^-}$
for all such $k$, and so also
    $\W_k = \U_k$
for these $k$.

If 
    $\rkkw=\omega$, 
$k_i$ is defined for each $i\in\omega$, $\W_k$ is defined for 
each $k\in\omega$,
    $\ens{\W_k:k\in\omega} = \ens{\U_k:k\in\omega}$,
and 
    $\rllw = \rkkw = \omega$.

\parSkip(b.2)
Finally, still assuming that 
    $\tp_m = 1$
for some least $m$, suppose that the second case of tail relocation
lemma applies, i.e., for some unique $k = k^0\in\rkkw$ (ii) holds. 

If $k^0=0$ then
    $\U_0 = \M^{(1)}_1$,
while for all $m$,
    $\U_1 = \M^{(m+1)}_0$,
and
    $\U_{m+2} 
    = \M^{(m+1)}_{\lfs_{m+1}^-}$,
so
    $\W_0 = \U_0$
and for all $k\ge 1$,
    $\W_k = \U_{k+1}$,
while $\U_1$ never appears in the tail array.
In particular,
    $\rllw = \rkkw = \omega$.

So suppose $k_0 > 0$. For all $k < k_0$ case (i) of tail relocation
lemma applies, so the definition starts as under (b.1), i.e., there 
exists a sequence 
     $k_i,i=0,\ldots,i^0$ in $(k,k^0)$ 
such that
    $\W_k = \U_{k+1}$
for
    $k\in\bigcup_i [k_{2i},k_{2i+1})$,
where $k_0:=0$ and
    $\W_{k_{2i+1}} = \U_{k_{2i}}$,
and if $k_{2i^0+1}<k^0$, then
    $\W_k = \U_k$
for $k\in (k_{2i^0+1},k^0)$. The remainder of the definition 
continues as in the case $k^0=0$, i.e.,
    $\U_{k^0+1} = \M^{(k^0+1)}_0 = \M^{(m+1)}_0$
for all $m \ge k^0+1$, while for all $k > k^0$
    $\U_{k} 
    = \M^{(n_k+1)}_{\lfs_{n_k+1}^-}$.
Thus, for all $k > k^0$
    $\W_k = \U_{k+1}$,
while
    $\U_{k^0}$ never gets inserted into $\WW_{\rllw}$,
and also 
    $\rllw = \rkkw$.

\secSkip
\Section{[8.2]}
\Title{Tracking new building blocks in the tail array}

\parSkip\noindent
We close this section by formalising a few observations about
how new building blocks appear in tail arrays.



\parSkip
\Number{(8.5)}
\Title{Definition}
\dispM{
    \eqalign{
    \dsc_0(\kkw_{0}) 
    & := \dsc_0(0) = \news_0 \cr 
    \dsc_{n+1}(\kkw_{n+1}) 
    & := \cases{
    \news_{n+1} 
    & if $\news_n\not=\void$ and hence $\kkw_{n+1} = \kkw_n+1$. \cr
    \dsc_n(\kkw_n)
    & if $\news_n = \void$ and hence $\kkw_{n+1} = \kkw_n$ \cr
    } \cr
    }
}
and for $k<\kkw_{n+1}$ 
\dispM{
    \dsc_{n+1}(k) := \cases{
    \dsc_n(k)
    & if $\tp_n = \omega$ \cr
    \dsc_n(k)\cup
    \ens{i\in\Lfs_{n+1}^-:(\forall j\in\dsc(k))[\wts_n(i)=\wts_n(j)]}
    & if $\tp_n = 0$ \cr
    \ens{i\in\Lfs_{n+1}^+:(\forall j\in\dsc(k))[\wts_n(i)=\wts_n(j)]}
    \cup\ens{\lf(s_n)+i:i\in\dsc_n(k)}
    & if $\tp_n = 1$. \cr
    }
}

\parSkip
\Title{Remark}
    $\M^{(n)}_i \isom \M^{(n')}_{i'}$ 
whenever for some $k\le\kkw_{\max\ens{n,n'}}$
    $i,i' \in\dsc_{\max\ens{n,n'}}(k)$
and $k$ is minimal for the property that
    $(\exists n'')[i,i'\in\dsc_{n''}(k)]$.


\parSkip
\Number{(8.6)}
\Title{Definition}
\dispM{
    \chi(k) := 
    \inf\ens{n\in\omega:
    \dsc_n(k)\cap\Lfs_{n+1}^-\not=\emptyset}.
}
\par\noindent
That is, $\chi(k)$ is the least level $n$ at which
the element $\U_k\in\UU$ gets 'copied' as $\M^{(n)}_i$ 
in a position $i\in\Lfs_{n+1}^-$ of the 'tail-array'
    $\MM_{n+1}^-$,
if there is, else $\chi(k) = \omega$. 

\parSkip
\Number{(8.7)}
\Title{Lemma}
    $\chi(k) < \omega$
iff for some $n$ such that $\kkw_n\ge k$,
    $\tp_n\in (0,\omega)$.

\Proof
Fix $\U_k$ and let $n_k=\nnw_k$ be the least $n$ 
such that $\kkw_{n+1}=k$.
Then either $k=n=0$, or 
    $\tp_{n_k}\in\ens{0,1}$ 
and
    $\lfs_{n_k+1} = \lfs_{n_k+1}^++\lfs_{n_k+1}^-$.
Let $i_k$ such that 
    $\U_k = \M^{(n_k)}_{i_k}$.
If 
    $\tp_{n_k} = 0$,
then 
    $i_k = \lfs_{n_k} = \lfs_{n_k+1}^+$,
and 
    $\chi(k) = n_k+1$. 
So suppose 
    $\tp_{n_k} = 1$ 
and 
    $i_k = 0$.
If for all $m\ge n_k+1$, 
    $\tp_m = 0$,
or 
    $\tp_m = \omega$, 
then for no $m$ with $\kkw_{m}\ge k$ is $\ps_m > 0$. 
and $\chi(k) = \omega$, since for all such $m$,
    $\MM_{m+1}^-
    = \MM_m$,
i.e., for all $m\ge n_k+1$,
    $\U_k = \M^{m}_0 = \M^{m+1}_0$.
Otherwise, there is a least $m\ge n_k+1$ such 
that $\ps_m\in (0,\omega)$, and for this $m$, 
$\kkw_{m}\ge k$ and
    $\MM_{m+1}^-
    = \MM_m$
so in particular
    $\U_k = \M^{m}_0 = \M^{m+1}_{\lfs_{m+1}^+}$.
\Qed




  
\def\cse{\gamma}
\def\thr{\theta}
\def\Cse{\Gamma}

\chapSkip
\Chapter{[9] Extensions}

\secSkip\noindent
We conclude the paper by pushing the construction 
principle exhibitd in section 7 to 'higher ranks'.

\secSkip
\Section{[9.1]}
\Title{Rank $\omega\cdot\omega=\omega^2$}

\parSkip\noindent
Climbing through levels $\omega\cdot n$ is now fairly 
straight-forward:
Write $\cse$ for $\tup{\ss,\ps}\in\CSE(\omega)$. By 5.12, 
every $A\in\UC(\omega)$ has a representation of form 
$\A(\cse)$, for some $\cse\in\CSE(\omega)$; by 7.7 combined 
with 5.12, every $A\in\UC(\omega+\omega)$ has a representation 
of form $\A(\cse)\rarr{\A(\cse_k)}_{k\in\rkkw}$ for $\rkkw
=\rkkw(\cse)\subeq\omega$ and some $\cse,\cse_k\in\CSE(\omega)$
(only depending on $A$, and therefore satsfying 
   $\f\ss_k=\f\ss_{k'}$, 
but otherwise arbitrary). 
By the same token, every $A\in\UC(\omega+\omega+\omega)$ 
is of form 
$$
\A(\cse)\rarr{
   \A(\cse_k)\rarr{
      \A(\cse_{kl})
   }_{l\in\rkkw(\cse_k)}
}_{k\in\rkkw(\cse)},
\eqno{(\hbox{a})}
$$ 
for $\cse,\cse_k,\cse_{kl}$ (again, only depending on $A$, and 
therefore satsfying 
   $\f\ss_k= \f\ss_{k'}$
and
   $\f\ss_{kl}= \f\ss_{kl'}$
but otherwise arbitrary). Let's abbreviate (a) as
$$
\A(\cse;\cse_k;\cse_{kl})
=
\A(\cse;\cse_k;\cse_{kl})_{k\in\rkkw(\cse),l\in\rkkw(\cse_k)}
$$
Formally, by induction,

\parSkip
\Number{(9.1)}
\Title{Definition}
$$
\eqalign{
\A(\cse;
   \cse_{k_1};
   \cse_{k_1k_2};
   \ldots;
   \cse_{k_1\ldots k_n})
&:=
\A(\cse;
   \cse_{k_1};
   \cse_{k_1k_2};
   \ldots;
   \cse_{k_1\ldots k_n})_{
   k_1\in\rkkw(\cse),
   k_2\in\rkkw(\cse_{k_1}),
   \ldots,
   k_n\in\rkkw(\cse_{k_1\ldots k_{n-1}})}
\cr
&:=
\A(\cse)\rarr{
   \A(\cse_{k_1};
      \cse_{k_1k_2};
      \ldots;
      \cse_{k_1\ldots k_n})_{
      k_2\in\rkkw(\cse_{k_1}),
      \ldots,
      k_n\in\rkkw(\cse_{k_1l\ldots k_{n-1}})}}_{
   k_1\in\rkkw(\cse)}
\cr
}
$$
 
Generalizing the preceding observations: 

\parSkip
\Number{(9.2)}
\Title{Corollary}
Every $A\in\UC(\omega\cdot n)$ is of form
$$
\A(\cse;
   \cse_{k_1};
   \cse_{k_1k_2};
   \ldots;
   \cse_{k_1\ldots k_n})
=
\A(\cse;
   \cse_{k_1};
   \cse_{k_1k_2};
   \ldots;
   \cse_{k_1\ldots k_n})_{
   k_1\in\rkkw(\cse),
   k_2\in\rkkw(\cse_{k_1}),
   \ldots,
   k_n\in\rkkw(\cse_{k_1\ldots k_{n-1}})}
$$ 
for some $\cse,\cse_{k_1},\ldots,\cse_{k_1\ldots k_n}$ 
(only depending on $A$, and therefore satisfying
   $\f\ss_{k_1}= \f\ss_{k_1'}$
   $\f\ss_{k_1k_2}= \f\ss_{k_1k_2'}$
   $\ldots$
   $\f\ss_{k_1\ldots k_{n-1}k_n}= \f\ss_{k_1\ldots k_{n-1}k_n'}$
   $\ldots$,
but otherwise arbitrary). 
\qed

\parSkip
Now let 
   $x\in A\in\UC(\omega^2)$; 
For each $n$,
$$
\cmp{x}{\omega\cdot n}{A}
=\A(\cse^{(n)};
    \cse^{(n)}_{k_1};
    \cse^{(n)}_{k_1k_2};
    \ldots;
    \cse^{(n)}_{k_1\ldots k_n})_{
    k_1\in\rkkw(\cse^{(n)}),
    k_2\in\rkkw(\cse^{(n)}_{k_1}),
    \ldots,
    k_n\in\rkkw(\cse^{(n)}_{k_1\ldots k_{n-1}})}
$$
for some
$
 \cse^{(n)},
 \cse^{(n)}_{k_1},
 \cse_{k_1k_2},
 \ldots,
 \cse^{(n)}_{k_1\ldots k_n}\in\CSE(\omega)
$, 
such that for each $n$ there is some 
   $\thr(n):=h\in\rkkw(\cse^{n+1})$, 
such that
$$
\eqalign{
&
\A(\cse^{(n)};
   \cse^{(n)}_{k_1};
   \cse^{(n)}_{k_1k_2};
   \ldots;
   \cse^{(n)}_{k_1\ldots k_n})_{
   k_1\in\rkkw(\cse^{(n)}),
   k_2\in\rkkw(\cse^{(n)}_{k_1}),
   \ldots,
   k_n\in\rkkw(\cse^{(n)}_{k_1\ldots k_{n-1}})}
\cr&\qquad
=
\A(\cse^{(n+1)}_{h};
   \cse^{(n+1)}_{hk_2};
   \ldots;
   \cse^{(n+1)}_{h,k_2\ldots k_n})_{
   k_2\in\rkkw(\cse^{(n+1)}_{h}),
   \ldots,
   k_n\in\rkkw(\cse^{(n+1)}_{h,k_2\ldots k_{n-1}})}
\cr
}
\eqno{(9.1)}
$$
and hence such that
$$
\eqalign{
&\cse^{(n)}=\cse^{(n+1)}_{h}\cr
&\cse^{(n)}_{k_1}=\cse^{(n+1)}_{hk_1}\cr
&\ldots\cr
&\cse^{(n)}_{k_1\ldots k_n}=\cse^{(n+1)}_{hk_1\ldots k_n}\cr
}
\eqno{(9.2)}
$$

\parSkip
Define an 'embedding thread' between the $\omega$-limit levels:

\parSkip
\Number{(9.3)}
\Title{Definition}
Let $\thr$ be the function defined via
$\thr(n)=h$ just if (9.1), and hence (9.2), holds.
Let 
$$A(\cse^{(0)};
    \cse^{(1)};
    \cse^{(1)}_{k^{(1)}_1};
    \cse^{(2)};
    \cse^{(2)}_{k^{(2)}_1};
    \cse^{(2)}_{k^{(2)}_1k^{(2)}_2}
    \cse^{(3)};
    \cse^{(3)}_{k^{(3)}_1};
    \cse^{(3)}_{k^{(3)}_1k^{(3)}_2};
    \cse^{(3)}_{k^{(3)}_1k^{(3)}_2k^{(3)}_3};
\ldots)
\eqno{(9.3)}
$$
be the unique 
   $A\in\UC(\omega^2)$ 
such that
$$
\eqalign{
\exists h_0&\;\cse^{(0)}=\cse^{(1)}_{h_0};
\cr
\exists h_1&\;\cse^{(1)}=\cse^{(2)}_{h_1}
\cr
   \&\; &
   \tup{\cse^{(1)}_{k^{(1)}_1} :
      k^{(1)}_1\in\rkkw(\cse^{(1)})
   }
  =\tup{\cse^{(2)}_{h_1k^{(2)}_2} :
      k^{(2)}_2\in\rkkw(\cse^{(2)}_{h_1})
   };
\cr
\exists h_2&\;\cse^{(2)}=\cse^{(3)}_{h_2}
\cr
   \&\; &
   \tup{\cse^{(2)}_{k^{(2)}_1} :
      k^{(2)}_1\in\rkkw(\cse^{(2)})
   }
  =\tup{\cse^{(3)}_{h_2k^{(3)}_2} :
      k^{(3)}_2\in\rkkw(\cse^{(3)}_{h_2})
   }
\cr
   \&\; &
   \tup{\cse^{(2)}_{k^{(2)}_1k^{(2)}_2} :
      k^{(2)}_2\in\rkkw(\cse^{(2)}_{k^{(2)}_1})
   }
  =\tup{\cse^{(3)}_{h_2k^{(3)}_2k^{(3)}_3} :
      k^{(3)}_3\in\rkkw(\cse^{(3)}_{h_2k^{(3)}_2})
   };
\cr
&\ldots\cr
}
\eqno{(9.4)}
$$

\parSkip
\Number{(9.4)}
\Title{Corollary}
\Subtitle{Rank $\omega^2$ representation lemma}
Every 
   $A\in\UC(\omega^2)$ 
has a representation of the
form (9.3).
\Qed

For completeness, note that the $n$th block of $\cse$ tokens
in representation (9.3) of $A$ is 
$$
\ldots
\cse^{(n)}_{k^{(n)}_1};
   \cse^{(n)}_{k^{(n)}_1k^{(n)}_2};
   \cse^{(n)}_{k^{(n)}_1k^{(n)}_2k^{(n)}_3};
   \ldots;
   \cse^{(n)}_{k^{(n)}_1\ldots k^{(n)}_n};
\ldots
\eqno{(9.3a)}
$$
and the 'embedding conditions' on the $n$th and $n+1$st
block are:
$$
\eqalign{
&\ldots\cr
\exists h_n&\;\cse^{(n)}=\cse^{(n+1)}_{h_n}
\cr
   \&\; &
   \tup{\cse^{(n)}_{k^{(n)}_1} :
      k^{(n)}_1\in\rkkw(\cse^{(n)})
   }
  =\tup{\cse^{(n+1)}_{h_nk^{(n+1)}_2} :
      k^{(n+1)}_2\in\rkkw(\cse^{(n+1)}_{h_n})
   }
\cr
   \&\; &
   \ldots
\cr
   \&\; &
   \tup{\cse^{(n)}_{k^{(n)}_1\ldots k^{(n)}_{i+1}} :
      k^{(n)}_{i+1}\in\rkkw(\cse^{(n)}_{k^{(n)}_1\ldots k^{(n)}_i})
   }
  =\tup{\cse^{(n+1)}_{h_nk^{(n+1)}_2\ldots k^{(n+1)}_{i+1}} :
      k^{(n+1)}_{i+1}\in\rkkw(\cse^{(n+1)}_{h_nk^{(n)}_2
                                            \ldots k^{(n+1)}_i})
   };
\cr
   \&\; &
   \ldots
\cr
   \&\; &
   \tup{\cse^{(n)}_{k^{(n)}_1\ldots k^{(n)}_n} :
      k^{(n)}_n\in\rkkw(\cse^{(n)}_{k^{(n)}_1\ldots k^{(n)}_{n-1}})
   }
  =\tup{\cse^{(n+1)}_{h_nk^{(n+1)}_2\ldots k^{(n+1)}_n} :
      k^{(n+1)}_n\in\rkkw(\cse^{(n+1)}_{h_nk^{(n)}_2
                                           \ldots k^{(n+1)}_{n-1}})
   }.
\cr
&\ldots\cr
}
\eqno{(9.4a)}
$$
The 'implicit' conditions on the
   $\cse^{(n)}_{i_1\ldots i_n}
   =\tup{\ss^{(n)}_{i_1\ldots i_n},\ps^{(n)}_{i_1\ldots i_n}}$
are:
   $$\f\ss^{(n)}_{k_1\ldots k_{n-1}k_n}
   = \f\ss^{(n)}_{k_1\ldots k_{n-1}k_n'}.$$
Let's call a system of 
   $\cse^{(n)}_{k^{(n)}_1\ldots k^{(n)}_n}\sub\CSE(\omega)$
satisfying the 'embedding conditions' (9.4) (as well as the
'implicit' conditions) a 
   {\it $\CSE(\omega^2)$ system}
and denote such by $\Cse$. Corollary 9.4 says 
that every 
   $A\in\UC(\omega^2)$ 
is represented by a
   $\CSE(\omega^2)$ 
system. The following notation denotes a 
   $\CSE(\omega^2)$ 
system, the arrows expressing the 'embedding conditions' (9.4):
$$
\matrix{
\omega&:&
\cse^{(0)}&\cr
&&
          &\searrow\cr
\omega\cdot 2&:&
\cse^{(1)}&&\cse^{(1)}_{k^{(1)}_1}\cr
&&
          &\searrow&&\searrow\cr
\omega\cdot 3&:&
\cse^{(2)}&&\cse^{(2)}_{k^{(2)}_1}&&\cse^{(2)}_{k^{(2)}_1k^{(2)}_2}\cr
&&
          &\searrow&&\searrow&&\searrow\cr
&&
\ldots    &\ldots&\ldots&\ldots&\ldots&\ldots\cr
&&
          &\searrow&&&&\searrow\cr
\omega\cdot n&:&
\cse^{(n)}&&\cse^{(n)}_{k^{(n)}_1}&&\ldots
          &&\cse^{(n)}_{k^{(n)}_1\ldots k^{(n)}_n}\cr
&&
          &\searrow&&\searrow&&&&\searrow\cr
&&
\ldots    &\ldots&\ldots&\ldots&\ldots&\ldots&\ldots&\ldots\cr
}
\eqno{(9.5)}
$$

\secSkip
\Section{[9.2]}
\Title{Rank $\omega^3$}

\par\noindent
Beyond $\omega^2$ the pace begins to look slow, since
we still only progress by $\omega$ steps, using the
construction principle of 7.7. Schematically:
$$
\matrix{
\omega^2&:&
\Cse^{(0)}&\cr
&&
          &\searrow\cr
\omega^2+\omega&:&
\cse^{(1)}&&\Cse^{(1)}_{k^{(1)}_1}\cr
&&
          &\searrow&&\searrow\cr
\omega^2+\omega\cdot 2&:&
\cse^{(2)}&&\cse^{(2)}_{k^{(2)}_1}&&\Cse^{(2)}_{k^{(2)}_1k^{(2)}_2}\cr
&&
          &\searrow&&\searrow&&\searrow\cr
&&
\ldots    &\ldots&\ldots&\ldots&\ldots&\ldots\cr
&&
          &\searrow&&&&\searrow\cr
\omega^2+\omega\cdot n&:&
\cse^{(n)}&&\cse^{(n)}_{k^{(n)}_1}&&\ldots
          &&\Cse^{(n)}_{k^{(n)}_1\ldots k^{(n)}_n}\cr
&&
          &\searrow&&\searrow&&&&\searrow\cr
&&
\ldots    &\ldots&\ldots&\ldots&\ldots&\ldots&\ldots&\ldots\cr
}
\eqno{(9.6)}
$$
where the embedding conditions for the 
   $\cse^{(n)}_{i_1\ldots i_n}$
translate verbatim into embedding conditions for the
   $\Cse^{(n)}_{i_1\ldots i_n}$.
(The 'implicit' conditions on the
   $\cse^{(n)}_{i_1\ldots i_{n-1}}$
are as before; the 'implicit' conditions on the
   $\Cse^{(n)}_{i_1\ldots i_n}$
are 'implicit' in the 'implicit' conditions of their
constituent $\cse^{(n)}_{i_1\ldots i_n}$.)
Again, given $A\in\UC(\omega^2+\omega^2)$ and $x\in A$, 
there is a system as in (9.6) such that for each $n$,
$$
\cmp{x}{\omega^2+\omega\cdot n}{A}
=\A(\cse^{(n)};
    \cse^{(n)}_{k_1};
    \cse^{(n)}_{k_1k_2};
    \ldots;
    \cse^{(n)}_{k_1\ldots k_{n-1}};
    \Cse^{(n)}_{k_1\ldots k_n})_{
    k_1\in\rkkw(\cse^{(n)}),
    k_2\in\rkkw(\cse^{(n)}_{k_1}),
    \ldots,
    k_n\in\rkkw(\cse^{(n)}_{k_1\ldots k_{n-1}})},
$$
so every $A\in\UC(\omega^2+\omega^2)$ has a representation
by a $\CSE(\omega^2\cdot 2)$ system of from (9.6). Now iterate. 
Replace the $\CSE(\omega^2)$ systems $\Cse^{(n)}_{i_1\ldots i_n}$ 
in (9.6) by $\CSE(\omega^2\cdot 2)$ systems (again denoted 
$\Cse^{(n)}_{ i_1\ldots i_n}$) and obtain a $\CSE(\omega^2\cdot 3)$ 
system representing $A\in\UC(\omega^2\cdot 3)$.

Formally, if (for some ordinal $\alpha$) the concept of
a $\CSE(\alpha)$ system is defined, a $\CSE(\alpha+\omega^2)$
system is a system for form
$$
\matrix{
\alpha&:&
\Cse^{(0)}&\cr
&&
          &\searrow\cr
\alpha+\omega&:&
\cse^{(1)}&&\Cse^{(1)}_{k^{(1)}_1}\cr
&&
          &\searrow&&\searrow\cr
\alpha+\omega\cdot 2&:&
\cse^{(2)}&&\cse^{(2)}_{k^{(2)}_1}&&\Cse^{(2)}_{k^{(2)}_1k^{(2)}_2}\cr
&&
          &\searrow&&\searrow&&\searrow\cr
&&
\ldots    &\ldots&\ldots&\ldots&\ldots&\ldots\cr
&&
          &\searrow&&&&\searrow\cr
\alpha+\omega\cdot n&:&
\cse^{(n)}&&\cse^{(n)}_{k^{(n)}_1}&&\ldots
          &&\Cse^{(n)}_{k^{(n)}_1\ldots k^{(n)}_n}\cr
&&
          &\searrow&&\searrow&&&&\searrow\cr
&&
\ldots    &\ldots&\ldots&\ldots&\ldots&\ldots&\ldots&\ldots\cr
}
$$
where each $\Cse^{(n)}_{k^{(n)}_1\ldots k^{(n)}_n}$ is a
$\CSE(\alpha)$ system, and the arrows represent the embedding
conditions (9.4):
$$
\eqalign{
\exists h_0&\;\Cse^{(0)}=\Cse^{(1)}_{h_0};
\cr
\exists h_1&\;\cse^{(1)}=\cse^{(2)}_{h_1}
\cr
   \&\; &
   \tup{\Cse^{(1)}_{k^{(1)}_1} :
      k^{(1)}_1\in\rkkw(\cse^{(1)})
   }
  =\tup{\Cse^{(2)}_{h_1k^{(2)}_2} :
      k^{(2)}_2\in\rkkw(\cse^{(2)}_{h_1})
   };
\cr
\exists h_2&\;\cse^{(2)}=\cse^{(3)}_{h_2}
\cr
   \&\; &
   \tup{\cse^{(2)}_{k^{(2)}_1} :
      k^{(2)}_1\in\rkkw(\cse^{(2)})
   }
  =\tup{\cse^{(3)}_{h_2k^{(3)}_2} :
      k^{(3)}_2\in\rkkw(\cse^{(3)}_{h_2})
   }
\cr
   \&\; &
   \tup{\Cse^{(2)}_{k^{(2)}_1k^{(2)}_2} :
      k^{(2)}_2\in\rkkw(\cse^{(2)}_{k^{(2)}_1})
   }
  =\tup{\Cse^{(3)}_{h_2k^{(3)}_2k^{(3)}_3} :
      k^{(3)}_3\in\rkkw(\cse^{(3)}_{h_2k^{(3)}_2})
   };
\cr
&\ldots\cr
\exists h_n&\;\cse^{(n)}=\cse^{(n+1)}_{h_n}
\cr
   \&\; &
   \tup{\cse^{(n)}_{k^{(n)}_1} :
      k^{(n)}_1\in\rkkw(\cse^{(n)})
   }
  =\tup{\cse^{(n+1)}_{h_nk^{(n+1)}_2} :
      k^{(n+1)}_2\in\rkkw(\cse^{(n+1)}_{h_n})
   }
\cr
   \&\; &
   \ldots
\cr
   \&\; &
   \tup{\cse^{(n)}_{k^{(n)}_1\ldots k^{(n)}_{i+1}} :
      k^{(n)}_{i+1}\in\rkkw(\cse^{(n)}_{k^{(n)}_1\ldots k^{(n)}_i})
   }
  =\tup{\cse^{(n+1)}_{h_nk^{(n+1)}_2\ldots k^{(n+1)}_{i+1}} :
      k^{(n+1)}_{i+1}\in\rkkw(\cse^{(n+1)}_{h_nk^{(n)}_2
                                            \ldots k^{(n+1)}_i})
   };
\cr
   \&\; &
   \ldots
\cr
   \&\; &
   \tup{\Cse^{(n)}_{k^{(n)}_1\ldots k^{(n)}_n} :
      k^{(n)}_n\in\rkkw(\cse^{(n)}_{k^{(n)}_1\ldots k^{(n)}_{n-1}})
   }
  =\tup{\Cse^{(n+1)}_{h_nk^{(n+1)}_2\ldots k^{(n+1)}_n} :
      k^{(n+1)}_n\in\rkkw(\cse^{(n+1)}_{h_nk^{(n)}_2
                                           \ldots k^{(n+1)}_{n-1}})
   };
\cr
&\ldots\cr
}
\eqno{(9.4b)}
$$

Iterating systems
   $\CSE(\omega^2\cdot n+\omega^2)$ 
through all $n$, we obtain a 
   $\CSE(\omega^2\cdot\omega)=\CSE(\omega^3)$ 
system (comprising within it $\CSE(\omega^2\cdot n)$ systems 
as 'building-blocks'), and every $A\in\UC(\omega^3)$ is 
represented by one such by the inductive reasoning of before, 
i.e., since each component $\cmp{x}{\omega^2\cdot n}{A}$
is represented by a constituent $\CSE(\omega^2\cdot n)$ system.

\secSkip
\Section{[9.3]}
\Title{Rank $\omega^{\omega}$}

\par\noindent
Now, we'd like to define a 
   $\CSE(\alpha+\omega^3)$ 
as the iteration just described, starting with a $\CSE(\alpha)$
system $\Cse^{(0)}_0$, Iterating yet again such systems
   $\CSE(\omega^3\cdot n+\omega^3)$
through all $n$, we should obtain a
   $\CSE(\omega^4)$
system. Iterating through all 
   $\CSE(\omega^n)$
systems, we should obtain a 
   $\CSE(\omega^{\omega})$
system. {\it And so on.}
But we are not ready to carry this out here, and will contend 
ourselves with merely having indicated how this 'iteration' 
might be continued through an initial segment of the ordinals.

\parSkip
\Number{(9.5)}
\Title{Challenge}
Devise a recursive system of notation for uniform chains of
constructible rank.





  \vfill\eject
  \centerline{\bf References:}
  \hfill\break [P] Francis Pastijn, {\sl A class of uniform chains},
  in: S.M.Goberstein and P.M.Higgins [eds.] Semigroups
  and their applications.
  \vfill\eject

  \vfill\eject
  \centerline{\bf Revision history:}

  \par\noindent
  Research for this paper was carried out between 1996 and 1999 by
  the first author. A first draft was completed December 2000. 
  A substantially revised second draft, begun almost a year later, 
  was completed January 2003. Some streamlining and improvements to 
  the diagrammes resulted in the first submission to the 'arxiv' 
  (December 2003). The ideas of the conclusion added for this final 
  revision represent unfinished business.
  [Contact address: 94 Spencer Place, Leeds, LS7 4DT, UK]
  

  \end\bye


\end\bye